\documentclass[11pt]{article}
\usepackage[utf8]{inputenc}
\usepackage{mathrsfs}
\usepackage{amsmath}
\usepackage{amssymb}
\allowdisplaybreaks
\usepackage[unicode=true,bookmarks=false,breaklinks=false,pdfborder={0 0 1},backref=section,colorlinks=true]{hyperref}
\hypersetup{citecolor=blue,urlcolor=blue,linkcolor=blue}

\makeatletter

\usepackage{graphicx}
\usepackage{amsfonts}
\usepackage{epstopdf}
\usepackage{cite}
\usepackage{todonotes}
\usepackage{graphicx}
\usepackage{subcaption}
\usepackage{caption}
\usepackage{xcolor}
\usepackage{fullpage}
\usepackage[misc]{ifsym}
\usepackage{algorithm}
\usepackage[noend]{algpseudocode}
\usepackage{mathrsfs}
\usepackage{fontawesome5}
\usepackage{booktabs}
\makeatletter
\newenvironment{breakablealgorithm}{%
  \begin{center}
    \refstepcounter{algorithm}%
    \renewcommand{\caption}[1]{%
      \addcontentsline{loa}{algorithm}{\protect\numberline{\thealgorithm}##1}%
      \parbox{\textwidth}{%
        \hrule height.8pt depth0pt \kern2pt
        {\raggedright\textbf{\fname@algorithm~\thealgorithm} ##1\par}%
        \kern2pt\hrule\kern2pt}}%
}{%
  \kern2pt\hrule\relax
  \end{center}}
\makeatother
\graphicspath{{figures/}{./}}

\providecommand{\U}[1]{\protect\rule{.1in}{.1in}}

\newtheorem{theorem}{Theorem}[section]

\newtheorem{lemma}{Lemma}[section]
\newtheorem{definition}{Definition}[section]

\newtheorem{corollary}{Corollary}[section]
\newtheorem{remark}{Remark}[section]
\newtheorem{example}{Example}[section]
\newenvironment{proof}[1][Proof]{\textbf{#1.} }{\ \rule{0.5em}{0.5em}}
\newtheorem{proposition}{Proposition}[section]
\newtheorem{assumption}{Assumption}

\makeatother

\begin{document}
\title{A Unified Bregman Projection Framework for Solving Multi-Split Multiple-Set Implicit Convex Feasibility Problems}
\author{Soumitra Dey$^{1,*}$, Aviv Gibali$^{2}$, Simeon Reich$^{3}$\medskip{}\\
$^1$ Einstein Institute of Mathematics, Hebrew University of Jerusalem, \\Jerusalem 9190401, Israel\\
$^2$ Department of Applied Mathematics, Holon Institute of Technology,\\ Holon 5810201, Israel\\
$^3$ Department of Mathematics, The Technion -- Israel Institute of Technology, \\Haifa 3200003, Israel\\
% \faEnvelope\, deysoumitra2012@gmail.com (S. Dey)\\
% \faEnvelope\, avivgi@hit.ac.il (A. Gibali)\\
% \faEnvelope\, sreich@technion.ac.il (S. Reich)\\
% $^*$ Corresponding author\\
}
\date{\today}
\maketitle

\begingroup
\renewcommand{\thefootnote}{\fnsymbol{footnote}}
\footnotetext[0]{%
\begin{minipage}{0.95\textwidth}
\small
\noindent
Soumitra Dey\\
\Letter\ deysoumitra2012@gmail.com
\medskip{}

\noindent Aviv Gibali\\
\Letter\ avivgi@hit.ac.il
\medskip{}
%\faEnvelope

\noindent
Simeon Reich\\
\Letter\ sreich@technion.ac.il

\medskip{}
\noindent
$^{*}$ Corresponding author
\end{minipage}
}
\endgroup

\begin{abstract}
We study the multi-split multiple-set implicit convex feasibility problem (MSMS-ICFP), in which the constraint sets are generated by point-to-set operators and the feasibility conditions are imposed both in the original space and in several image spaces connected by bounded linear operators.
We analyze the generalized proximity function, built from weighted Bregman distances and Bregman projections, introduced by Censor and Dey. % \cite{SDYCEN2026}.
For translation-type point-to-set operators, we establish the continuity and Lipschitz continuity of the corresponding Bregman projections, the continuous differentiability of the proximity function with an explicit formula for its gradient, and global and local Lipschitz estimates for this gradient.
We then propose a generalized simultaneous projection gradient algorithm. In the nonconvex case, the gradient norms of the iterates converge to zero and every accumulation point is stationary. When the proximity function is convex and the solution set is nonempty, the iterates converge in norm to a solution, with an $O(1/k)$ rate for the proximity values, an $O(1/\sqrt{k})$ rate for the gradient norms, and a linear rate under a quadratic growth condition.
Cyclic and block-iterative variants are also analyzed. Numerical experiments, including a cooperative sensor network localization problem, illustrate the theory and show that, in contrast with frozen-set implicit projection methods, the proposed gradient schemes remain convergent when the translation operators are expansive.
\end{abstract}

\medskip{}


\section{Introduction}
\label{Sec:1:Introduction}

The \textit{split inverse problem} (SIP) deals with two different inverse problems in two different spaces. Indeed, let $X$ and $Y$ be two vector spaces, and $A:X\rightarrow Y$ be
a bounded linear operator.
% It was first introduced in \cite[Section 2]{YCEN22012}. Below we recall the formulation of a split inverse problem.
Let $\text{IP}_{1}$ be a first inverse problem, formulated in $X$ and the second
one is $\text{IP}_{2}$, formulated in $Y$. With these given data,
the \textit{split inverse problem} (SIP) is defined as follows:
\[
\begin{array}{r@{\;}l}
\text{Find a point } & x^{*}\in X\text{ that solves}\text{ IP}_{1}\\
 & \text{such that}\\
\text{the point }y^{*} & =Ax^{*}\in Y\text{ solves}\text{ IP}_{2}.\tag{SIP}
\end{array}
\]

It was first introduced in \cite[Section 2]{YCEN22012}. 
Different choices of $\text{IP}_{1}$ and $\text{IP}_{2}$ lead to a large class of split problems, for instance, the split
variational inequality problem (SVIP) \cite{YCEN22012}, and the split multivalued variational inequality problem \cite{SDEY22024}. We refer to \cite{SDEY2024,GLOP2012,AMOU2011,AMOU2014,YSHE2015,SREI22022,CBYR2002} and references therein, for an extensive study of different types of split problems.
Various types of SIPs have recently attracted considerable attention due to their significant applications in, e.g., the intensity-modulated radiation therapy treatment problem \cite{YCEN2005} and in economic theory \cite{JLI2019}.
Motivated by this, in the present paper we consider an important SIP, namely, \textit{multi-split multiple-set implicit convex feasibility problem} and give an application to \textit{cooperative sensor network localization}.
%\eqref{MSMSICFP} considered in this paper is also an important instance of a SIP.

\subsection{Previous work and motivation}
\label{subsec:1:General-Background-and-Motivation}

The {\em convex feasibility problem} (CFP) is a fundamental framework of many inverse problems in various fields of mathematics and physical sciences (see, e.g., \cite{EREM1969,YCEN1982,DEUT1992} and references therein).
Owing to its applications in fields such as medical imaging (e.g., computerized tomography), radiation therapy treatment planning, and electron microscopy, the CFP has attracted the attention of many researchers and has been generalized in various directions.
% Due to its applicability in various fields, over the past decades, a wide group of researchers have focused their attention to it and generalized it in the various directions.
It consists of finding a common point of a finite or (countably) infinite collection of nonempty, closed, and convex subsets of a normed linear space.
Let $\{C_{i}\}_{i=1}^{r}$ be a finite collection of nonempty, closed, and convex subsets of $\mathbb{R}^{n}$.
The CFP is formulated below.
\begin{align}
\mathrm{\mathrm{Find}\:a\:point\:}x^{*}\in\mathbb{R}^{n}\mathrm{\:such\:that\:}x^{*}\in\bigcap_{i=1}^{r}C_{i}.\tag{CFP}\label{CFP:Prob}
\end{align}

Fundamental references for the problem \eqref{CFP:Prob} are \cite{LMBR1965, LMBR1967}.
Thereafter, many researchers extensively explored the problem \eqref{CFP:Prob} in various normed spaces and proposed algorithms to solve it (see, for instance, \cite{BYRN2004,YCEN1982,HHBA1996,BUTN1990,GUBI1967,HANS1988, YCEN1988,RAHA1989,ARDP1992} and references therein).
Later, the problem \eqref{CFP:Prob} was extended in \cite{YCEN2016} to point-to-set operators, under the name {\em implicit convex feasibility problem} (ICFP), and two types of projection algorithms were studied there. Additionally, the authors demonstrated its application in image denoising within a  medical imaging situation.

Projection algorithms open a new avenue for solving a large class of problems in the area of operations research. They employ projections onto nonempty, closed, and convex sets in numerous ways. %They are capable of solving a wide class of different problems in the field of operations research.
Specifically, simultaneous algorithms allow parallel computations, while sequential algorithms process one step at a time and move to the next step after completing the computations at the current step. We refer the reader to \cite{YCEN1997,ACEG2012} for more details on projection algorithms.
Recent progress on \eqref{CFP:Prob} in nonconvex and inconsistent settings, based on cyclic relaxed Douglas--Rachford splitting, can be found in \cite{TLDI2026}. Furthermore, \cite{XWAN2026} studied \eqref{CFP:Prob} in Hilbert spaces, where each convex set is a sublevel set of a real-valued convex function, and established the convergence rates for a family of subgradient algorithms. The performance of the algorithms was illustrated numerically. In this direction \cite{SDEYSR2026} could open a avenue for future research. 

The notion of a proximity (a.k.a. gap) function plays a crucial role in guiding the iterative process toward solutions of various problems in nonlinear analysis.
Such functions provide a quantitative measure of constraint violation and therefore reformulate various nonlinear problems as optimization problems (see, for instance, \cite{ACEG2012,FUKU1992,KUBO2009} and references therein).
In \cite{CBYR2001}, the orthogonal projections and the Euclidean distance in a proximity function for the CFP were replaced by Bregman projections and Bregman distances, which led to a simultaneous multiprojection algorithm based on Bregman projections. Additionally, a multidistance successive generalized projection algorithm is studied in \cite{CBYR22001} for solving the CFP using Bregman projections.

Motivated by the above discussions, in this paper, we consider the \textit{multi-split multiple-set implicit convex feasibility problem} introduced in \cite{SDYCEN2026} as a further generalization of the implicit convex feasibility problem. We recall its formulation below.

Let $\mathbb{R}^{n}$ and $\mathbb{R}^{m_{j}}$ for $j=1,2,\dots,t$
be Euclidean spaces. For each $i=1,2,\dots,r$, let $C_{i}:\mathbb{R}^{n}\rightarrow2^{\mathbb{R}^{n}}$
be a point-to-set operator such that for each $x\in\mathbb{R}^{n}$,
the set $C_{i}(x)$ is nonempty, closed, and convex. For each $j=1,2,\dots,t$
and $h=1,2,\dots,\ell_{j}$, let $Q_{h}^{j}:\mathbb{R}^{m_{j}}\rightarrow2^{\mathbb{R}^{m_{j}}}$
be a point-to-set operator such that for each $x\in\mathbb{R}^{m_{j}}$,
the set $Q_{h}^{j}(x)$ is nonempty, closed, and convex. Let $A_{j}:\mathbb{R}^{n}\rightarrow\mathbb{R}^{m_{j}}$
be a bounded linear operator for each $j=1,2,\dots,t$. With these
given data the {\em multi-split multiple-set implicit convex feasibility
problem} (MSMS-ICFP) is defined as follows: find a point $x^{*}\in\mathbb{R}^{n}$
such that
\begin{equation}
x^{*}\in\bigcap\limits_{i=1}^{r}C_{i}(x^{*})\text{ and }A_{j}x^{*}\in\bigcap\limits_{h=1}^{\ell_{j}}Q_{h}^{j}(A_{j}x^{*}),\;\text{for each }j=1,2,\dots,t.\tag{MSMS-ICFP}\label{MSMSICFP}
\end{equation}

In particular, if we set $t=r$, $\mathbb{R}^{m_{j}}=\mathbb{R}^{n}$, $A_{j}:=I$,
the identity operator, $Q_{h}^{j}(x)=C_{j}(x)$, for all $x\in\mathbb{R}^{n}$, $h=1,2,\dots,\ell_{j}$, $j=1,2,\dots,t$, %$h=i=j$, and $t=r$
in the above formulation, then the problem \eqref{MSMSICFP}
reduces to the classical {\em implicit convex feasibility problem}
(ICFP) of \cite{YCEN2016} which consists of finding a point $x^{*}\in\mathbb{R}^{n}$
such that
\begin{equation}
x^{*}\in\bigcap_{i=1}^{r}C_{i}(x^{*}).\tag{ICFP}\label{ICFP}
\end{equation}

This problem originated in \cite{MRGH2013} in the context of the cooperative wireless sensor network positioning problem, and it was later generalized in \cite{YCEN2016}. We also note that the explicit counterpart of \eqref{MSMSICFP}, namely, the problem \eqref{MSMSCFP} studied in Section~\ref{subsec:MSMSCFP}, can be reduced, under an appropriate constraint qualification, to the multiple-set split convex feasibility problem studied in infinite-dimensional Hilbert spaces in \cite{EMAS2007}.

\subsection{Main contributions}
\label{subsec:Our-Contribution}

The main purpose of this paper is to develop the Bregman proximity-function framework introduced in \cite{SDYCEN2026} and the associated gradient-type algorithms for the problem \eqref{MSMSICFP}. The principal contributions are summarized as follows.

\begin{enumerate}
\item We establish an explicit formula for the Bregman projection for translation-type point-to-set operators of the form
\begin{eqnarray}\label{ineq1:contributions}
    C(x)=g(x)+\Omega.
\end{eqnarray}
We compute the Bregman projection for the above mentioned translation-type operator with the general quadratic Bregman function to support this explicit formula.
Additionally, we also demonstrate, by means of an example, that a standard Lipschitz estimate available for orthogonal projections does not extend directly to general Bregman projections.

\item We derive continuity properties of the operators
\begin{eqnarray}\label{ineq2:contributions}
    x\mapsto P_{C(x)}^{f}(x)
    \;\;\text{and}\;\;
    x\mapsto P_{Q(Ax)}^{g}(Ax)
\end{eqnarray}
under suitable assumptions on the Bregman functions and the translation operators. We further provide examples showing that either inner semi-continuity or outer semi-continuity alone is not sufficient to guarantee the continuity of the implicit Bregman projection operator. Moreover, we establish the Lipschitz continuity of the above mentioned operators under appropriate assumptions.

\item Using Danskin's theorem, we establish the continuous differentiability of the generalized proximity function and derive an explicit formula for its gradient in the case of translation-type point-to-set operators. The formula involves the Jacobians of the translation operators, the Hessians of the Bregman functions, the Bregman projections, and the bounded linear operators connecting the underlying spaces.

\item We provide sufficient conditions ensuring the global Lipschitz continuity of the gradient of the generalized proximity function for the translation-type point-to-set operators using the results established for the operators in \eqref{ineq2:contributions} and derive an explicit estimate of its Lipschitz constant.

\item We propose a generalized simultaneous projection gradient algorithm for solving the problem \eqref{MSMSICFP}. When the proximity function is not necessarily convex, we prove that the gradient norms of the generated iterates converge to zero, that every accumulation point is a stationary point, and that $\min_{0\leq k\leq N}\|\nabla p_{f,g}(x^k)\|^2=O(1/N)$. In the convex setting, we establish norm convergence of the entire sequence to a solution of the problem \eqref{MSMSICFP}, the rates $p_{f,g}(x^k)=O(1/k)$ and $\|\nabla p_{f,g}(x^k)\|=O(1/\sqrt{k})$, and a linear (geometric) rate under a quadratic growth condition.

\item We introduce a cyclic projection gradient algorithm with steered step sizes and prove its norm convergence to a solution of \eqref{MSMSICFP} when the components of the proximity function are convex. We give verifiable sufficient conditions for this convexity assumption, which cover the applications considered in this paper, as well as an example with a non-quadratic Bregman function.

\item We establish a block-iterative algorithm for the problem \eqref{MSMSCFP}, which constitutes a particular case of the problem \eqref{MSMSICFP}, and prove its convergence by means of a dominating Bregman distance. Whether Algorithm \ref{Algo:GPBCPA} can be extended to the more general problem \eqref{MSMSICFP} is left as an open question.

\item We also derive weaker sufficient conditions guaranteeing the Lipschitz continuity of the gradient on a suitable residual-bounded subset. This local analysis reduces the need to impose global Lipschitz assumptions on all the component operators and leads to a projected gradient-type algorithm, for which we prove convergence and an $O(1/k)$ rate.

\item We present numerical experiments on an academic example with nonlinear and expansive translation operators and on a cooperative wireless sensor network localization problem, and we compare our algorithms with frozen-set implicit projection methods of the type studied in \cite{YCEN2016,MRGH2013}.
\end{enumerate}

\subsection{The organization of this paper}
This paper is organized as follows. In Section \ref{Sec:2:Preliminaries}, we recall some preliminary results that we have used in our analysis. In Section \ref{Sec:3:AlgorithmsANDBregamnGe}, we establish an explicit form of the Bregman projection onto translation-type sets and characterize its analytical properties under appropriate assumptions. Using these results we study a couple of algorithms to solve the problem \eqref{MSMSICFP} and analyze the convergence of the sequences generated by the algorithms. We provide several examples to illustrate our theoretical results. Additionally, we study a block-iterative algorithm for the problem \eqref{MSMSCFP}, which is a particular case of the problem \eqref{MSMSICFP}, and establish its convergence analysis with appropriate assumptions.
Moreover, we study a projected algorithm under local assumptions.
In Section \ref{Sec:numerics}, we report numerical experiments, and in Section \ref{Sec:conclusion}, we draw several conclusions of our work.

\section{Preliminaries}
\label{Sec:2:Preliminaries}

In this section, we recall some basic definitions and results that
will be used later in the paper. Throughout this paper, $I$ denotes the identity operator on the Euclidean space under consideration, and $\|\cdot\|$ denotes the Euclidean
norm induced by the inner product $\left\langle \cdot,\cdot\right\rangle$; for a matrix, $\|\cdot\|_2$ denotes the spectral norm. An operator $T:\mathbb{R}^n\to\mathbb{R}^n$ is $\nu$-inverse strongly monotone ($\nu$-ism, or $\nu$-cocoercive) if $\left\langle T(x)-T(y), x-y\right\rangle\geq \nu\|T(x)-T(y)\|^2$ for all $x,y\in\mathbb{R}^n$. For a nonempty, closed, and convex set $C\subseteq\mathbb{R}^n$, $\mathrm{dist}(x,C):=\min_{y\in C}\|x-y\|$, and the metric projection $P_C(x):=\arg\min_{y\in C}\|x-y\|$ is well defined and characterized by \cite{HHBA2011}
\begin{equation}
z=P_{C}(x)\Leftrightarrow z\in C\;\text{and}\;\left\langle x-z,y-z\right\rangle \leq0\;\mathrm{for\;all\;}y\in C.\label{proj:inq1}
\end{equation}
A point-to-set operator $C:\mathbb{R}^n\to2^{\mathbb{R}^n}$ is \emph{outer semi-continuous} at $x$ if $x_k\to x$, $y_k\in C(x_k)$ and $y_k\to y$ imply $y\in C(x)$, and \emph{inner semi-continuous} at $x$ if for every $y\in C(x)$ and every sequence $x_k\to x$ there exist $y_k\in C(x_k)$ with $y_k\to y$ (see \cite{RTRO2009}).

\begin{lemma}\cite{SKAR1969}\label{lemma:strongconvex}
    A differentiable function $f:\mathbb{R}^n\rightarrow\mathbb{R}$ is $\sigma$-strongly convex if and only if
        $\left\langle \nabla f(y)-\nabla f(x),y-x\right\rangle\geq\sigma\|y-x\|^2$ for all $x,y\in\mathbb{R}^n$.
\end{lemma}

\begin{theorem}[Danskin's theorem \cite{ASHA2009}]\label{Danskin's:Thm}
Let $C\subseteq\mathbb{R}^n$ be nonempty and compact, and let $\psi:\mathbb{R}^n\times C\rightarrow\mathbb{R}$ be such that $\psi(\cdot, y)$ is differentiable for every $y\in C$ and $\nabla_x \psi$ is continuous on $\mathbb{R}^n\times C$.
If $\psi(x_0,\cdot)$ attains its maximum over $C$ at a unique point $\bar y(x_0)$, then $g(x):=\max_{y\in C}\psi(x,y)$ is differentiable at $x_0$ and $\nabla g(x_0)=\nabla_x \psi(x_0,\bar{y}(x_0))$.
\end{theorem}

\begin{lemma}[Baillon--Haddad \cite{HHBA2011}]\label{Baillon-Haddad}
Let $f:\mathbb{R}^n\rightarrow\mathbb{R}$ be convex and differentiable, and let $\beta>0$.
Then $\nabla f$ is $\beta$-Lipschitz continuous if and only if $\nabla f$ is $1/\beta$-ism.
\end{lemma}

\begin{lemma}[Descent lemma \cite{ABEC2017}]\label{lem:descent}
    If $f:\mathbb{R}^n\to\mathbb{R}$ is differentiable and $\nabla f$ is $L$-Lipschitz continuous on a convex set $S$, then
       $f(x)\leq f(y)+\left\langle\nabla f(y), x-y\right\rangle+\frac{L}{2}\|x-y\|^2$ for all $x,y\in S$.
\end{lemma}

\begin{definition}\cite{YCEN2004}
    A sequence of positive real numbers $\{\sigma_k\}$ is a \emph{steering sequence} if
        $\lim_{k\to\infty}\sigma_k=0$, $\lim_{k\to\infty}\sigma_{k+1}/\sigma_k=1$, and $\sum_{k\geq 0}\sigma_k=+\infty$. Given $q\in\mathbb{N}$, a steering sequence is a \emph{$q$-steering sequence} if $\lim_{k\to\infty}\sigma_{q k+s}/\sigma_{q k}=1$ for all $0\leq s\leq q -1$.
For instance, $\sigma_k=\sigma_0(1+k/q)^{-\gamma}$ with $\sigma_0>0$ and $\gamma\in(0,1]$ is a $q$-steering sequence.
\end{definition}

\section{Generalized algorithms for MSMS-ICFP based on the Bregman distance}
\label{Sec:3:AlgorithmsANDBregamnGe}

We denote $\mathscr{B}(\mathbb{R}^n)$ by the set of all Bregman functions with {\em zone} $\mathbb{R}^n$. Throughout this section, whenever we say $f\in\mathscr{B}(\mathbb{R}^n)$, we mean that $f:\mathbb{R}^n\rightarrow\mathbb{R}$ is a Bregman function. We refer the reader to \cite[Chapter 2]{YCEN1997} for a precise definition of a Bregman function with {\em zone} $\mathbb{R}^n$.

Let $f\in\mathscr{B}(\mathbb{R}^n)$. Then the Bregman
distance is given by
\begin{eqnarray}
\begin{aligned}\label{ineq:bregmain}
D_{f}(y,x):=f(y)-f(x)-\left\langle \nabla f(x),y-x\right\rangle .
\end{aligned}
\end{eqnarray}

Consider a nonempty, closed, and convex subset $C$ of $\mathbb{R}^{n}$.
Then the Bregman projection of $x$ onto $C$ with respect to $f\in\mathscr{B}(\mathbb{R}^n)$ is given by
\begin{eqnarray}\label{Breg:projection}
\begin{aligned}P_{C}^{f}(x):=\arg\min_{y\in C}D_{f}(y,x).
\end{aligned}
\end{eqnarray}

It is well-known that the Bregman projection exists, and that it is
unique when $f\in\mathscr{B}(\mathbb{R}^n)$. Different types of Bregman distances and proximity functions have been considered in various contexts in the literature. See, e.g., \cite{HHBA2006,RSBU2021,CBYR2001}.

We recall the following result for Bregman projections.

\begin{theorem}\cite{YCEN1997}\label{Breg:Proj:Characterization}
Let $f\in\mathscr{B}(\mathbb{R}^n)$. Then for any $x\in\mathbb{R}^{n}$,
$x^{*}=P_{C}^{f}(x)$ if and only if
\begin{eqnarray}
\begin{aligned}
\left\langle \nabla f(x)-\nabla f(x^{*}),y-x^{*}\right\rangle \leq0,\;\text{for all}\;y\in C.
\end{aligned}
\end{eqnarray}
\end{theorem}

Unless otherwise stated, by $g_j\in\mathscr{B}(\mathbb{R}^{m_j})$, we mean that $g_j:\mathbb{R}^{m_j}\rightarrow\mathbb{R}$ is a Bregman function with zone $\mathbb{R}^{m_j}$, $j=1,2,\dots, t$. We denote $g:=(g_j)_{j=1}^t=(g_1, g_2,\dots,g_t).$
In Propositions~\ref{propo1:tran:breg}--\ref{propo3:tran:breg:diff} and Corollary~\ref{cor1:tran:breg:cont}, where a single point-to-set operator is considered, the letter $g$ denotes a \emph{translation operator} $g:\mathbb{R}^n\to\mathbb{R}^n$; this should not be confused with the Bregman functions $g_j$.
We also recall that, for $f\in\mathscr{B}(\mathbb{R}^n)$, a nonempty, closed, and convex set $C$, $x\in\mathbb{R}^n$ and $z\in C$, the Bregman projection satisfies the generalized Pythagorean inequality \cite[Chapter 2]{YCEN1997}
\begin{equation}\label{eq:three:point}
D_f\big(P_C^f(x),x\big)+D_f\big(z,P_C^f(x)\big)\leq D_f(z,x),
\end{equation}
which follows from Theorem~\ref{Breg:Proj:Characterization} and the three-point identity $D_f(z,x)=D_f(z,y)+D_f(y,x)+\langle\nabla f(y)-\nabla f(x),z-y\rangle$.

\subsection{Analytical properties of the generalized proximity function and an algorithm}
\label{subsec:Our-Bregman-proximity}

In this section, we aim to develop a gradient descent-type algorithm to solve the original problem \eqref{MSMSICFP}.
The construction of the algorithm is based on a generalized proximity function associated with the problem \eqref{MSMSICFP} introduced in \cite[Subsection 3.3.1]{SDYCEN2026}.
Specifically, for $f\in\mathscr{B}(\mathbb{R}^n)$ and $g_j\in\mathscr{B}(\mathbb{R}^{m_j})$, $j=1,2,\dots, t$, the function $p_{f,g}:\mathbb{R}^{n}\rightarrow\mathbb{R}$ defined by
\begin{eqnarray}\label{Breg:MainProxFunciton}
\begin{aligned}p_{f,g}(x):=\sum_{i=1}^{r}\alpha_{i}D_{f}(P_{C_{i}(x)}^{f}(x),x)+\sum_{j=1}^{t}\sum_{h=1}^{\ell_{j}}\beta_{h}^{j}D_{g_j}\bigg(P_{Q_{h}^{j}(A_{j}x)}^{g_j}(A_{j}x),A_{j}x\bigg),\;x\in\mathbb{R}^{n},
\end{aligned}
\end{eqnarray}
where $\alpha_{i}>0,\beta_{h}^{j}>0$ are parameters, and $C_{i},Q_{h}^{j}$
are the point-to-set operators defined above in the problem \eqref{MSMSICFP} for $h,j,i$,
is a generalized proximity function for it, in the sense of the following lemma.

\begin{lemma}\label{lem:proximity}
Let $f\in\mathscr{B}(\mathbb{R}^n)$, let $C\subseteq\mathbb{R}^n$ be nonempty, closed, and convex, and let $x\in\mathbb{R}^n$. Then $D_f(P^f_C(x),x)\geq0$, and $D_f(P^f_C(x),x)=0$ if and only if $x\in C$. Consequently, $p_{f,g}\geq0$, and $p_{f,g}(x)=0$ if and only if $x$ solves \eqref{MSMSICFP}.
\end{lemma}

\begin{proof}
Since $f$ is strictly convex and differentiable, $D_f(y,x)\geq0$ for all $y$, with equality if and only if $y=x$. If $x\in C$, then $x$ minimizes $D_f(\cdot,x)$ over $C$, so $P^f_C(x)=x$ and $D_f(P^f_C(x),x)=0$. Conversely, if $D_f(P^f_C(x),x)=0$, then $x=P^f_C(x)\in C$. The second assertion follows by applying the first one to $C=C_i(x)$ and to $C=Q_h^j(A_jx)$ (with $g_j$ and $A_jx$ in place of $f$ and $x$), since all the weights are positive.
\end{proof}

To facilitate the development of the algorithm, we first present several key results.

\begin{proposition}\label{propo1:tran:breg}
    Let $C:\mathbb{R}^n\rightarrow 2^{\mathbb{R}^n}$ be a point-to-set operator such that $C(x)$ is a nonempty, closed, and convex subset of $\mathbb{R}^n$ for every $x\in\mathbb{R}^n$. Let $f\in\mathscr{B}(\mathbb{R}^n)$ be $\sigma$-strongly convex. Assume that there exists an operator $g:\mathbb{R}^n\rightarrow\mathbb{R}^n$ and a nonempty, closed, and convex subset $\varOmega$ of $\mathbb{R}^n$ such that
    \begin{eqnarray}\label{ineq1:propo1:tran:breg}
        C(x)=g(x)+\varOmega,\;\text{for all}\;x\in\mathbb{R}^n.
    \end{eqnarray}
    Then for $z\in\mathbb{R}^n$
    \begin{eqnarray}\label{ineq2:propo1:tran:breg}
        P_{C(x)}^f(z)=g(x)+\omega^*,
    \end{eqnarray}
    where
    \begin{eqnarray}\label{ineq3:propo1:tran:breg}
        \omega^*=\arg\min_{\omega\in\varOmega}D_f(g(x)+\omega,z).
    \end{eqnarray}
\end{proposition}

\begin{proof}
    Let $x\in\mathbb{R}^n$ be fixed. For every $z\in\mathbb{R}^n$, the point $x^*= P_{C(x)}^f(z)$ exists and is unique since $f\in\mathscr{B}(\mathbb{R}^n)$ (see \cite[Chapter 2]{YCEN1997}). Since $x^*\in C(x)$, we have $x^*=g(x)+\omega^*$ with $\omega^*:=x^*-g(x)\in\varOmega$, and, since $\omega\mapsto g(x)+\omega$ maps $\varOmega$ bijectively onto $C(x)$, $\omega^*$ is the unique minimizer of $H(\omega):=D_f(g(x)+\omega, z)$ over $\varOmega$. (Independently of the theory of Bregman projections, the existence and uniqueness of this minimizer also follow from the facts that $H$ is continuous, hence closed, and, since $f$ is strongly convex, coercive and strictly convex; see \cite[Definition 2.13 and Theorem 2.14]{ABEC2017}.)
    
    By Theorem \ref{Breg:Proj:Characterization}, we obtain
    \begin{eqnarray}\label{ineq4:propo1:tran:breg}
        \left\langle\nabla f(z)-\nabla f(x^*),y-x^*\right\rangle\leq 0,\;\forall y\in C(x).
    \end{eqnarray}
    This implies that
    \begin{eqnarray}\label{ineq5:propo1:tran:breg}
    \left\langle\nabla f(z)-\nabla f(g(x)+\omega^*),g(x)+v-(g(x)+\omega^*)\right\rangle\leq 0,\;\forall v\in \varOmega,
    \end{eqnarray}
    \begin{eqnarray}\label{ineq6:propo1:tran:breg}
    \left\langle\nabla f(z)-\nabla f(g(x)+\omega^*),v-\omega^*\right\rangle\leq 0,\;\forall v\in \varOmega,
    \end{eqnarray}

    Define $h(\omega):=f(g(x)+\omega)$, $\omega\in\mathbb{R}^n$, which is a Bregman function, being a translate of $f$. Then $\nabla h(\omega)=\nabla f(g(x)+\omega)$. Therefore, from \eqref{ineq6:propo1:tran:breg}, we get
    \begin{eqnarray}\label{ineq7:propo1:tran:breg}
    \left\langle\nabla h(z-g(x))-\nabla h(\omega^*),v-\omega^*\right\rangle\leq 0,\;\forall v\in \varOmega
    \end{eqnarray}
    and by Theorem \ref{Breg:Proj:Characterization}, we obtain
    \begin{eqnarray}\label{ineq8:propo1:tran:breg}
    \omega^*=P_\varOmega^h(z-g(x)).
    \end{eqnarray}
    Hence,
    \begin{eqnarray}\label{ineq9:propo1:tran:breg}
        P_{g(x)+\varOmega}^f(z)=g(x)+P_\varOmega^h(z-g(x)),
    \end{eqnarray}
    where
    \begin{eqnarray}\label{ineq10:propo1:tran:breg}
        P_\varOmega^h(z-g(x))&=&\arg\min_{\omega\in\varOmega}D_h(\omega,z-g(x))\nonumber\\
        &=&\arg\min_{\omega\in\varOmega}D_f(g(x)+\omega,z).
    \end{eqnarray}
\end{proof}

The following example computes the Bregman projection for the translation-type point-to-set operator for the general quadratic Bregman function.
\begin{example}
    Let $Q$ be a symmetric and positive definite matrix of order $n$ and $\varOmega=\mathbb{R}^n_+$. Consider the general quadratic function
    \begin{eqnarray}
        f(x)=\frac{1}{2}\left\langle Qx,x\right\rangle+\left\langle c,x\right\rangle+d,
    \end{eqnarray}
    where $c\in\mathbb{R}^n$ and $d\in\mathbb{R}$. Then it is easy to check that $f\in\mathscr{B}(\mathbb{R}^n)$ and $\nabla f(z)=Qz+c$. Assume that $C(x)=g(x)+\varOmega$, where $g:\mathbb{R}^n\rightarrow\mathbb{R}^n$ is an operator. Define $\|u\|_Q^2:=\left\langle Qu,u\right\rangle$ for $u\in\mathbb{R}^n$.
    Now,
    \begin{eqnarray}
        D_f(g(x)+\omega,z)&=&\frac{1}{2}\left\langle Q(g(x)+\omega), g(x)+\omega\right\rangle+\left\langle c, g(x)+\omega\right\rangle\nonumber\\&&-\frac{1}{2}\left\langle Qz,z\right\rangle-\left\langle c,z\right\rangle- \left\langle Qz+c, g(x)+\omega-z\right\rangle\\
        &=&\frac{1}{2}\|\omega-(z-g(x))\|^2_Q.\nonumber
    \end{eqnarray}
    The gradient of $D_f(g(x)+\omega,z)$ is
    \begin{eqnarray}
         \nabla_\omega D_f(g(x)+\omega,z)=Q(\omega+g(x))-Qz.
    \end{eqnarray}
    Now, $\nabla_\omega D_f(g(x)+\omega,z)=0$ implies $\omega=z-g(x)$ and $\nabla^2_\omega D_f(g(x)+\omega,z)=Q$.
    The function $\omega\mapsto\frac{1}{2}\|\omega-(z-g(x))\|^2_Q$ is continuous and coercive. Therefore, it attains a minimum on $\mathbb{R}^n_+$. Since $Q$ is symmetric and positive definite, the minimizer is unique.
    Therefore,
    \begin{eqnarray}
        \omega^*\in\arg\min_{\omega\in\mathbb{R}^n_+}D_f(g(x)+\omega,z)=\arg\min_{\omega\in\mathbb{R}^n_+}\frac{1}{2}\|\omega-(z-g(x))\|^2_Q:=P_{\mathbb{R}^n_+}^Q(z-g(x)).
    \end{eqnarray}
    This implies that $\omega^*=P_{\mathbb{R}^n_+}^Q(z-g(x))$.
    Hence,
    \begin{eqnarray}
        P_{C(x)}^f(z)=g(x)+P_{\mathbb{R}^n_+}^Q(z-g(x)).
    \end{eqnarray}
\end{example}

\begin{remark}\label{rmk1:tran:breg}
    If we set $f(x)=(1/2)\|x\|^2$, in Proposition \ref{propo1:tran:breg}, then the formula in \eqref{ineq2:propo1:tran:breg} reduces to
    \begin{eqnarray}
        P_{g(x)+\varOmega}(z)=g(x)+P_\varOmega(z-g(x)).
    \end{eqnarray}
    This is a well-known formula (see, for instance, \cite{YNES2011}).
\end{remark}

\begin{remark} Let $C:\mathbb{R}^n\rightarrow 2^{\mathbb{R}^n}$ be a point-to-set operator such that $C(x)$ is a nonempty, closed, and convex subset of $\mathbb{R}^n$ for each $x\in\mathbb{R}^n$ and $g:\mathbb{R}^n\rightarrow {\mathbb{R}^n}$ be an operator.  Assume that $C(x)=g(x)+\varOmega$, where $\varOmega$ is a nonempty, closed, and convex subset of $\mathbb{R}^n$. Then it is well-known that in the case the orthogonal projection for any $x,y,z\in\mathbb{R}^n$
    \begin{eqnarray}
        \|P_{C(x)}(z)-P_{C(y)}(z)\|\leq L\|x-y\|,
    \end{eqnarray}
    where $L$ is the Lipschitz constant of $g$; indeed, $P_{C(x)}(z)=z-(I-P_\varOmega)(z-g(x))$ and $I-P_\varOmega$ is nonexpansive. This estimate fails for Bregman projections, as the following example shows.
\end{remark}

\begin{example}
    Let $\varOmega=\{0\}\times [0,1]\subset\mathbb{R}^2$, $g(x)=x$, and $f(x)=(1/2)(x_1^2+x_1x_2+x_2^2)$. Then $C(x)=\{(x_1,x_2+u)\vert u\in [0,1]\}$ and $f\in\mathscr{B}(\mathbb{R}^n)$. For $z=(0,0)$, we obtain $D_f(\omega,z)=(1/2)(x_1^2+x_1(x_2+u)+(x_2+u)^2)$ for $\omega\in C(x)$. Suppose $x=(1/6,-1/2)$. Then the minimizer of $D_f((1/6,-1/2+u),0)$ is $(0,5/12)$. Therefore,
    \begin{eqnarray}
        P^f_{C(x)}(z)=(1/6,-1/2)+(0,5/12)=(1/6,-1/12).
    \end{eqnarray}
    Suppose $y=(0,-1/2)$. Then the minimizer of $D_f((0,-1/2+u),0)$ is $(0,1/2)$. Therefore,
    \begin{eqnarray}
        P^f_{C(y)}(z)=(0,-1/2)+(0,1/2)=(0,0).
    \end{eqnarray}
    Therefore,
    \begin{eqnarray}
        \| P^f_{C(x)}(z)- P^f_{C(y)}(z)\|=\sqrt{5}/12>1/6=\|x-y\|.
    \end{eqnarray}
    However, $g$ is Lipschitz continuous with constant $L=1$.
\end{example}

\begin{proposition}\label{propo2:tran:breg:cont}
    Let $C:\mathbb{R}^n\rightarrow 2^{\mathbb{R}^n}$ be a continuous point-to-set operator, that is, $C$ is both inner and outer semi-continuous at every point of $\mathbb{R}^n$, such that $C(x)$ is a nonempty, closed, and convex subset of $\mathbb{R}^n$ for each $x\in\mathbb{R}^n$. Let $f\in\mathscr{B}(\mathbb{R}^n)$ be $\sigma$-strongly convex. Then
    \begin{eqnarray}\label{ineq1:propo2:tran:breg:cont}
        x\mapsto P_{C(x)}^f(x)
    \end{eqnarray}
    is continuous.
\end{proposition}

\begin{proof}
    Let $x_k\rightarrow x^*$ as $k\rightarrow\infty$. Define $y_k:= P_{C(x_k)}^f(x_k)$. Let $w^*\in C(x^*)$. Since $C$ is inner semi-continuous,
    there exists $w_k\in C(x_k)$ such that $w_k\rightarrow w^*$. Since $y_k$ is the Bregman projection of $x_k$ onto $C(x_k)$,
    \begin{eqnarray}\label{ineq2:propo2:tran:breg:cont}
        D_f(y_k,x_k)\leq D_f(w_k,x_k),\;\text{for all}\; k.
    \end{eqnarray}
    Since $f$ is $\sigma$-strongly convex,
    \begin{eqnarray}\label{ineq2.1:propo2:tran:breg:cont}
        \frac{\sigma}{2}\|y_k-x_k\|^2\leq D_f(y_k,x_k).
    \end{eqnarray}
    Since $x_k\rightarrow x^*$ and $w_k\rightarrow w^*$, the sequence $\{D_f(w_k,x_k)\}$ is bounded.
    Moreover, from \eqref{ineq2:propo2:tran:breg:cont} and \eqref{ineq2.1:propo2:tran:breg:cont}, we guarantee that $\{\|y_k-x_k\|\}$ is bounded.
    Thus, the sequence $\{y_k\}$ is bounded. Therefore, it has a convergent subsequence. Let $\{y_{k_\ell}\}$ be a subsequence of $\{y_k\}$ such that $y_{k_\ell}\rightarrow y^*$ as $\ell\rightarrow\infty$. Since $C$ is outer semi-continuous and $y_{k_\ell}\in C(x_{k_\ell})$ with $x_{k_\ell}\rightarrow x^*$ as $\ell\rightarrow\infty$, $y^*\in C(x^*)$. Let $u\in C(x^*)$. Since $C$ is inner semi-continuous, there exists $u_{k_\ell}\in C(x_{k_\ell})$ such that $u_{k_\ell}\rightarrow u$ as $\ell\rightarrow\infty$. Now,
    \begin{eqnarray}\label{ineq3:propo2:tran:breg:cont}
        D_f(y_{k_\ell},x_{k_\ell})\leq D_f(u_{k_\ell},x_{k_\ell}).
    \end{eqnarray}
    Taking the limit as $\ell\rightarrow\infty$ on both sides, we get
    \begin{eqnarray}\label{ineq4:propo2:tran:breg:cont}
        D_f(y^*,x^*)\leq D_f(u,x^*),\;\text{for all}\; u\in C(x^*).
    \end{eqnarray}
    This implies $y^*=P_{C(x^*)}^f(x^*)$.
    The same argument shows that every subsequence of $\{y_k\}$ has a further subsequence which converges to $P_{C(x^*)}^f(x^*)$. Therefore, $y_k\rightarrow P_{C(x^*)}^f(x^*)$ as $k\rightarrow\infty$.
\end{proof}

\begin{remark}\label{remark1:propo2:tran:breg:cont}
    It is important to note that only one-sided continuity of $C$ in Proposition \ref{propo2:tran:breg:cont}, i.e., either inner semi-continuity or outer semi-continuity, does not guarantee the continuity of $x\mapsto P_{C(x)}^f(x)$. The following two examples support this claim.
\end{remark}

\begin{example}\label{example1:remark1:propo2:tran:breg:cont}
    Let $f(x_1,x_2)=(1/2)(x_1^2+x_1x_2+x_2^2)$ on $\mathbb{R}^2$ and
    \begin{equation}\label{ineq1:example1:remark1:propo2:tran:breg:cont}
    C(x_1,x_2)=
    \begin{cases}
        \{(1,0)\},& (x_1,x_2)= (0,0)\\
        [0,1]\times\{0\},& (x_1,x_2)\neq (0,0).
    \end{cases}
    \end{equation}
    One can check that it is inner semi-continuous on $\mathbb{R}^2.$ Indeed,
    let $x_k\rightarrow x^*$ as $k\rightarrow\infty$ and $y^*\in C(x^*)$. If $x^*=(0,0)$ then take $y_k=y^*=(1,0)$ for all $k$. Otherwise, for any $y^*\in C(x^*)$ take $y_k=y^*$ for all $k$.
    However, it is not outer semi-continuous at $(0,0)$. Indeed, let $x_k=(1/(k+1),0)$. Then $x_k\rightarrow (0,0)$ as $k\rightarrow \infty.$ Taking $y_k=(0,0)\in C(x_k)$ for all $k$, we observe that $(0,0)\notin C(0,0).$

    Moreover,
    \begin{eqnarray}\label{ineq2:example1:remark1:propo2:tran:breg:cont}
        D_f(y,x)=f(y)-f(x)-\left\langle \nabla f(x),y-x\right\rangle,
    \end{eqnarray}
    where
    \begin{equation}\label{ineq2b:example1:remark1:propo2:tran:breg:cont}
    \nabla f(x)=
    \begin{pmatrix}
    x_1+\frac{x_2}{2}\\
    x_2+\frac{x_1}{2}
    \end{pmatrix}.
    \end{equation}
    Therefore,
    \begin{eqnarray}\label{ine3:example1:remark1:propo2:tran:breg:cont}
        D_f((y_1,0),x)=\frac{1}{2}((y_1-x_1)^2-x_2(y_1-x_1)+x_2^2).
    \end{eqnarray}
    This implies that $P_{C(x)}^f(x)=(\min\{1,\max\{0,x_1+(x_2/2)\}\},0)$, i.e.,
    \begin{equation}\label{ineq4:example1:remark1:propo2:tran:breg:cont}
    P_{C(x)}^f(x)=
    \begin{cases}
        (0,0), \;& x_1+\frac{x_2}{2}< 0\\
        (x_1+\frac{x_2}{2},0), \;& 0\leq x_1+\frac{x_2}{2}\leq 1\\
        (1,0),& x_1+\frac{x_2}{2}> 1.
    \end{cases}
    \end{equation}
    Let $x=(0,0)$ and $x_k=(1/(k+1),0)$. Then $x_k\rightarrow x$ as $k\rightarrow\infty$. However,
    \begin{eqnarray}\label{ineq5:example1:remark1:propo2:tran:breg:cont}
        P_{C(x_k)}^f(x_k)=(1/(k+1),0)\rightarrow (0,0)\neq (1,0)=P_{C(x)}^f(x)\;\text{as}\; k\rightarrow\infty.
    \end{eqnarray}
\end{example}

\begin{example}\label{example2:remark1:propo2:tran:breg:cont}
    Let $f(x_1,x_2)=(1/2)(x_1^2+x_1x_2+x_2^2)$ on $\mathbb{R}^2$ and
    \begin{equation}\label{ineq1:example2:remark1:propo2:tran:breg:cont}
    C(x_1,x_2)=
    \begin{cases}
        [0,1]\times\{0\},& (x_1,x_2)= (0,0)\\
        \{(1,0)\},& (x_1,x_2)\neq (0,0).
    \end{cases}
    \end{equation}
    One can check that it is outer semi-continuous on $\mathbb{R}^2.$ Indeed,
    let $x_k\rightarrow x^*$ as $k\rightarrow\infty$ and $y_k\in C(x_k)$ such that $y_k\rightarrow y^*$. Let $x^*=(0,0)$. If $x_k\neq (0,0)$ then $C(x_k)=\{(1,0)\}$ for all $k$. Therefore, $y^*\in C(x^*)$. Otherwise, if for infinitely many $k$, $x_k$ equals $(0,0)$, the corresponding $y_k$ belongs to $[0,1]\times\{0\}$. In either case, the possible limit $y^*$ belongs to $[0,1]\times\{0\}=C(x^*).$
    However, it is not inner semi-continuous at $x^*=(0,0)$. Indeed, let $x_k=(1/(k+1),0)$. Then $x_k\rightarrow (0,0)$ as $k\rightarrow \infty.$ Let $y^*=(0,0)$. Then $y^*\in C(x^*)$. But $y_k\in C(x_k)$ does not converge to $y^*$.

    Moreover,
    \begin{eqnarray}\label{ineq2:example2:remark1:propo2:tran:breg:cont}
        D_f((y_1,0),x)=\frac{1}{2}((y_1-x_1)^2-x_2(y_1-x_1)+x_2^2).
    \end{eqnarray}
    This implies that
    \begin{equation}\label{ineq3:example2:remark1:propo2:tran:breg:cont}
    P_{C(x)}^f(x)=
    \begin{cases}
        (0,0), \;& (x_1,x_2)=(0,0)\\
        (1,0), \;& (x_1,x_2)\neq (0,0).
    \end{cases}
    \end{equation}
    Let $x=(0,0)$ and $x_k=(1/(k+1),0)$. Then $x_k\rightarrow x$ as $k\rightarrow\infty$. However,
    \begin{eqnarray}\label{ineq4:example2:remark1:propo2:tran:breg:cont}
        P_{C(x_k)}^f(x_k)=(1,0)\nrightarrow (0,0)=P_{C(x)}^f(x)\;\text{as}\; k\rightarrow\infty.
    \end{eqnarray}
\end{example}

\begin{corollary}\label{cor1:tran:breg:cont}
    Let $C:\mathbb{R}^n\rightarrow 2^{\mathbb{R}^n}$ be a point-to-set operator such that $C(x)$ is a nonempty, closed, and convex subset of $\mathbb{R}^n$ for every $x\in\mathbb{R}^n$. Let $f\in\mathscr{B}(\mathbb{R}^n)$ be $\sigma$-strongly convex. Assume that there exists a continuous operator $g:\mathbb{R}^n\rightarrow {\mathbb{R}^n}$ and $\varOmega$ is a nonempty, closed, and convex subset of $\mathbb{R}^n$ such that
    \begin{eqnarray}\label{ineq1:cor1:tran:breg:cont}
        C(x)=g(x)+\varOmega,\;\text{for all}\;x\in\mathbb{R}^n.
    \end{eqnarray}
    Then the function
    \begin{eqnarray}\label{ineq2:cor1:tran:breg:cont}
        x\mapsto P_{C(x)}^f(x)
    \end{eqnarray}
    is continuous.
\end{corollary}

\begin{proof}
Since $g$ is continuous, the point-to-set operator $x\mapsto C(x)$ is continuous.
Indeed, let $x_k\rightarrow x$ as $k\rightarrow\infty$ and $y_k\in C(x_k)$ such that $y_k\rightarrow y$ as $k\rightarrow\infty$. Since $g$ is continuous, $y_k-g(x_k)\rightarrow y-g(x)$ as $k\rightarrow\infty$. Again since $y_k-g(x_k)\in\varOmega$ and $\varOmega$ is closed, $y-g(x)\in \varOmega$. This implies $y\in C(x)$. Thus, $x\mapsto C(x)$ is outer semi-continuous.

To show the inner semi-continuity of $x\mapsto C(x)$, let $x_k\rightarrow x$ as $k\rightarrow\infty$ and $y\in C(x)$. Then $\omega:=y-g(x)\in\varOmega$, and $y_k:=g(x_k)+\omega\in C(x_k)$. Since $g$ is continuous, $y_k\rightarrow g(x)+\omega=y$ as $k\rightarrow\infty$.
Thus, $x\mapsto C(x)$ is inner semi-continuous.
Therefore, the proof directly follows from Proposition \ref{propo2:tran:breg:cont}.
\end{proof}

\begin{proposition}\label{propo1:tran:breg:lipschitz}
    Let $C:\mathbb{R}^n\rightarrow 2^{\mathbb{R}^n}$ be a point-to-set operator such that $C(x)$ is a nonempty, closed, and convex subset of $\mathbb{R}^n$ for every $x\in\mathbb{R}^n$. Let $f\in\mathscr{B}(\mathbb{R}^n)$ be $\sigma$-strongly convex and its gradient is $L_f$-Lipschitz continuous. Assume that there exists a Lipschitz continuous operator $g:\mathbb{R}^n\rightarrow {\mathbb{R}^n}$ with constant $L_g$ and $\varOmega$ is a nonempty, closed, and convex subset of $\mathbb{R}^n$ such that
    \begin{eqnarray}\label{ineq1:propo1:tran:breg:lipschitz}
        C(x)=g(x)+\varOmega,\;\text{for all}\;x\in\mathbb{R}^n.
    \end{eqnarray}
    Then the function
    \begin{eqnarray}\label{ineq2:propo1:tran:breg:lipschitz}
        x\mapsto P_{C(x)}^f(x)
    \end{eqnarray}
    is Lipschitz continuous with constant
    \begin{eqnarray}\label{ineq3:propo1:tran:breg:lipschitz}
        L_g+\frac{L_f}{\sigma}(L_g+1).
    \end{eqnarray}
\end{proposition}

\begin{proof}
    Fix $x,y\in\mathbb{R}^n$. Let $u_x=P_{C(x)}^f(x)=g(x)+\omega_x$ and $u_y=P_{C(y)}^f(y)=g(y)+\omega_y$ for some $\omega_x,\omega_y\in\varOmega$. Using Theorem \ref{Breg:Proj:Characterization}, we get
    \begin{eqnarray}\label{ineq4:propo1:tran:breg:lipschitz}
        \left\langle \nabla f(x)-\nabla f(u_x), \omega-\omega_x\right\rangle\leq 0,\;\text{for all}\;\omega\in\varOmega
    \end{eqnarray}
    and
    \begin{eqnarray}\label{ineq5:propo1:tran:breg:lipschitz}
        \left\langle \nabla f(y)-\nabla f(u_y), \omega^\prime-\omega_y\right\rangle\leq 0,\;\text{for all}\;\omega^\prime\in\varOmega.
    \end{eqnarray}
    Taking $\omega=\omega_y$ in \eqref{ineq4:propo1:tran:breg:lipschitz} and $\omega^\prime=\omega_x$ in \eqref{ineq5:propo1:tran:breg:lipschitz} and adding, we get
    \begin{eqnarray}\label{ineq6a:propo1:tran:breg:lipschitz}
        \left\langle \nabla f(x)-\nabla f(y)+\nabla f(u_y)-\nabla f(u_x), \omega_y-\omega_x\right\rangle\leq 0.
    \end{eqnarray}
   Now, writing $\nabla f(u_x)-\nabla f(u_y)=\big(\nabla f(g(x)+\omega_x)-\nabla f(g(x)+\omega_y)\big)+\big(\nabla f(g(x)+\omega_y)-\nabla f(g(y)+\omega_y)\big)$, we obtain from \eqref{ineq6a:propo1:tran:breg:lipschitz}
    \begin{eqnarray}\label{ineq6:propo1:tran:breg:lipschitz}
    \begin{aligned}
        &\left\langle \nabla f(g(x)+\omega_x)-\nabla f(g(x)+\omega_y), \omega_x-\omega_y\right\rangle \\&\hspace{0.5cm} \leq \left\langle \nabla f(g(y)+\omega_y)-\nabla f(g(x)+\omega_y),\omega_x-\omega_y\right\rangle\\&\hspace{1cm} +\left\langle \nabla f(x)-\nabla f(y), \omega_x-\omega_y\right\rangle.
    \end{aligned}
    \end{eqnarray}
    Since $\nabla f$ is $\sigma$-strongly monotone (Lemma~\ref{lemma:strongconvex}), the left-hand side of \eqref{ineq6:propo1:tran:breg:lipschitz} is bounded below by $\sigma\|\omega_x-\omega_y\|^2$, while, by the Cauchy--Schwarz inequality and the Lipschitz continuity of $\nabla f$ and $g$, its right-hand side is bounded above by $L_f(L_g+1)\|x-y\|\|\omega_x-\omega_y\|$. Hence
    \begin{eqnarray}\label{ineq7:propo1:tran:breg:lipschitz}
    \begin{aligned}
        \sigma\|\omega_x-\omega_y\|^2 \leq L_f(L_g+1)\|x-y\|\|\omega_x-\omega_y\|.
    \end{aligned}
    \end{eqnarray}
    This implies that
    \begin{eqnarray}\label{ineq8:propo1:tran:breg:lipschitz}
    \begin{aligned}
        \|\omega_x-\omega_y\| \leq \frac{L_f}{\sigma}(L_g+1)\|x-y\|.
    \end{aligned}
    \end{eqnarray}
    Therefore,
    \begin{eqnarray}\label{ineq9:propo1:tran:breg:lipschitz}
    \begin{aligned}
        \|P_{C(x)}^f(x)-P_{C(y)}^f(y)\|&\leq\|g(x)-g(y)\|+\|\omega_x-\omega_y\|\\&
        \leq \Big(L_g+\frac{L_f}{\sigma}(L_g+1)\Big)\|x-y\|.
    \end{aligned}
    \end{eqnarray}
\end{proof}

\begin{proposition}\label{propo2:tran:breg:lipschitz}
    Let $A:\mathbb{R}^n\rightarrow\mathbb{R}^m$ be a bounded linear operator, let $\tilde{g}\in\mathscr{B}(\mathbb{R}^m)$ be $\sigma_{\tilde{g}}$-strongly convex with $L_{\tilde{g}}$-Lipschitz continuous gradient, and let $Q(y)=q(y)+\tilde{\varOmega}$ for all $y\in\mathbb{R}^m$, where $q:\mathbb{R}^m\rightarrow\mathbb{R}^m$ is $L_q$-Lipschitz continuous and $\tilde{\varOmega}\subseteq\mathbb{R}^m$ is nonempty, closed, and convex. Then the function
    \begin{eqnarray}\label{ineq2:propo2:tran:breg:lipschitz}
        x\mapsto P_{Q(Ax)}^{\tilde{g}}(Ax)
    \end{eqnarray}
    is Lipschitz continuous with constant
    \begin{eqnarray}\label{ineq3:propo2:tran:breg:lipschitz}
        \Big(L_q+\frac{L_{\tilde{g}}}{\sigma_{\tilde{g}}}(L_q+1)\Big)\|A\|_2.
    \end{eqnarray}
\end{proposition}

\begin{proof}
    By Proposition~\ref{propo1:tran:breg:lipschitz} (applied in $\mathbb{R}^m$), the map $\Phi(y):=P^{\tilde g}_{Q(y)}(y)$ is Lipschitz continuous with constant $L_q+\frac{L_{\tilde{g}}}{\sigma_{\tilde{g}}}(L_q+1)$. Since $\|Ax-Ay\|\leq\|A\|_2\|x-y\|$, the composition $x\mapsto\Phi(Ax)$ is Lipschitz continuous with the constant \eqref{ineq3:propo2:tran:breg:lipschitz}.
\end{proof}

\begin{proposition}\label{propo3:tran:breg:diff}
    Let $C:\mathbb{R}^n\rightarrow 2^{\mathbb{R}^n}$ be a point-to-set operator such that $C(x)$ is a nonempty, closed, and convex subset of $\mathbb{R}^n$ for every $x\in\mathbb{R}^n$. Let $f\in\mathscr{B}(\mathbb{R}^n)$ be $\sigma$-strongly convex and suppose $f\in \mathcal{C}^2(\mathbb{R}^n)$. Assume that there exists a $g\in\mathcal{C}^1(\mathbb{R}^n;\mathbb{R}^n)$, the space of all continuously differentiable operators on $\mathbb{R}^n$ (that is, the Jacobian matrix $J_g(x)$ exists and is continuous on $\mathbb{R}^n$), and a nonempty, closed, and convex subset $\varOmega$ of $\mathbb{R}^n$ such that
    \begin{eqnarray}\label{ineq11:propo3:tran:breg:cont}
        C(x)=g(x)+\varOmega,\;\text{for all}\;x\in\mathbb{R}^n.
    \end{eqnarray}
    Then the function
    \begin{eqnarray}\label{ineq1:propo3:tran:breg:cont}
        x\mapsto D_f(P_{C(x)}^f(x),x)
    \end{eqnarray}
    is continuously differentiable and its gradient is given by
    \begin{eqnarray}\label{ineq2:propo3:tran:breg:cont}
        \nabla D_f(P_{C(x)}^f(x),x)=J_g(x)^T(\nabla f(P_{C(x)}^f(x))-\nabla f(x))-\nabla^2 f(x)(P_{C(x)}^f(x)-x),
    \end{eqnarray}
    where $J_g$ is the Jacobian of $g$ and $\nabla^2 f$ is the Hessian of $f$.
\end{proposition}

\begin{proof}
    Let $\phi(x,\omega):=D_f(g(x)+\omega,x)=f(g(x)+\omega)-f(x)-\left\langle\nabla f(x),g(x)+\omega-x\right\rangle$ for $x,\omega\in\mathbb{R}^n$. Then
    \begin{eqnarray}\label{ineq4:propo3:tran:breg:cont}
        D_f(P_{C(x)}^f(x),x)=\min_{y\in C(x)} D_f(y,x)=\min_{\omega\in \varOmega} \phi(x,\omega).
    \end{eqnarray}
    Since $f$ is strongly convex, $\phi(x,\cdot)$ is strongly convex, so it has a unique minimizer $\bar\omega(x)$ over $\varOmega$, and $P_{C(x)}^f(x)=g(x)+\bar{\omega}(x)$ by Proposition~\ref{propo1:tran:breg}. Since $f\in\mathcal{C}^2$ and $g\in\mathcal{C}^1$, the function $\phi(\cdot,\omega)$ is differentiable, and a direct computation gives
    \begin{eqnarray}\label{ineq7:propo3:tran:breg:cont}
       \nabla_x\phi(x,\omega)= J_g(x)^T\big(\nabla f(g(x)+\omega)-\nabla f(x)\big)-\nabla^2 f(x)\big(g(x)+\omega-x\big),
    \end{eqnarray}
    which is continuous on $\mathbb{R}^n\times\mathbb{R}^n$. Fix $x_0\in\mathbb{R}^n$ and an open ball $U$ centered at $x_0$. By Corollary~\ref{cor1:tran:breg:cont}, $x\mapsto\bar\omega(x)=P^f_{C(x)}(x)-g(x)$ is continuous, so there is $\rho>0$ with $\|\bar\omega(x)\|<\rho$ for all $x\in\overline{U}$. The set $\varOmega_\rho:=\{\omega\in\varOmega\mid\|\omega\|\leq\rho\}$ is nonempty, compact, and convex, and, for $x\in U$, $\min_{\omega\in\varOmega}\phi(x,\omega)=\min_{\omega\in\varOmega_\rho}\phi(x,\omega)$, the unique minimizer over $\varOmega_\rho$ being $\bar\omega(x)$. Applying Theorem~\ref{Danskin's:Thm} to $\psi:=-\phi$ on $\mathbb{R}^n\times\varOmega_\rho$, we conclude that $x\mapsto D_f(P_{C(x)}^f(x),x)$ is differentiable at $x_0$ with gradient $\nabla_x\phi(x_0,\bar\omega(x_0))$, which is \eqref{ineq2:propo3:tran:breg:cont} by \eqref{ineq7:propo3:tran:breg:cont}. Finally, the gradient is continuous by Corollary~\ref{cor1:tran:breg:cont} and the continuity of $J_g$, $\nabla f$ and $\nabla^2 f$.
\end{proof}

\begin{corollary}\label{cor2:tran:breg:diff}
    Let $A:\mathbb{R}^n\rightarrow {\mathbb{R}^m}$ be a bounded linear operator, let $\tilde g\in\mathscr{B}(\mathbb{R}^m)\cap\mathcal{C}^2(\mathbb{R}^m)$ be strongly convex, and let $Q(y)=q(y)+\tilde{\varOmega}$ for all $y\in\mathbb{R}^m$, where $q\in\mathcal{C}^1(\mathbb{R}^m;\mathbb{R}^m)$ and $\tilde{\varOmega}\subseteq\mathbb{R}^m$ is nonempty, closed, and convex. Set $\Phi(y):=D_{\tilde g}\big(P^{\tilde g}_{Q(y)}(y),y\big)$ for $y\in\mathbb{R}^m$. Then the function $\Phi\circ A:x\mapsto D_{\tilde g}\big(P_{Q(Ax)}^{\tilde g}(Ax),Ax\big)$ is continuously differentiable on $\mathbb{R}^n$ and, writing $u(x):=P_{Q(Ax)}^{\tilde g}(Ax)$,
    \begin{eqnarray}\label{ineq2:cor2:tran:breg:diff}
         \nabla(\Phi\circ A)(x)=A^T\Big[J_q(Ax)^T\big(\nabla \tilde g(u(x))-\nabla \tilde g(Ax)\big)-\nabla^2 \tilde g(Ax)\big(u(x)-Ax\big)\Big],
    \end{eqnarray}
    where $J_q$ is the Jacobian of $q$ and $\nabla^2\tilde g$ is the Hessian of $\tilde g$.
\end{corollary}

\begin{proof}
    Proposition~\ref{propo3:tran:breg:diff}, applied in $\mathbb{R}^m$ with the Bregman function $\tilde g$ in place of $f$ and the translation operator $q$ in place of $g$, shows that $\Phi$ is continuously differentiable with $\nabla\Phi(y)=J_q(y)^T\big(\nabla \tilde g(P^{\tilde g}_{Q(y)}(y))-\nabla \tilde g(y)\big)-\nabla^2\tilde g(y)\big(P^{\tilde g}_{Q(y)}(y)-y\big)$. Since $A$ is linear, the chain rule gives $\nabla(\Phi\circ A)(x)=A^T\nabla\Phi(Ax)$, which is \eqref{ineq2:cor2:tran:breg:diff}, and $\nabla(\Phi\circ A)$ is continuous as a composition of continuous maps.
\end{proof}

\begin{assumption}\label{ass:A}
\begin{itemize}
\item[(i)] $f\in\mathscr{B}(\mathbb{R}^n)\cap\mathcal{C}^2(\mathbb{R}^n)$ is $\sigma$-strongly convex and, for each $j$, $g_j\in\mathscr{B}(\mathbb{R}^{m_j})\cap\mathcal{C}^2(\mathbb{R}^{m_j})$ is $\sigma_j$-strongly convex;
\item[(ii)] for each $i$, $C_i(x)=\tilde g_i(x)+\varOmega_i$ for all $x\in\mathbb{R}^n$, where $\tilde g_i\in\mathcal{C}^1(\mathbb{R}^n;\mathbb{R}^n)$ and $\varOmega_i\subseteq\mathbb{R}^n$ is nonempty, closed, and convex;
\item[(iii)] for each $j$ and $h$, $Q_h^j(y)=p_h^j(y)+\tilde\varOmega_h^j$ for all $y\in\mathbb{R}^{m_j}$, where $p_h^j\in\mathcal{C}^1(\mathbb{R}^{m_j};\mathbb{R}^{m_j})$ and $\tilde\varOmega_h^j\subseteq\mathbb{R}^{m_j}$ is nonempty, closed, and convex;
\item[(iv)] the weights $\alpha_i$ and $\beta_h^j$ are positive.
\end{itemize}
\end{assumption}

\begin{proposition}\label{propo4:final:diffe:sum}
    Under Assumption~\ref{ass:A}, the function $p_{f,g}$ defined in \eqref{Breg:MainProxFunciton} is continuously differentiable and
    \begin{eqnarray}\label{ineq2:propo4:final:diffe:sum}
        \begin{aligned}
            \nabla p_{f,g}(x)&=\sum_{i=1}^{r}\alpha_{i}\Big(J_{\tilde{g}_i}(x)^T(\nabla f\Big(P_{C_i(x)}^f(x)\Big)-\nabla f(x))-\nabla^2 f(x)\Big(P_{C_i(x)}^f(x)-x\Big)\Big)\\&\hspace{0.5cm}+\sum_{j=1}^{t}\sum_{h=1}^{\ell_{j}}\beta_{h}^{j}\Big(A_j^TJ_{p_h^j}(A_jx)^T(\nabla g_j\Big(P_{Q_h^j(A_jx)}^{g_j}(A_jx)\Big)\\&\hspace{1cm}-\nabla g_j(A_jx))-A_j^T\nabla^2 g_j(A_jx)\Big(P_{Q_h^j(A_jx)}^{g_j}(A_jx)-A_jx\Big)\Big).
        \end{aligned}
    \end{eqnarray}
\end{proposition}

\begin{proof}
    Apply Proposition~\ref{propo3:tran:breg:diff} to each term $D_f(P^f_{C_i(x)}(x),x)$ and Corollary~\ref{cor2:tran:breg:diff} to each term $D_{g_j}(P^{g_j}_{Q_h^j(A_jx)}(A_jx),A_jx)$, and sum.
\end{proof}

\begin{proposition}\label{propo5:final:diffe:Lipschitz}
    Suppose that Assumption~\ref{ass:A} holds. Assume that
    \begin{eqnarray}\label{ineq1:propo5:final:diffe:Lipschitz}
        x\mapsto J_{\tilde{g}_i}(x)^T(\nabla f(P_{C_i(x)}^f(x))-\nabla f(x)),
    \end{eqnarray}
    \begin{eqnarray}\label{ineq2:propo5:final:diffe:Lipschitz}
        x\mapsto  \nabla^2 f(x)(P_{C_i(x)}^f(x)-x),
    \end{eqnarray}
    \begin{eqnarray}\label{ineq3:propo5:final:diffe:Lipschitz}
        x\mapsto  J_{p_h^j}(A_jx)^T(\nabla g_j(P_{Q_h^j(A_jx)}^{g_j}(A_jx))-\nabla g_j(A_jx)),
    \end{eqnarray}
    and
    \begin{eqnarray}\label{ineq4:propo5:final:diffe:Lipschitz}
        x\mapsto  \nabla^2 g_j(A_jx)(P_{Q_h^j(A_jx)}^{g_j}(A_jx)-A_jx),
    \end{eqnarray}
    are $L_1^i$, $M_1^i$, $L_2^{j,h}$, and $M_2^{j,h}$-Lipschitz continuous, respectively.
    Then the function
    \begin{eqnarray}\label{ineq5:propo5:final:diffe:Lipschitz}
        x\mapsto \nabla p_{f,g}(x)
    \end{eqnarray}
    is Lipschitz continuous with constant
    \begin{eqnarray}\label{ineq6:propo5:final:diffe:Lipschitz}
        \begin{aligned}
           L_{\nabla p_{f,g}}:=\sum_{i=1}^{r}\alpha_{i}(L_1^i+M_1^i)+\sum_{j=1}^{t}\sum_{h=1}^{\ell_{j}}\beta_{h}^{j}\|A_j\|_2(L_2^{j,h}+M_2^{j,h}).
        \end{aligned}
    \end{eqnarray}
\end{proposition}

\begin{proof}
    By \eqref{ineq2:propo4:final:diffe:sum}, $\nabla p_{f,g}$ is the sum of the maps $\alpha_i\times$\eqref{ineq1:propo5:final:diffe:Lipschitz}, $-\alpha_i\times$\eqref{ineq2:propo5:final:diffe:Lipschitz}, $\beta_h^jA_j^T\times$\eqref{ineq3:propo5:final:diffe:Lipschitz} and $-\beta_h^jA_j^T\times$\eqref{ineq4:propo5:final:diffe:Lipschitz}; since $\|A_j^Tz\|\leq\|A_j\|_2\|z\|$, the triangle inequality yields the constant \eqref{ineq6:propo5:final:diffe:Lipschitz}.
\end{proof}

\subsubsection{The simultaneous algorithm and its convergence}

\begin{breakablealgorithm}\caption{Generalized Simultaneous Projection Gradient Algorithm (GSPGA)}
\label{Algo:GSPGA} \begin{algorithmic}[1]\vspace{-0.2cm}
\State \textbf{Input}: $x^{0}\in\mathbb{R}^{n}$.
\State Choose $\lambda\in(0,2/L_{\nabla p_{f,g}})$.
\For{$k=0,1,2,\dots$}
\State Compute
\begin{eqnarray}\label{grad:breg:gen:algo}
\begin{aligned}
x^{k+1}:=x^k-\lambda \nabla p_{f,g}(x^k).
\end{aligned}
\end{eqnarray}
\EndFor
\end{algorithmic}
\end{breakablealgorithm}

\begin{remark}\label{rmk1:Breg:AfterAlgo2}
By Lemma~\ref{lem:proximity}, the problem \eqref{MSMSICFP} is equivalent, whenever its solution set $\varTheta$ is nonempty, to the
unconstrained minimization problem
\begin{eqnarray}
\begin{aligned}\min_{x\in\mathbb{R}^{n}}p_{f,g}(x).\end{aligned}
\label{min:rmk1:Breg:AfterAlgo2}
\end{eqnarray}
in the sense that $\varTheta=\arg\min_{x\in\mathbb{R}^n}p_{f,g}(x)$ and $\min p_{f,g}=0$. If $p_{f,g}$ is convex, then $\varTheta$ is closed and convex.
\end{remark}

\begin{theorem}\label{MainThm1:Bregman} Assume that all the hypotheses of Proposition \ref{propo5:final:diffe:Lipschitz} hold. Then $\lim_{k\to\infty}\|\nabla p_{f,g}(x^k)\|=0$ and every accumulation point of the
sequence $\{x^{k}\}_{k=0}^{\infty}$ generated by Algorithm \ref{Algo:GSPGA}
is a stationary point of $p_{f,g}$.
\end{theorem}

\begin{proof}
By Proposition \ref{propo5:final:diffe:Lipschitz}, the function $p_{f,g}$ is continuously
differentiable and its gradient is $L_{\nabla p_{f,g}}$-Lipschitz continuous; moreover, $p_{f,g}\geq0$ is bounded below. Since $\lambda\in(0,2/L_{\nabla p_{f,g}})$, Algorithm \ref{Algo:GSPGA} is the gradient method with a constant step size applied to the (possibly nonconvex) function $p_{f,g}$. Therefore, by Theorem~\ref{thm:rate:nonconvex} below (see also \cite{ABEC2017}), $\sum_k\|\nabla p_{f,g}(x^k)\|^2<\infty$, and hence
\begin{eqnarray}\label{ineq2:MainThm1:Bregman}
    \lim_{k\rightarrow\infty}\|\nabla p_{f,g}(x^k)\|=0.
\end{eqnarray}
Let $\{x^{k_\ell}\}$ be a subsequence of $\{x^k\}$ such that $x^{k_\ell}\rightarrow x^*$ as $\ell\rightarrow\infty$. Since $\nabla p_{f,g}$ is continuous, \eqref{ineq2:MainThm1:Bregman} yields $\nabla p_{f,g}(x^*)=0$.
\end{proof}

\begin{remark}
    Theorem~\ref{MainThm1:Bregman} only guarantees that the limit of any convergent subsequence is a stationary point of $p_{f,g}$. Since $p_{f,g}$ is not convex in general, a stationary point need not be a global minimizer, i.e., it need not solve \eqref{MSMSICFP}.
\end{remark}

\begin{theorem}\label{MainThm:2} Assume that all the hypotheses of Proposition \ref{propo5:final:diffe:Lipschitz} hold. Additionally, let $p_{f,g}$ be convex. Suppose the solution set $\varTheta$
of \eqref{MSMSICFP} is nonempty. Then any
sequence $\{x^{k}\}_{k=0}^{\infty}$ generated by Algorithm \ref{Algo:GSPGA}
converges in norm to a point in $\varTheta$.
\end{theorem}

\begin{proof}
    By Remark~\ref{rmk1:Breg:AfterAlgo2}, $\varTheta=\arg\min p_{f,g}$, so $\nabla p_{f,g}(\bar x)=0$ for every $\bar x\in\varTheta$. By Lemma~\ref{Baillon-Haddad}, $\nabla p_{f,g}$ is $1/L_{\nabla p_{f,g}}$-cocoercive, and inequality \eqref{eq:fejer:GSPGA} in the proof of Theorem~\ref{thm:rate:convex} below shows that $\{x^k\}$ is Fej\'er monotone with respect to $\varTheta$; in particular, it is bounded. By Theorem~\ref{MainThm1:Bregman}, every accumulation point $\bar x$ of $\{x^k\}$ satisfies $\nabla p_{f,g}(\bar x)=0$, hence $\bar x\in\varTheta$ by convexity. Since $\{\|x^k-\bar x\|\}$ is nonincreasing and a subsequence of it tends to zero, the whole sequence converges to $\bar x$.
\end{proof}

\subsubsection{Convergence rates of Algorithm~\ref{Algo:GSPGA}}\label{subsec:rates}

The estimates of this subsection are classical for the gradient method applied to a function with Lipschitz continuous gradient (see, e.g., \cite{ABEC2017,NEST2018,NNG2019}); we include short proofs adapted to $p_{f,g}$ for completeness and in order to make the constants explicit.
Throughout this subsection we write $L:=L_{\nabla p_{f,g}}$, assume $\lambda\in(0,2/L)$, and set
\begin{equation}\label{eq:clambda}
c_\lambda:=\lambda\Big(1-\frac{\lambda L}{2}\Big)>0\quad\text{and}\quad d_\lambda:=\lambda\Big(\frac{2}{L}-\lambda\Big)>0.
\end{equation}
Both constants are maximized at $\lambda=1/L$, where $c_{1/L}=1/(2L)$ and $d_{1/L}=1/L^2$.

\begin{theorem}[Nonconvex case]\label{thm:rate:nonconvex}
Assume that all the hypotheses of Proposition \ref{propo5:final:diffe:Lipschitz} hold and let $\{x^k\}$ be generated by Algorithm~\ref{Algo:GSPGA}. Then $\{p_{f,g}(x^k)\}$ is nonincreasing, $\sum_{k=0}^\infty\|\nabla p_{f,g}(x^k)\|^2\leq p_{f,g}(x^0)/c_\lambda$, and for every $N\geq0$,
\begin{equation}\label{eq:rate:nonconvex}
\min_{0\leq k\leq N}\|\nabla p_{f,g}(x^k)\|\leq\sqrt{\frac{p_{f,g}(x^0)}{c_\lambda(N+1)}}.
\end{equation}
\end{theorem}

\begin{proof}
By the descent lemma and \eqref{grad:breg:gen:algo},
\begin{equation}\label{eq:descent:GSPGA}
p_{f,g}(x^{k+1})\leq p_{f,g}(x^k)-\lambda\|\nabla p_{f,g}(x^k)\|^2+\frac{L\lambda^2}{2}\|\nabla p_{f,g}(x^k)\|^2=p_{f,g}(x^k)-c_\lambda\|\nabla p_{f,g}(x^k)\|^2.
\end{equation}
Summing \eqref{eq:descent:GSPGA} from $k=0$ to $N$ and using $p_{f,g}\geq0$, we get $c_\lambda\sum_{k=0}^N\|\nabla p_{f,g}(x^k)\|^2\leq p_{f,g}(x^0)-p_{f,g}(x^{N+1})\leq p_{f,g}(x^0)$. Bounding the sum from below by $(N+1)\min_{0\leq k\leq N}\|\nabla p_{f,g}(x^k)\|^2$ yields \eqref{eq:rate:nonconvex}.
\end{proof}

\begin{theorem}[Convex case]\label{thm:rate:convex}
Assume the hypotheses of Theorem~\ref{MainThm:2} and let $\{x^k\}$ be generated by Algorithm~\ref{Algo:GSPGA}. Set $R_0:=\mathrm{dist}(x^0,\varTheta)$. Then the following assertions hold.
\begin{itemize}
\item[(a)] For every $k\geq1$,
\begin{equation}\label{eq:rate:convex:a}
p_{f,g}(x^k)\leq\frac{R_0^2}{c_\lambda\,k}.
\end{equation}
\item[(b)] The sequence $\{\|\nabla p_{f,g}(x^k)\|\}$ is nonincreasing and, for every $k\geq0$,
\begin{equation}\label{eq:rate:convex:b}
\|\nabla p_{f,g}(x^k)\|\leq\frac{R_0}{\sqrt{d_\lambda(k+1)}}.
\end{equation}
\item[(c)] If, in addition, $p_{f,g}$ satisfies the quadratic growth condition
\begin{equation}\label{eq:QG}
p_{f,g}(x)\geq\frac{\mu}{2}\,\mathrm{dist}^2(x,\varTheta)\;\text{for all }x\text{ with }\mathrm{dist}(x,\varTheta)\leq R_0,
\end{equation}
for some $\mu>0$, then $\mu\leq L$ and, with $\rho:=\sqrt{1-d_\lambda\mu^2/4}\in (0,1)$,
\begin{equation}\label{eq:rate:convex:c}
\mathrm{dist}(x^{k},\varTheta)\leq\rho^k R_0\;\text{and}\;\|x^k-x^*\|\leq2\rho^kR_0\;\text{for all }k\geq0,
\end{equation}
where $x^*\in\varTheta$ is the limit of $\{x^k\}$ given by Theorem~\ref{MainThm:2}.
\end{itemize}
\end{theorem}

\begin{proof}
Denote $g_k:=\nabla p_{f,g}(x^k)$. Let $\bar x\in\varTheta$. Since $\nabla p_{f,g}(\bar x)=0$ and $\nabla p_{f,g}$ is $1/L$-cocoercive (Lemma~\ref{Baillon-Haddad}), we have $\langle g_k,x^k-\bar x\rangle\geq\frac1L\|g_k\|^2$, and therefore
\begin{equation}\label{eq:fejer:GSPGA}
\|x^{k+1}-\bar x\|^2=\|x^k-\bar x\|^2-2\lambda\langle g_k,x^k-\bar x\rangle+\lambda^2\|g_k\|^2\leq\|x^k-\bar x\|^2-d_\lambda\|g_k\|^2.
\end{equation}

(a) Since $p_{f,g}$ is convex and continuous, $\varTheta=\arg\min p_{f,g}$ is nonempty, closed, and convex. Take $\bar x:=P_\varTheta(x^0)$. By \eqref{eq:fejer:GSPGA}, $\|x^k-\bar x\|\leq R_0$ for all $k$. By convexity and the Cauchy--Schwarz inequality,
$\delta_k:=p_{f,g}(x^k)=p_{f,g}(x^k)-p_{f,g}(\bar x)\leq\langle g_k,x^k-\bar x\rangle\leq R_0\|g_k\|$.
Combined with \eqref{eq:descent:GSPGA}, this gives $\delta_{k+1}\leq\delta_k-\frac{c_\lambda}{R_0^2}\delta_k^2$. If $\delta_m=0$ for some $m$, then $x^m\in\varTheta$, $g_m=0$, and $x^k=x^m$ for all $k\geq m$, so \eqref{eq:rate:convex:a} is trivial for $k\ge m$; we may therefore assume $\delta_k>0$ for the indices under consideration (and $R_0>0$). Since $\delta_{k+1}\leq\delta_k$,
\begin{eqnarray}
\frac1{\delta_{k+1}}-\frac1{\delta_k}=\frac{\delta_k-\delta_{k+1}}{\delta_k\delta_{k+1}}\geq\frac{c_\lambda}{R_0^2}\,\frac{\delta_k^2}{\delta_k\delta_{k+1}}\geq\frac{c_\lambda}{R_0^2}.
\end{eqnarray}
Summing from $0$ to $k-1$ gives $1/\delta_k\geq c_\lambda k/R_0^2$, which is \eqref{eq:rate:convex:a}.

(b) Let $\Delta_k:=g_{k+1}-g_k$. By cocoercivity, $\langle\Delta_k,x^{k+1}-x^k\rangle\geq\frac1L\|\Delta_k\|^2$, i.e., $\langle\Delta_k,g_k\rangle\leq-\frac{1}{\lambda L}\|\Delta_k\|^2$. Hence
$\|g_{k+1}\|^2=\|g_k\|^2+2\langle\Delta_k,g_k\rangle+\|\Delta_k\|^2\leq\|g_k\|^2-\big(\tfrac{2}{\lambda L}-1\big)\|\Delta_k\|^2\leq\|g_k\|^2$, as $\lambda L<2$. Summing \eqref{eq:fejer:GSPGA} (with $\bar x=P_\varTheta(x^0)$) from $0$ to $k$ and using the monotonicity established above, $d_\lambda(k+1)\|g_k\|^2\leq d_\lambda\sum_{i=0}^k\|g_i\|^2\leq R_0^2$.

(c) Let $\bar x_k:=P_\varTheta(x^k)$ and $r_k:=\mathrm{dist}(x^k,\varTheta)$. Applying \eqref{eq:fejer:GSPGA} with $\bar x=\bar x_k$,
$r_{k+1}^2\leq\|x^{k+1}-\bar x_k\|^2\leq r_k^2-d_\lambda\|g_k\|^2$; in particular, $r_k\leq R_0$ for all $k$. By \eqref{eq:QG}, convexity, and the Cauchy--Schwarz inequality,
$\frac{\mu}{2}r_k^2\leq p_{f,g}(x^k)-p_{f,g}(\bar x_k)\leq\langle g_k,x^k-\bar x_k\rangle\leq\|g_k\|\,r_k$, so that $\|g_k\|\geq\frac{\mu}{2}r_k$. Consequently, $r_{k+1}^2\leq(1-d_\lambda\mu^2/4)r_k^2$. If $R_0=0$, then $x^k=x^0\in\varTheta$ for all $k$ and there is nothing to prove; let $R_0>0$. The descent lemma at $\bar x_0$, together with $\nabla p_{f,g}(\bar x_0)=0$, gives $p_{f,g}(x^0)\leq\frac L2R_0^2$, so \eqref{eq:QG} at $x^0$ yields $\mu\leq L$. Since $d_\lambda\leq1/L^2$, we get $1-d_\lambda\mu^2/4\geq3/4>0$, so $\rho$ is well defined and $r_k\leq\rho^kR_0$ by induction. Finally, by Fej\'er monotonicity, $\|x^m-\bar x_k\|\leq\|x^k-\bar x_k\|=r_k$ for all $m\geq k$; letting $m\to\infty$ gives $\|x^*-\bar x_k\|\leq r_k$, whence $\|x^k-x^*\|\leq\|x^k-\bar x_k\|+\|\bar x_k-x^*\|\leq2r_k\leq2\rho^kR_0$.
\end{proof}

\begin{remark}\label{rem:rates}
\begin{itemize}
\item[(i)] For $\lambda=1/L$, \eqref{eq:rate:convex:a} reads $p_{f,g}(x^k)\leq 2L\,\mathrm{dist}^2(x^0,\varTheta)/k$ and the linear factor in \eqref{eq:rate:convex:c} becomes $\rho=\sqrt{1-\mu^2/(4L^2)}$. The estimates hold for every $\lambda\in(0,2/L)$, i.e., for the whole step-size range of Algorithm~\ref{Algo:GSPGA}; for $\lambda\in(0,1/L]$ the sharper bound $p_{f,g}(x^k)\leq R_0^2/(2\lambda k)$ also follows from the argument used in the proof of Theorem~\ref{MainThm1:Algo:GPSPGA} below (with $\mathcal K=\mathbb{R}^n$).
\item[(ii)] Condition \eqref{eq:QG} is a quadratic growth (error bound) condition, closely related to the (bounded) linear regularity of the collection of constraint sets, which is classical in the analysis of projection methods for convex feasibility problems \cite{HHBA1996}; see \cite{NNG2019} for the role of such conditions in the linear convergence of first order methods. Identifying verifiable sufficient conditions for \eqref{eq:QG} in terms of the data of \eqref{MSMSICFP} is an interesting topic for future research. In the experiments of Section~\ref{Sec:numerics} we observe the linear behaviour predicted by \eqref{eq:rate:convex:c}.
\item[(iii)] The estimate of $L$ in Proposition~\ref{propo5:final:diffe:Lipschitz} is a sum of the constants of all components and may be conservative (see Section~\ref{Sec:numerics}); sharper constants directly improve all the bounds above.
\end{itemize}
\end{remark}

\subsubsection{The cyclic projection gradient algorithm and its convergence}

\begin{breakablealgorithm} \caption{Cyclic Projection Gradient Algorithm (CPGA)}
\label{Algo:CPGA}
\begin{algorithmic}[1] \vspace{-0.2cm}
\State \textbf{Input}: $x^{0}\in\mathbb{R}^{n}$ and a $q$-steering sequence $\{\sigma_{k}\}_{k\geq 0}$. (The functions $\varphi_i$ and $\psi_h^j$ are defined in \eqref{eq:components} below.)
\For{$k=0,1,\dots$}
\State $i(k)=(k\;\text{mod}\; r)+1$
\State $j(k)=(k\;\text{mod}\; t)+1$
\State $h(k)=\left(\left\lfloor \frac{k}{t}\right\rfloor
\bmod \ell_{j(k)}\right)+1$
\State Compute
\begin{eqnarray}\label{inq1:Algo:CPGA}
\begin{aligned}
 x^{k+1}& :=x^{k}-\sigma_k\Big(\alpha_{i(k)}\nabla\varphi_{i(k)}(x^k)+\beta_{h(k)}^{j(k)}\nabla\psi_{h(k)}^{j(k)}(x^k)\Big),
\end{aligned}
\end{eqnarray}
where
\begin{eqnarray}
\begin{aligned}
    \nabla\varphi_{i(k)}(x^k)&:=J_{\tilde{g}_{i(k)}}(x^k)^T\Big(\nabla f(P_{C_{i(k)}(x^k)}^f(x^k))-\nabla f(x^k)\Big)-\nabla^2 f(x^k)\Big(P_{C_{i(k)}(x^k)}^f(x^k)-x^k\Big),\\
    \nabla\psi_{h(k)}^{j(k)}(x^k)&:=A_{j(k)}^TJ_{p_{h(k)}^{j(k)}}(A_{j(k)}x^k)^T\Big(\nabla g_{j(k)}\Big(P_{Q_{h(k)}^{j(k)}(A_{j(k)}x^k)}^{g_{j(k)}}(A_{j(k)}x^k)\Big)-\nabla g_{j(k)}(A_{j(k)}x^k)\Big)\\&\hspace{1cm}-A_{j(k)}^T\nabla^2 g_{j(k)}(A_{j(k)}x^k)\Big(P_{Q_{h(k)}^{j(k)}(A_{j(k)}x^k)}^{g_{j(k)}}(A_{j(k)}x^k)-A_{j(k)}x^k\Big).
\end{aligned}
\end{eqnarray}
\EndFor
\end{algorithmic}
\end{breakablealgorithm}

The control sequences in Algorithm~\ref{Algo:CPGA} are periodic: with $q:=\operatorname{lcm}\{r,\,t\operatorname{lcm}\{\ell_1,\dots,\ell_t\}\}$ we have $i(k+q)=i(k)$, $j(k+q)=j(k)$, and $h(k+q)=h(k)$ for all $k\geq0$. Within one period, every index $i$ is visited $N_i:=q/r$ times and every pair $(j,h)$ is visited $N_{j,h}:=q/(t\ell_j)$ times. We set
\begin{equation}\label{eq:components}
\varphi_i(x):=D_f\big(P^f_{C_i(x)}(x),x\big),\qquad \psi_h^j(x):=D_{g_j}\big(P^{g_j}_{Q_h^j(A_jx)}(A_jx),A_jx\big),
\end{equation}
$F_{i,j,h}:=\alpha_i\varphi_i+\beta_h^j\psi_h^j$, $F_s:=F_{i(s),j(s),h(s)}$ for $s=0,1,\dots,q-1$, and introduce the \emph{period-averaged proximity function}
\begin{equation}\label{eq:phat}
\widehat{p}(x):=\sum_{s=0}^{q-1}F_s(x)=\sum_{i=1}^r N_i\alpha_i\varphi_i(x)+\sum_{j=1}^t\sum_{h=1}^{\ell_j}N_{j,h}\beta_h^j\psi_h^j(x),\quad x\in\mathbb{R}^n.
\end{equation}
Since all the multiplicities are positive, $\widehat p$ is the proximity function $p_{f,g}$ with the weights $N_i\alpha_i$ and $N_{j,h}\beta_h^j$; hence $\widehat p\geq0$ and $\widehat p(x)=0$ if and only if $x\in\varTheta$. We emphasize that one period of Algorithm~\ref{Algo:CPGA} is an (approximate) gradient step for $\widehat p$ and \emph{not} for $p_{f,g}$, unless the multiplicities are all equal.

\begin{theorem}\label{MainThm1:Algo:CPGA}
Assume that all the hypotheses of Proposition \ref{propo5:final:diffe:Lipschitz} hold, that each function $\varphi_i$ and $\psi_h^j$ is convex, and that the solution set $\varTheta$ of \eqref{MSMSICFP} is nonempty. For every triple $(i,j,h)$, let
$L_{i,j,h}:=\alpha_i(L_1^i+M_1^i)+\beta_h^j\|A_j\|_2(L_2^{j,h}+M_2^{j,h})$ (a Lipschitz constant of $\nabla F_{i,j,h}$, see Proposition~\ref{propo5:final:diffe:Lipschitz}), and set $L:=\max_{i,j,h}L_{i,j,h}$. Let $\{\sigma_k\}_{k\geq0}$ be a $q$-steering sequence such that $0<\sigma_k<2/L$ for all $k\geq0$. Then any sequence $\{x^{k}\}_{k=0}^{\infty}$ generated by Algorithm \ref{Algo:CPGA} converges in norm to a point in $\varTheta$.
\end{theorem}

The proof of Theorem \ref{MainThm1:Algo:CPGA} follows the strategy of \cite{YCEN2004} and is divided into several steps.

\begin{lemma}\label{Lemma1:MainThm1:Algo:CPGA}
Under the hypotheses of Theorem~\ref{MainThm1:Algo:CPGA}, for every $x^*\in\varTheta$ and every $k\geq0$,
\begin{equation}\label{ineq5:Lemma1:MainThm1:Algo:CPGA}
\|x^{k+1}-x^*\|^2\leq\|x^k-x^*\|^2-\sigma_k\Big(\frac{2}{L}-\sigma_k\Big)\big\|\nabla F_{i(k),j(k),h(k)}(x^k)\big\|^2 .
\end{equation}
In particular, $\{x^k\}$ is Fej\'er monotone with respect to $\varTheta$ and hence bounded.
\end{lemma}

\begin{proof}
Fix $x^*\in\varTheta$ and write $F:=F_{i(k),j(k),h(k)}$. Since $F$ is convex and nonnegative, and $F(x^*)=0$, the point $x^*$ is a global minimizer of $F$, so that $\nabla F(x^*)=0$. Since $\nabla F$ is Lipschitz continuous with a constant not exceeding $L$, Lemma~\ref{Baillon-Haddad} yields
$\langle\nabla F(x^k),x^k-x^*\rangle=\langle\nabla F(x^k)-\nabla F(x^*),x^k-x^*\rangle\geq\frac1L\|\nabla F(x^k)\|^2.$
Since $x^{k+1}=x^k-\sigma_k\nabla F(x^k)$ by \eqref{inq1:Algo:CPGA}, we get
\begin{eqnarray}
\begin{aligned}
\|x^{k+1}-x^*\|^2&=\|x^k-x^*\|^2-2\sigma_k\langle\nabla F(x^k),x^k-x^*\rangle+\sigma_k^2\|\nabla F(x^k)\|^2\\
&\leq\|x^k-x^*\|^2-\sigma_k\Big(\frac2L-\sigma_k\Big)\|\nabla F(x^k)\|^2.
\end{aligned}
\end{eqnarray}
Since $\sigma_k<2/L$, the sequence $\{\|x^k-x^*\|\}$ is nonincreasing.
\end{proof}

In what follows, $\mathcal{A}:=\overline{\operatorname{conv}}\{x^k\mid k\geq 0\}$, which is compact by Lemma~\ref{Lemma1:MainThm1:Algo:CPGA}, $y^k:=x^{kq}$ for $k\geq0$, and $\widehat L:=\sum_{s=0}^{q-1}L_{i(s),j(s),h(s)}$ is a Lipschitz constant of $\nabla\widehat p$.

\begin{proposition}\label{pro1:cyclic:Algo:CPGA}
Assume the hypotheses of Theorem~\ref{MainThm1:Algo:CPGA}. Then, for every $\mu>0$, there exists $\delta_\mu>0$ such that
\begin{eqnarray}\label{ineq3:pro1:cyclic:Algo:CPGA}
x\in\mathcal{A},\quad \widehat p(x)\geq\mu\;\Rightarrow\;\|\nabla \widehat p(x)\|\geq\delta_\mu.
\end{eqnarray}
\end{proposition}

\begin{proof}
Suppose, to the contrary, that there exist $\mu>0$ and a sequence $\{\tilde{x}^k\}\subset \mathcal{A}$ such that $\widehat p(\tilde{x}^k)\geq\mu$ and $\nabla \widehat p(\tilde{x}^k)\to 0$ as $k\to\infty$. Since $\mathcal{A}$ is compact, a subsequence of $\{\tilde x^k\}$ converges to some $\tilde{x}\in \mathcal{A}$. By continuity, $\widehat p(\tilde{x})\geq\mu$ and $\nabla \widehat p(\tilde{x})=0$. Since $\widehat p$ is convex and differentiable, $\tilde{x}$ is a global minimizer of $\widehat p$, and therefore $\widehat p(\tilde{x})=\min\widehat p=0$, because $\varTheta\neq\varnothing$. This contradicts $\widehat p(\tilde{x})\geq\mu>0$.
\end{proof}

\begin{proposition}\label{pro2:cyclic:Algo:CPGA}
Assume the hypotheses of Theorem~\ref{MainThm1:Algo:CPGA}. Then $\max_{0\leq s\leq q}\|x^{kq+s}-y^k\|\to0$ as $k\to\infty$, and there exists a sequence of vectors $\{r^k\}_{k\geq 0}$ with $\lim_{k\to\infty} r^k=0$ such that
\begin{eqnarray}\label{ineq0:pro2:cyclic:Algo:CPGA}
    y^k-y^{k+1}=\sigma_{kq}\big(\nabla \widehat p(y^k)+r^k\big)\quad\text{for all }k\geq0.
\end{eqnarray}
\end{proposition}

\begin{proof}
By periodicity, the step $x^{kq+s}\mapsto x^{kq+s+1}$ uses the function $F_s$. Since $\{x^k\}\subset\mathcal{A}$ and $\nabla F_s$ is continuous, there exists $M>0$ such that $\|\nabla F_s(x^k)\|\leq M$ for all $s$ and $k$. Hence, for $0\leq s\leq q$,
\begin{eqnarray}
\|x^{kq+s}-y^k\|\leq \sum_{u=0}^{s-1}\|x^{kq+u+1}-x^{kq+u}\|\leq M\sum_{u=0}^{q-1}\sigma_{kq+u}\longrightarrow0\quad\text{as }k\to\infty,
\end{eqnarray}
since $q$ is fixed and $\sigma_k\to0$. By the uniform continuity of $\nabla F_s$ on $\mathcal A$,
$\eta_k^s:=\nabla F_s(x^{kq+s})-\nabla F_s(y^k)\to0$, and, by the $q$-steering property,
$\xi_k^s:=\sigma_{kq+s}/\sigma_{kq}-1\to0$ as $k\to\infty$, for each $s=0,\dots,q-1$. Therefore,
\begin{eqnarray}
\begin{aligned}
y^k-y^{k+1}&=\sum_{s=0}^{q-1}\big(x^{kq+s}-x^{kq+s+1}\big)=\sum_{s=0}^{q-1}\sigma_{kq+s}\nabla F_s(x^{kq+s})\\
&=\sigma_{kq}\sum_{s=0}^{q-1}(1+\xi_k^s)\big(\nabla F_s(y^k)+\eta_k^s\big)=\sigma_{kq}\big(\nabla\widehat p(y^k)+r^k\big),
\end{aligned}
\end{eqnarray}
where
\begin{eqnarray}\label{ineq12:pro2:cyclic:Algo:CPGA}
r^k:=\sum_{s=0}^{q-1}\Big(\eta_k^s+\xi_k^s\big(\nabla F_s(y^k)+\eta_k^s\big)\Big)\longrightarrow0\quad\text{as }k\to\infty,
\end{eqnarray}
because $\{\nabla F_s(y^k)\}$ is bounded.
\end{proof}

\begin{proposition}\label{pro4:cyclic:Algo:CPGA}
Assume the hypotheses of Theorem~\ref{MainThm1:Algo:CPGA}. Then, for every $\mu>0$,
there exist $k_0\in\mathbb{N}$ and $\bar{\varepsilon}>0$ such that
\begin{eqnarray}\label{ineq1:pro4:cyclic:Algo:CPGA}
\widehat p(y^{k+1})\leq \widehat p(y^k)-\bar{\varepsilon}\sigma_{kq}\quad\text{for all }k\geq k_0\text{ with }\widehat p(y^k)\geq\mu.
\end{eqnarray}
\end{proposition}

\begin{proof}
Let $\delta_\mu$ be given by Proposition~\ref{pro1:cyclic:Algo:CPGA}. Since $r^k\to0$ and $\sigma_{kq}\to0$, there exists $k_0$ such that $\|r^k\|\leq\delta_\mu/4$ and $\widehat L\sigma_{kq}\leq 8/25$ for all $k\geq k_0$. Let $k\geq k_0$ with $\widehat p(y^k)\geq\mu$, and write $G:=\nabla\widehat p(y^k)$. Then $\|G\|\geq\delta_\mu$, so $\|r^k\|\leq\frac14\|G\|$, which implies
$\langle G,G+r^k\rangle\geq\frac34\|G\|^2$ and $\|G+r^k\|^2\leq\frac{25}{16}\|G\|^2$. By the descent lemma and \eqref{ineq0:pro2:cyclic:Algo:CPGA},
\begin{eqnarray}
\begin{aligned}
\widehat p(y^{k+1})&\leq\widehat p(y^k)-\sigma_{kq}\langle G,G+r^k\rangle+\frac{\widehat L}{2}\sigma_{kq}^2\|G+r^k\|^2\\
&\leq\widehat p(y^k)-\sigma_{kq}\|G\|^2\Big(\frac34-\frac{25}{32}\widehat L\sigma_{kq}\Big)\leq\widehat p(y^k)-\frac{\sigma_{kq}}{2}\|G\|^2\leq\widehat p(y^k)-\frac{\delta_\mu^2}{2}\sigma_{kq}.
\end{aligned}
\end{eqnarray}
Thus \eqref{ineq1:pro4:cyclic:Algo:CPGA} holds with $\bar\varepsilon:=\delta_\mu^2/2$.
\end{proof}

\noindent\textbf{Proof of Theorem \ref{MainThm1:Algo:CPGA}}.
\emph{Step 1: $\liminf_{k\to\infty}\widehat p(y^k)=0$.} Note first that $\sum_{k\geq0}\sigma_{kq}=+\infty$. Indeed, by the $q$-steering property, $\sum_{s=0}^{q-1}\sigma_{kq+s}\leq2q\,\sigma_{kq}$ for all sufficiently large $k$, while $\sum_{k}\sum_{s=0}^{q-1}\sigma_{kq+s}=\sum_k\sigma_k=+\infty$. Suppose now that $\liminf_{k\to\infty}\widehat p(y^k)>0$. Then there exist $\nu>0$ and $k_1$ such that $\widehat p(y^k)\geq\nu$ for all $k\geq k_1$. By Proposition~\ref{pro4:cyclic:Algo:CPGA} (with $\mu=\nu$), there are $k_0$ and $\bar\varepsilon>0$ such that, for $K:=\max\{k_0,k_1\}$ and every $m\geq K$,
\begin{eqnarray*}
0\leq\widehat p(y^{m+1})\leq\widehat p(y^{K})-\bar\varepsilon\sum_{k=K}^{m}\sigma_{kq}\longrightarrow-\infty\quad\text{as }m\to\infty,
\end{eqnarray*}
a contradiction.

\emph{Step 2: $\lim_{k\to\infty}\widehat p(y^k)=0$.} Suppose that $\limsup_{k\to\infty}\widehat p(y^k)>0$ and choose $0<\nu_1<\nu_2<\limsup_{k\to\infty}\widehat p(y^k)$. By Proposition~\ref{pro2:cyclic:Algo:CPGA}, $\|y^{k+1}-y^k\|\to0$, and since $\widehat p$ is uniformly continuous on $\mathcal{A}$, there exists $k_2$ such that $|\widehat p(y^{k+1})-\widehat p(y^k)|<\nu_2-\nu_1$ for all $k\geq k_2$. By Proposition~\ref{pro4:cyclic:Algo:CPGA} (with $\mu=\nu_1$), there exists $k_0$ such that $\widehat p(y^{k+1})<\widehat p(y^k)$ whenever $k\geq k_0$ and $\widehat p(y^k)\geq\nu_1$. By Step~1, there exists $k_3\geq\max\{k_0,k_2\}$ such that $\widehat p(y^{k_3})<\nu_1$. We claim that $\widehat p(y^k)<\nu_2$ for all $k\geq k_3$. This holds for $k=k_3$. If it holds for some $k\geq k_3$, then either $\widehat p(y^k)<\nu_1$, in which case $\widehat p(y^{k+1})<\widehat p(y^k)+(\nu_2-\nu_1)<\nu_2$, or $\nu_1\leq\widehat p(y^k)<\nu_2$, in which case $\widehat p(y^{k+1})<\widehat p(y^k)<\nu_2$. This proves the claim, which contradicts $\limsup_{k\to\infty}\widehat p(y^k)>\nu_2$.

\emph{Step 3: convergence of $\{x^k\}$.} By Proposition~\ref{pro2:cyclic:Algo:CPGA}, $\max_{0\leq s\leq q-1}\|x^{kq+s}-y^k\|\to0$, and by the uniform continuity of $\widehat p$ on $\mathcal{A}$ and Step~2, $\widehat p(x^{kq+s})\to0$ for every $s\in\{0,\dots,q-1\}$. Since these finitely many residue classes cover the whole sequence, $\lim_{k\to\infty}\widehat p(x^k)=0$. Let $\bar x$ be an accumulation point of the bounded sequence $\{x^k\}$, say $x^{k_\ell}\to\bar{x}$. By continuity, $\widehat p(\bar{x})=0$, i.e., $\bar{x}\in\varTheta$. By Lemma~\ref{Lemma1:MainThm1:Algo:CPGA}, the sequence $\{\|x^k-\bar x\|\}$ is nonincreasing; since its subsequence $\{\|x^{k_\ell}-\bar x\|\}$ tends to zero, $x^k\to\bar x$ as $k\to\infty$.
\hfill\rule{0.5em}{0.5em}

\begin{remark}
Because of the diminishing steering steps, Theorem~\ref{MainThm1:Algo:CPGA} does not provide a rate of convergence and, in our experiments (Section~\ref{Sec:numerics}), Algorithm~\ref{Algo:CPGA} is considerably slower than Algorithm~\ref{Algo:GSPGA} in terms of iterations. Its advantage is that each iterative step requires only two Bregman projections, which is attractive when the number of constraints is very large.
\end{remark}

\begin{remark}
    In Theorems~\ref{MainThm:2}, \ref{thm:rate:convex}, and \ref{MainThm1:Algo:CPGA}, we assume that $p_{f,g}$, respectively, each of its components, is convex. This assumption is restrictive in general. Lemma~\ref{lem:convex:components} below provides verifiable sufficient conditions in the Euclidean case, which cover both applications of Section~\ref{Sec:numerics}, and Example~\ref{example1:cyclic:convex:bregman} shows that the assumption can also hold for a non-quadratic Bregman function.
\end{remark}

\begin{lemma}\label{lem:convex:components}
Let $f=\frac12\|\cdot\|^2$, let $\varOmega\subseteq\mathbb{R}^n$ be nonempty, closed, and convex, let $\tilde g:\mathbb{R}^n\to\mathbb{R}^n$ be continuous, $C(x)=\tilde g(x)+\varOmega$, and $s(x):=x-\tilde g(x)$. Then
\begin{equation}\label{eq:phi:dist}
\varphi(x):=D_f\big(P_{C(x)}^f(x),x\big)=\frac12\,\mathrm{dist}^2\big(s(x),\varOmega\big),\qquad x\in\mathbb{R}^n,
\end{equation}
and the function $x\mapsto\mathrm{dist}(s(x),\varOmega)$, and hence $\varphi$, is convex in each of the following cases:
\begin{itemize}
\item[(a)] $s$ is affine;
\item[(b)] $\varOmega=\prod_{k=1}^n[a_k,b_k]$ with $-\infty\leq a_k\leq b_k\leq+\infty$ and, for each $k$, the $k$-th component $s_k$ of $s$ is affine, or convex with $s_k(x)\geq a_k$ for all $x$, or concave with $s_k(x)\leq b_k$ for all $x$.
\end{itemize}
The same conclusions hold for $\psi(x):=D_{g}\big(P^{g}_{Q(Ax)}(Ax),Ax\big)=\frac12\mathrm{dist}^2\big(Ax-q(Ax),\tilde\varOmega\big)$, where $g=\frac12\|\cdot\|^2$, $Q(y)=q(y)+\tilde\varOmega$, $A$ is linear, and $s(x):=Ax-q(Ax)$. In particular, in all these cases the sets $\{x\mid\mathrm{dist}(s(x),\varOmega)\leq R\}$, $R>0$, are convex.
\end{lemma}

\begin{proof}
By Remark~\ref{rmk1:tran:breg}, $P_{C(x)}(x)-x=P_\varOmega(s(x))-s(x)$, which gives \eqref{eq:phi:dist}. (a) The function $\mathrm{dist}(\cdot,\varOmega)$ is convex, and its composition with an affine map is convex. (b) We have $\mathrm{dist}(s(x),\varOmega)=\big(\sum_{k=1}^n d_k(x)^2\big)^{1/2}$ with $d_k(x):=\mathrm{dist}(s_k(x),[a_k,b_k])=[s_k(x)-b_k]_++[a_k-s_k(x)]_+$ (terms with infinite bounds being omitted). If $s_k$ is affine, $d_k$ is convex. If $s_k$ is convex and $s_k\geq a_k$, then $d_k=[s_k-b_k]_+$ is the composition of the convex nondecreasing function $t\mapsto[t-b_k]_+$ with a convex function, hence convex; the concave case is symmetric. Since the Euclidean norm is convex and nondecreasing in each coordinate on $\mathbb{R}^n_+$, the function $x\mapsto\|(d_1(x),\dots,d_n(x))\|$ is convex. Finally, $\varphi$ is the composition of this nonnegative convex function with the convex function $t\mapsto t^2/2$, which is nondecreasing on $[0,\infty)$, and is therefore convex. The proof for $\psi$ is identical, and the last assertion holds because sublevel sets of convex functions are convex.
\end{proof}

\begin{example}\label{example1:cyclic:convex:bregman}
    Consider $\Omega=[0,1]^n$, and for $x\in\mathbb{R}^n$
    \begin{eqnarray}
        f(x)&=& \sum_{i=1}^n\Big(\frac{17}{12}x_i^2+\frac{x_i}{2}\Big(\frac{x_i^2}{3}-1\Big) \arctan(x_i)+\Big(\frac{1}{12}-\frac{x_i^2}{4}\Big)\log (1+x_i^2)\Big),\\
        g(x)&=& (x_i+d_i)_{i=1}^n=x+d,
    \end{eqnarray}
    where $d\in\mathbb{R}^n_+.$
    Note that $\Omega$ is compact and let $C(x):=g(x)+\Omega$.
    Since $f$ is separable, $\varphi'$ is increasing (see below) and $x_i+d_i+\omega_i\geq x_i$ for $\omega\in\Omega$, the derivative $\varphi'(x_i+d_i+\omega_i)-\varphi'(x_i)$ of the $i$-th term of $\omega\mapsto D_f(g(x)+\omega,x)$ is nonnegative, so this term is nondecreasing in $\omega_i\geq0$; hence, by Proposition \ref{propo1:tran:breg}, $P_{C(x)}^f(x)=g(x)+\omega^*$ with $\omega^*=0$.
    It is easy to check that $J_{g}(x)=I$ and
    \begin{eqnarray}
        \nabla^2 f(x)=\mathrm{diag}(2+x_1 \arctan(x_1)-\frac{1}{2}\log(1+x_1^2),\dots,2+x_n \arctan(x_n)-\frac{1}{2}\log(1+x_n^2))
    \end{eqnarray}
    Let $u(x_i):=2+x_i \arctan(x_i)-\frac{1}{2}\log(1+x_i^2)$. Then $u^\prime(x_i)=\arctan(x_i)$, and therefore,
    $\vert \arctan(x_i)\vert\leq \pi/2$. Thus, $u$ is $\pi/2$-Lipschitz.
    Now, let
    \begin{eqnarray}
        h_1(x):=J_{g}(x)^T(\nabla f(P_{C(x)}^f(x))-\nabla f(x))=\nabla f(x+d)-\nabla f(x).
    \end{eqnarray}
    Its Jacobian is $\mathrm{diag}\big(u(x_i+d_i)-u(x_i)\big)_{i=1}^n$, and $|u(x_i+d_i)-u(x_i)|\leq\frac{\pi}{2}d_i$, so
    \begin{eqnarray}
        \|J_{h_1}(x)\|_2\leq \frac{\pi}{2}\|d\|_\infty,
    \end{eqnarray}
    and $h_1$ is $\frac{\pi}{2}\|d\|_\infty$-Lipschitz by the mean value inequality.
    Again, let
    \begin{eqnarray}
        h_2(x):=\nabla^2 f(x)(P_{C(x)}^f(x)-x)=\nabla^2 f(x)d.
    \end{eqnarray}
    Since $h_2(x)=(u(x_i)d_i)_{i=1}^n$ and $u$ is $\frac{\pi}{2}$-Lipschitz, we have
    \begin{eqnarray}
       \| h_2(x)-h_2(y)\|_2\leq \frac{\pi}{2}\|d\|_\infty\|x-y\|_2.
    \end{eqnarray}
    It is observed that $v(x_i):=x_i \arctan(x_i)-\frac{1}{2}\log(1+x_i^2)=v(-x_i)$. For $x_i\geq 0$, $v^\prime(x_i)\geq 0$. Thus, $v$ is increasing for $x_i\geq 0$. Therefore, $v(x_i)\geq v(0)=0$.
    Hence, $u(x_i)\geq 2$. This further implies that $\nabla^2f\succeq 2I$; therefore, $f$ is $2$-strongly convex on $\mathbb{R}^n$ and, being of class $\mathcal{C}^2$, $f\in\mathscr{B}(\mathbb{R}^n)$.
\begin{eqnarray}
    D_f\left(P_{C(x)}^f(x),x\right)=f(x+d)-f(x)-\left\langle \nabla f(x),d\right\rangle=
\sum_{i=1}^{n} c_{i}(x_i),
\end{eqnarray}
where $\varphi:\mathbb{R}\to\mathbb{R}$ denotes the common one-dimensional summand of $f$, i.e., $f(x)=\sum_{i=1}^n\varphi(x_i)$, and
\begin{equation*}
c_{i}(s):=D_{\varphi}(s+d_i,s)=\varphi(s+d_i)-\varphi(s)-\varphi'(s)d_i,\qquad s\in\mathbb{R}.
\end{equation*}
Differentiating twice, $c_{i}''(s)=\varphi''(s+d_i)-\varphi''(s)-\varphi'''(s)d_i$. Here $\varphi''(s)=u(s)=2+s\arctan(s)-\frac{1}{2}\log(1+s^2)$, $\varphi'''(s)=\arctan(s)$ and $\varphi^{(4)}(s)=1/(1+s^2)>0$. Hence $\varphi'''$ is increasing, that is, $\varphi''$ is convex, and the tangent-line inequality for the convex function $\varphi''$ gives $\varphi''(s+d_i)\geq\varphi''(s)+\varphi'''(s)d_i$. Therefore $c_i''(s)\geq0$ for all $s\in\mathbb{R}$, each $c_i$ is convex, and $D_f\big(P_{C(x)}^f(x),x\big)=\sum_{i=1}^nc_i(x_i)$ is convex on $\mathbb{R}^n$.
Note that, since the translation in this example is a constant shift, the inclusion $x\in C(x)$ holds if and only if $-d\in\varOmega$, i.e., $d=0$; the example therefore only illustrates the convexity of an individual component for a non-quadratic Bregman function, while nontrivial feasibility problems with convex components are provided by Lemma~\ref{lem:convex:components} and Section~\ref{Sec:numerics}.
\end{example}

\subsection{Multi-split multiple-sets convex feasibility problem}\label{subsec:MSMSCFP}

Let $\mathbb{R}^{n}$ and $\mathbb{R}^{m_{j}}$ for $j=1,2,\dots,t$
be Euclidean spaces. Let $\{C_{i}\}_{i=1}^r$
be a family of nonempty, closed, and convex subsets of $\mathbb{R}^{n}$.
Let $\{Q_{h}^{j}\}_{j=1, h=1}^{t, \ell_j}$ be families of nonempty, closed, and convex subsets of $\mathbb{R}^{m_{j}}$. Let $A_{j}:\mathbb{R}^{n}\rightarrow\mathbb{R}^{m_{j}}$
be a bounded linear operator for each $j=1,2,\dots,t$. With these
given data the {\em multi-split multiple-sets convex feasibility problem} (MSMS-CFP) is defined as follows: find a point $x^{*}\in\mathbb{R}^{n}$
such that
\begin{equation}\label{MSMSCFP}
x^{*}\in\bigcap\limits_{i=1}^{r}C_{i}\text{ and }A_{j}x^{*}\in\bigcap\limits_{h=1}^{\ell_{j}}Q_{h}^{j},\;\text{for each }j=1,2,\dots,t.\tag{MSMS-CFP}
\end{equation}

It is important to note that the problem \eqref{MSMSCFP} is equivalent to finding a point $x^*\in \mathbb{R}^n$ such that
\begin{eqnarray}\label{eqivalence:MSMSCFP}
    x^{*}\in\bigcap\limits_{i=1}^{r}C_{i}\text{ and }x^{*}\in\bigcap_{j=1}^t\bigcap\limits_{h=1}^{\ell_{j}}A_{j}^{-1}\left(Q_{h}^{j}\right),
\end{eqnarray}
where $A_{j}^{-1}\left(Q_{h}^{j}\right)=\left\{x\in\mathbb{R}^n\;\vert\; A_jx\in Q_{h}^{j}\right\}.$

Given Bregman functions $f_i\in\mathscr{B}(\mathbb{R}^n)$ and $g_{j,h}\in\mathscr{B}(\mathbb{R}^{m_j})$, Lemma~\ref{lem:proximity} shows that
\begin{eqnarray}\label{ineq1:thm1:MSMSCFP}
p^\circ_{f,g}(x):=\sum_{i=1}^{r}\alpha_{i}D_{f_i}(P_{C_{i}}^{f_i}(x),x)+\sum_{j=1}^{t}\sum_{h=1}^{\ell_{j}}\beta_{h}^{j}D_{g_{j,h}}\big(P_{Q_{h}^{j}}^{g_{j,h}}(A_{j}x),A_{j}x\big)
\end{eqnarray}
is nonnegative and vanishes exactly at the solutions of \eqref{MSMSCFP}.
Let
\begin{eqnarray}\label{ineq3:thm1:MSMSCFP}
    \mathcal{I}_Q=\{(j,h)\;\vert\; 1\leq j\leq t,\; 1\leq h\leq \ell_j\}
\end{eqnarray}
be partitioned into $K$ subsets.
Suppose $\{\mathcal{B}_s^Q\}_{s=1}^K$ be the partition.
That is, $\mathcal{B}_m^Q\cap\mathcal{B}_n^Q=\varnothing$ for $m\neq n$ and $\mathcal{I}_Q=\cup_{s=1}^K\mathcal{B}_s^Q$.
Similarly, let
\begin{eqnarray}\label{ineq4:thm1:MSMSCFP}
    \mathcal{I}=\{i\;\vert\; 1\leq i\leq r\}
\end{eqnarray}
be partitioned into $\{\mathcal{B}_s\}_{s=1}^K$.

\begin{definition}\label{def1:thm1:MSMSCFP}
    Let $u$ be a super-coercive Bregman function on $\mathbb{R}^n$. Then we say that the Bregman distance $D_u$ is \textit{dominating} with respect to the problem \eqref{MSMSCFP} if, for each $s=1,2,\dots,K$, there exist $\eta_s>0$ and $\delta_s>0$ such that
    \begin{eqnarray}\label{ineq1:def1:thm1:MSMSCFP}
        D_u(x,y)\geq \eta_s\sum_{i\in \mathcal{B}_s} \alpha_i D_{f_i}(x,y)+\delta_s\sum_{(j,h)\in \mathcal{B}_s^Q} \beta_h^jD_{g_{j,h}}(A_jx,A_jy),\quad\text{for all }x,y\in\mathbb{R}^n.
    \end{eqnarray}
    Since all the functions involved are differentiable, \eqref{ineq1:def1:thm1:MSMSCFP} holds if and only if the function $u-\eta_s\sum_{i\in\mathcal{B}_s}\alpha_if_i-\delta_s\sum_{(j,h)\in\mathcal{B}_s^Q}\beta_h^j\,g_{j,h}\circ A_j$ is convex, for every $s$. For instance, if all $f_i$ and $g_{j,h}$ are $\frac12\|\cdot\|^2$, then $u=\frac{c}{2}\|\cdot\|^2$ is dominating (with $\eta_s=\delta_s=1$) whenever $c\geq\max_s\big\|\sum_{i\in\mathcal{B}_s}\alpha_iI+\sum_{(j,h)\in\mathcal{B}_s^Q}\beta_h^jA_j^TA_j\big\|_2$, and Algorithm~\ref{Algo:GPBCPA} below then reduces to a block-iterative CQ-type method with step size $1/c$.
\end{definition}

\begin{breakablealgorithm}\caption{Generalized Projected Block Cyclic Projection Algorithm (GPBCPA)}
\label{Algo:GPBCPA} \begin{algorithmic}[1]\vspace{-0.2cm}
\State \textbf{Input}: $x^{0}\in\mathbb{R}^n$.
\For{$k=0,1,2,\dots$}
\State $s=k\; (\text{mod }K)+1$
\State Compute
\begin{eqnarray}\label{ineq1:Algo:GPBCPA}
\begin{aligned}
x^{k+1}:=\arg\min_{x\in\mathbb{R}^n} H_k(x),
\end{aligned}
\end{eqnarray}
where
\begin{eqnarray}\label{ineq2:grad:breg:gen:Algo:GPSPGA}
\begin{aligned}
H_k(x)
&= D_u\left(x,x^{k}\right)+\eta_s \sum_{i\in \mathcal{B}_s}\alpha_i D_{f_i}\left(x,P_{C_i}^{f_i}x^{k}\right)-\eta_s \sum_{i\in \mathcal{B}_s}\alpha_i D_{f_i}\left(x,x^{k}\right)\\&
+\delta_s \sum_{(j,h)\in \mathcal{B}_s^Q}\beta_h^j D_{g_{j,h}}\left(A_jx,P_{Q_h^j}^{g_{j,h}}A_jx^{k}\right)-\delta_s \sum_{(j,h)\in \mathcal{B}_s^Q}\beta_h^j D_{g_{j,h}}\left(A_jx,A_jx^{k}\right).
\end{aligned}
\end{eqnarray}
\EndFor
\end{algorithmic}
\end{breakablealgorithm}

\begin{lemma}\label{lem1:Algo:GPBCPA}
    Let $f\in\mathscr{B}(\mathbb{R}^n)$.
    \begin{itemize}
    \item[(a)] If $x\in\mathbb{R}^n$ and $D_f(x,y^k)\to 0$, then $y^k\to x$ as $k\to  \infty$.
    \item[(b)] If $\{a^k\}$ and $\{b^k\}$ are bounded sequences, then $D_f(a^k,b^k)\to0$ if and only if $\|a^k-b^k\|\to0$.
    \end{itemize}
\end{lemma}

\begin{proof}
    (a) Since $D_f(x,y^k)\to 0$, there exists $M>0$ such that $D_f(x,y^k)\leq M$ for all $k$, i.e., $y^k\in L_2^f(x,M)=\{y\in\mathbb{R}^n\;\vert\; D_f(x,y)\leq M\}$, which is bounded because $f$ is a Bregman function. Let $z$ be an accumulation point of $\{y^k\}$, say $y^{k_\ell}\to z$. By the continuity of $D_f(x,\cdot)$ (recall that the zone of $f$ is $\mathbb{R}^n$), $D_f(x,z)=\lim_\ell D_f(x,y^{k_\ell})=0$, hence $z=x$ by the strict convexity of $f$. Since the bounded sequence $\{y^k\}$ has $x$ as its only accumulation point, $y^k\to x$.

    (b) Suppose that $D_f(a^k,b^k)\to0$ but $\|a^{k}-b^{k}\|\geq\varepsilon>0$ along a subsequence. Passing to a further subsequence, $a^k\to a$ and $b^k\to b$ with $\|a-b\|\geq\varepsilon$, while $D_f(a,b)=0$ by continuity, so $a=b$, a contradiction. Conversely, if $\|a^k-b^k\|\to0$, then $D_f(a^k,b^k)\to0$ because $D_f$ is uniformly continuous on compact sets and vanishes on the diagonal.
\end{proof}

\begin{lemma}\label{lem2:Algo:GPBCPA}
    Let $f\in\mathscr{B}(\mathbb{R}^n)$ and $C$ be a nonempty, closed, and convex subset of $\mathbb{R}^n$. Then the Bregman projection defined in \eqref{Breg:projection} is continuous on $\mathbb{R}^n$.
\end{lemma}

\begin{proof}
    Let $\{x^k\}$ be a sequence such that $x^k\to x^*$ as $k\to\infty$.
Write $p^k:=P_C^f(x^k)$. We first show that $\{p^k\}$ is bounded.
Suppose, to the contrary, that it is unbounded. Since $\{x^k\}$ is bounded, there exists a subsequence such that $\|p^{k_\ell}-x^{k_\ell}\|\to\infty$ as $\ell\to\infty$.
Fix $R>0$ and choose $\ell_0$ such that $\|p^{k_\ell}-x^{k_\ell}\|>R$ for all $\ell\geq\ell_0$.
For $\ell\geq\ell_0$, define $\lambda_{k_\ell}:=R/\|p^{k_\ell}-x^{k_\ell}\|\in(0,1)$; then $\lambda_{k_\ell}\to 0$ as $\ell\to\infty$.
Now for $\ell\geq \ell_0$, let $z^{k_\ell}=(1-\lambda_{k_\ell})x^{k_\ell}+\lambda_{k_\ell} p^{k_\ell}$.
Then $\|z^{k_\ell}-x^{k_\ell}\|=\lambda_{k_\ell}\|p^{k_\ell}-x^{k_\ell}\|=R$.
Thus, $\{z^{k_\ell}\}$ is bounded, and there is a subsequence $\{z^{k_{\ell_t}}\}$ such that $z^{k_{\ell_t}}\to z$ as $t\to\infty$; since $x^{k_{\ell_t}}\to x^*$, we have $\|z-x^*\|=R>0$.
By the convexity of the Bregman distance in its first argument,
\begin{eqnarray}\label{ineq1:lem2:Algo:GPBCPA}
\begin{aligned}
D_f(z^{k_{\ell_t}},x^{k_{\ell_t}})&\leq(1-\lambda_{k_{\ell_t}})D_f(x^{k_{\ell_t}},x^{k_{\ell_t}})+\lambda_{k_{\ell_t}} D_f(p^{k_{\ell_t}},x^{k_{\ell_t}})\\
&=\lambda_{k_{\ell_t}} D_f(p^{k_{\ell_t}},x^{k_{\ell_t}}).
\end{aligned}
\end{eqnarray}
Since \(x^{k_{\ell_t}}\to x^*\), $D_f(c,x^{k_{\ell_t}})\to D_f(c,x^*)$.
Therefore, for a fixed \(c\in C\), there exists \(M>0\) such that $D_f(c,x^{k_{\ell_t}})\leq M$ for all $t$.
Consequently, using the fact that for every \(c\in C\), $D_f(p^{k_\ell},x^{k_\ell})\leq D_f(c,x^{k_\ell})$, we get $0\leq D_f(z^{k_{\ell_t}},x^{k_{\ell_t}})\leq\lambda_{k_{\ell_t}} D_f(c,x^{k_{\ell_t}})\leq \lambda_{k_{\ell_t}} M.$
Thus, $D_f(z^{k_{\ell_t}},x^{k_{\ell_t}})\to0$ as $t\to\infty$.
By the joint continuity of \(D_f\), along the subsequences, we obtain $D_f(z,x^*)=0$, which implies $z=x^*$, a contradiction.
Therefore, $\{p^k\}$ is bounded.
Now, let \(\{p^{k_\nu}\}\) be any convergent subsequence. Assume that $p^{k_{\nu}}\to p^*$; then \(p^*\in C\) since $C$ is closed.
Since $D_f(p^{k_{\nu}},x^{k_{\nu}})\leq D_f(c,x^{k_{\nu}})$, for all $c\in C$, letting \(\nu\to\infty\) and using the continuity of \(D_f\), we obtain $D_f(p^*,x^*)\leq D_f(c,x^*)$, for all $c\in C$.
Hence, $p^*=P_C^f(x^*)$.
Since every convergent subsequence of the bounded sequence $\{p^k\}$ has the same limit $P_C^f(x^*)$, we conclude that $p^k\to P_C^f(x^*)$ as $k\to\infty$.
\end{proof}

\begin{theorem}\label{thm1:Algo:GPBCPA}
    Suppose $\{f_i\}_{i=1}^r$ and $\{g_{j,h}\}_{j=1,h=1}^{t,\ell_j}$ are Bregman functions with zones $\mathbb{R}^n$ and $\mathbb{R}^{m_j}$, respectively, and let $u$ be a super-coercive Bregman function on $\mathbb{R}^n$ such that $D_u$ is dominating with respect to \eqref{MSMSCFP}.
    Assume that the solution set $\Theta^\circ$ of the problem \eqref{MSMSCFP} is nonempty. Then the sequence generated by Algorithm \ref{Algo:GPBCPA} converges to a point in $\Theta^\circ$.
\end{theorem}

\begin{proof}
For brevity, we write $P_i:=P_{C_i}^{f_i}$ and $P_{j,h}:=P_{Q_h^j}^{g_{j,h}}$, and $s=s(k)=(k\bmod K)+1$ for the block used at iteration $k$.

\emph{Step 1: the iteration is well defined.} For any Bregman function $\theta$ and points $a,b$, the difference $D_\theta(x,a)-D_\theta(x,b)$ is an affine function of $x$. Hence $H_k(x)=u(x)-\langle c_k,x\rangle+\tilde R_k$, where $\tilde R_k$ is a constant and
\begin{eqnarray}\label{ineq2:thm1:Algo:GPBCPA}
\begin{aligned}
c_k&=\nabla u(x^{k})+\eta_s\sum_{i\in \mathcal{B}_s}\alpha_i\Big(\nabla f_i\big(P_i(x^{k})\big)-\nabla f_i(x^{k})\Big)\\&
+\delta_s\sum_{(j,h)\in \mathcal{B}_s^Q}\beta_h^jA_j^T\Big(\nabla g_{j,h}\big(P_{j,h}(A_jx^{k})\big)-\nabla g_{j,h}(A_jx^{k})\Big).
\end{aligned}
\end{eqnarray}
Since $u$ is strictly convex and super-coercive, $H_k$ has a unique minimizer $x^{k+1}$, characterized by $\nabla u(x^{k+1})=c_k$, and for all $x\in\mathbb{R}^n$,
\begin{eqnarray}\label{ineq3:thm1:Algo:GPBCPA}
H_k(x)-H_k(x^{k+1})=u(x)-u(x^{k+1})-\langle c_k,x-x^{k+1}\rangle=D_u(x,x^{k+1}).
\end{eqnarray}
Moreover, since $D_u$ is dominating, the sum of the first, third, and fifth terms in \eqref{ineq2:grad:breg:gen:Algo:GPSPGA} is nonnegative, and therefore
\begin{eqnarray}\label{ineq3b:thm1:Algo:GPBCPA}
H_k(x)\geq\eta_s\sum_{i\in \mathcal{B}_s}\alpha_iD_{f_i}\big(x,P_i(x^{k})\big)+\delta_s\sum_{(j,h)\in \mathcal{B}_s^Q}\beta_h^jD_{g_{j,h}}\big(A_jx,P_{j,h}(A_jx^{k})\big)\geq0,\quad x\in\mathbb{R}^n.
\end{eqnarray}

\emph{Step 2: Fej\'er-type inequality.} Let $x^*\in\Theta^\circ$. Evaluating \eqref{ineq3:thm1:Algo:GPBCPA} at $x=x^*$ and using the definition of $H_k$, we get
\begin{eqnarray}\label{ineq4:thm1:Algo:GPBCPA}
\begin{aligned}
D_u(x^*,x^{k})-D_u(x^*,x^{k+1})
&= H_k(x^{k+1})+\eta_s\sum_{i\in \mathcal{B}_s}\alpha_i\Big(D_{f_i}(x^*,x^{k})-D_{f_i}\big(x^*,P_i(x^{k})\big)\Big)\\&
+\delta_s\sum_{(j,h)\in \mathcal{B}_s^Q}\beta_h^j\Big(D_{g_{j,h}}(A_jx^*,A_jx^{k})-D_{g_{j,h}}\big(A_jx^*,P_{j,h}(A_jx^{k})\big)\Big)\\
&\geq H_k(x^{k+1})+\eta_s\sum_{i\in \mathcal{B}_s}\alpha_iD_{f_i}\big(P_i(x^k),x^k\big)\\&+\delta_s\sum_{(j,h)\in \mathcal{B}_s^Q}\beta_h^jD_{g_{j,h}}\big(P_{j,h}(A_jx^k),A_jx^k\big),
\end{aligned}
\end{eqnarray}
where we used the generalized Pythagorean inequality \eqref{eq:three:point} with $x^*\in C_i$ and $A_jx^*\in Q_h^j$. All terms on the right-hand side are nonnegative. Hence $\{D_u(x^*,x^{k})\}$ is decreasing and convergent, $\{x^k\}$ is bounded (it lies in the bounded level set $\{z\mid D_u(x^*,z)\leq D_u(x^*,x^0)\}$), and, summing \eqref{ineq4:thm1:Algo:GPBCPA} over $k$, as $k\to\infty$,
\begin{equation}\label{ineq5:thm1:Algo:GPBCPA}
\begin{aligned}
&H_k(x^{k+1})\to0,\qquad D_{f_i}\big(P_i(x^k),x^k\big)\to0\;\;\big(i\in\mathcal{B}_{s(k)}\big),\\
&D_{g_{j,h}}\big(P_{j,h}(A_jx^k),A_jx^k\big)\to0\;\;\big((j,h)\in\mathcal{B}^Q_{s(k)}\big).
\end{aligned}
\end{equation}

\emph{Step 3: asymptotic regularity.} By Lemma~\ref{lem2:Algo:GPBCPA}, the sequences $\{P_i(x^k)\}$ and $\{P_{j,h}(A_jx^k)\}$ are bounded, so Lemma~\ref{lem1:Algo:GPBCPA}(b) and \eqref{ineq5:thm1:Algo:GPBCPA} give $\|P_i(x^k)-x^k\|\to0$ and $\|P_{j,h}(A_jx^k)-A_jx^k\|\to0$ for the active indices, and then, again by Lemma~\ref{lem1:Algo:GPBCPA}(b), $D_{f_i}(x^k,P_i(x^k))\to0$ and $D_{g_{j,h}}(A_jx^k,P_{j,h}(A_jx^k))\to0$. Consequently, $H_k(x^k)=\eta_s\sum_{i\in\mathcal{B}_s}\alpha_iD_{f_i}(x^k,P_i(x^k))+\delta_s\sum_{(j,h)\in\mathcal{B}_s^Q}\beta_h^jD_{g_{j,h}}(A_jx^k,P_{j,h}(A_jx^k))\to0$. By \eqref{ineq3:thm1:Algo:GPBCPA} with $x=x^k$, $D_u(x^k,x^{k+1})=H_k(x^k)-H_k(x^{k+1})\to0$, and Lemma~\ref{lem1:Algo:GPBCPA}(b) yields $\|x^{k+1}-x^k\|\to0$.

\emph{Step 4: accumulation points solve \eqref{MSMSCFP}.} Let $y^*$ be an accumulation point of $\{x^k\}$, say $x^{k_\ell}\to y^*$. By Step~3, $x^{k_\ell+m}\to y^*$ for every fixed $m\in\{0,1,\dots,K-1\}$. Given a block index $s'\in\{1,\dots,K\}$, choose, for each $\ell$, $m_\ell\in\{0,\dots,K-1\}$ such that $s(k_\ell+m_\ell)=s'$; then $x^{k_\ell+m_\ell}\to y^*$. By the continuity of the Bregman projections (Lemma~\ref{lem2:Algo:GPBCPA}) and Step~3, $P_i(y^*)=y^*$ for every $i\in\mathcal{B}_{s'}$ and $P_{j,h}(A_jy^*)=A_jy^*$ for every $(j,h)\in\mathcal{B}^Q_{s'}$, i.e., $y^*\in C_i$ and $A_jy^*\in Q_h^j$. Since $s'$ is arbitrary and the blocks cover all the indices, $y^*\in\Theta^\circ$.

\emph{Step 5: convergence of the whole sequence.} By Step~2 with $x^*=y^*$, the sequence $\{D_u(y^*,x^k)\}$ is nonincreasing, and $D_u(y^*,x^{k_\ell})\to D_u(y^*,y^*)=0$ by continuity. Therefore $D_u(y^*,x^k)\to0$, and Lemma~\ref{lem1:Algo:GPBCPA}(a) yields $x^k\to y^*\in\Theta^\circ$.
\end{proof}

\begin{remark}
    The analysis carried out above can also be performed with different numbers of blocks for the constraints in the different spaces. For simplicity, we assumed that both partitions consist of the same number $K$ of blocks. The computational mechanism of Algorithm \ref{Algo:GPBCPA} for problem \eqref{MSMSCFP} is depicted in Figure~\ref{fig:blockiterative}.
    \begin{figure}[htbp]
    \centering
    \includegraphics[width=0.75\textwidth]{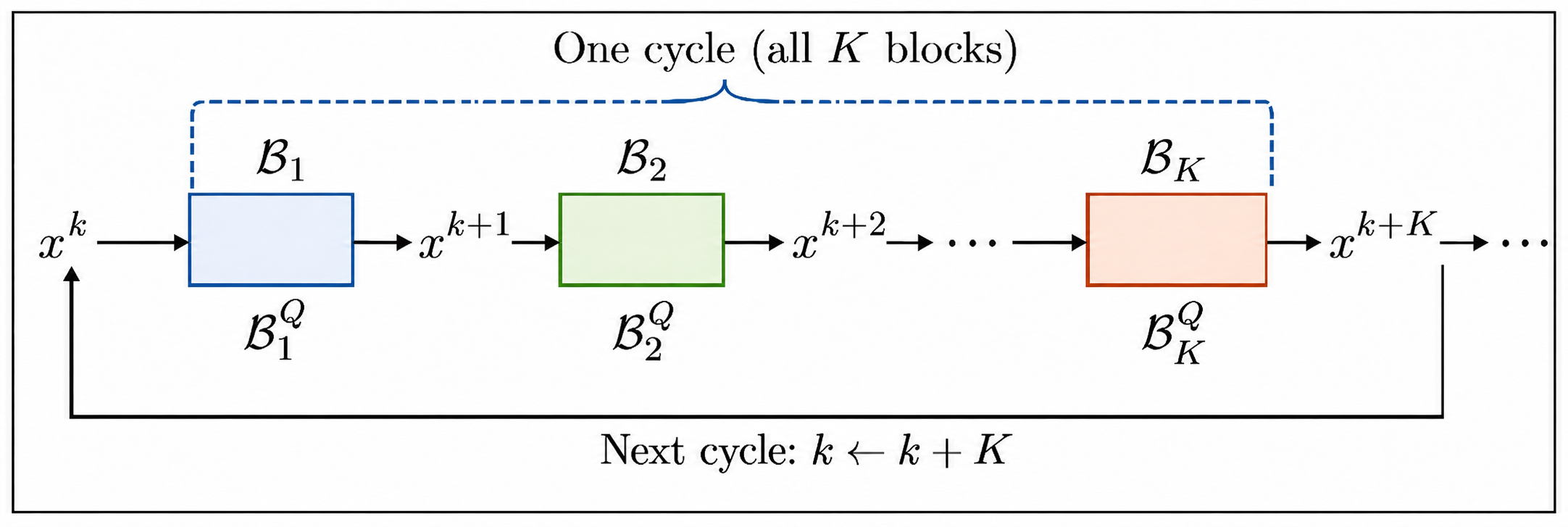}
    \caption{Computational mechanism of Algorithm \ref{Algo:GPBCPA}}
    \label{fig:blockiterative}
    \end{figure}
\end{remark}

\subsection{Local Lipschitz continuity and an algorithm}

In this subsection, the global Lipschitz continuity of the maps \eqref{ineq1:propo5:final:diffe:Lipschitz}--\eqref{ineq4:propo5:final:diffe:Lipschitz} in Proposition~\ref{propo5:final:diffe:Lipschitz} is replaced by bounds on the second derivatives of the Bregman functions and on the Jacobians of the translation operators. These conditions are imposed on the data rather than on composite maps, and they guarantee the Lipschitz continuity of $\nabla p_{f,g}$ on a residual-bounded set $\mathcal{K}$ (not necessarily on the whole space).

\begin{proposition}\label{propo5:final:diffe:Lipschitz:weak}
    Suppose that Assumption~\ref{ass:A} holds. Assume that there exist positive
    constants $\delta_{f}$ and $\Lambda_{f}$ such that
    \begin{eqnarray}\label{ineq1:propo5:final:diffe:Lipschitz:weak}
        \|\nabla^{2}f(x)\|_2\leq\delta_{f}\;\;\text{and}\;\;\|\nabla^{2}f(x)-\nabla^{2}f(y)\|_2\leq \Lambda_f\|x-y\|,\;\text{for all}\;x, y\in\mathbb{R}^n.
    \end{eqnarray}
    Suppose, for each $i=1,2,\dots,r$, there exist positive constants $\delta_{\tilde{g}_i}$ and $ \Lambda_{\tilde{g}_i}$ such that
    \begin{eqnarray}\label{ineq2:propo5:final:diffe:Lipschitz:weak}
        \|J_{\tilde{g}_i}(x)\|\leq \delta_{\tilde{g}_i}\;\;\text{and}\;\; \|J_{\tilde{g}_i}(x)-J_{\tilde{g}_i}(y)\|\leq \Lambda_{\tilde{g}_i}\|x-y\|,\; \text{for all}\; x,y\in\mathbb{R}^n.
    \end{eqnarray}
    Suppose, for each $j=1,2,\dots,t$ and $h=1,2,\dots,\ell_j$, there exist positive constants $\delta_{g_j}$, $\Lambda_{g_j}$, $\delta_{p_h^j}$ and $\Lambda_{p_h^j}$ such that, for all $y,y'\in\mathbb{R}^{m_j}$,
    \begin{eqnarray}\label{ineq3:propo5:final:diffe:Lipschitz:weak}
        \|\nabla^{2}g_j(y)\|_2\leq\delta_{g_j},\;\;\|\nabla^{2}g_j(y)-\nabla^{2}g_j(y')\|_2\leq \Lambda_{g_j}\|y-y'\|,
    \end{eqnarray}
    \begin{eqnarray}\label{ineq4:propo5:final:diffe:Lipschitz:weak}
            \|J_{p_h^j}(y)\|_2\leq \delta_{p_h^j},\;\;\|J_{p_h^j}(y)-J_{p_h^j}(y')\|_2\leq \Lambda_{p_h^j}\|y-y'\|.
    \end{eqnarray}
    Let $R_i>0$ and $K_i=\{x\in\mathbb{R}^n\vert \|x-P_{C_i(x)}^f(x)\|\leq R_i\}$. Additionally, let $T_h^j>0$ and $\tilde{K}_h^j=\{x\in\mathbb{R}^n\vert \|A_jx-P_{Q_h^j(A_jx)}^{g_j}(A_jx)\|\leq T_h^j\}$. Assume that $\mathcal{K}:=\Big(\cap_{i=1}^r K_i\Big)\cap\Big(\cap_{j=1}^{t}\cap_{h=1}^{\ell_j}\tilde{K}_h^j\Big)\neq\varnothing$.
    Then the function
    \begin{eqnarray}\label{ineq5:propo5:final:diffe:Lipschitz:weak}
        \mathcal{K}\ni x\mapsto  \nabla p_{f,g}(x)
    \end{eqnarray}
    is Lipschitz continuous on $\mathcal{K}$ with constant $L_{\mathcal{K}}$ given by
    \begin{eqnarray}\label{ineq6:propo5:final:diffe:Lipschitz:weak}
    \begin{aligned}
    &\hspace{-2cm}L_{\mathcal{K}}:=\sum_{i=1}^r\alpha_i\Big(\delta_{\tilde{g}_i}\delta_f\Big(1+\frac{\delta_f}{\sigma}\Big)(\delta_{\tilde{g}_i}+1)+\Lambda_{\tilde{g}_i} \delta_fR_i\Big)\\&
    \hspace{0.1cm}+\sum_{i=1}^{r}\alpha_{i}\Big(\Lambda_f R_i +\delta_f\Big(1+\Big(\delta_{\tilde{g}_i}+\frac{\delta_f}{\sigma}(\delta_{\tilde{g}_i}+1)\Big)\Big)\Big)\\&
    \hspace{0.1cm}+\sum_{j=1}^{t}\sum_{h=1}^{\ell_{j}}\beta_{h}^{j}\Big(\Lambda_{g_j} T_h^j +\delta_{g_j}\Big(1+\Big(\delta_{p_h^j}+\frac{\delta_{g_j}}{\sigma_j}(\delta_{p_h^j}+1)\Big)\Big)\Big)\|A_j\|_2^2\\&
    \hspace{0.1cm}+\sum_{j=1}^{t}\sum_{h=1}^{\ell_{j}}\beta_{h}^{j}\Big(\delta_{p_h^j}\delta_{g_j}\Big(1+\frac{\delta_{g_j}}{\sigma_j}\Big)(\delta_{p_h^j}+1)+\Lambda_{p_h^j}\delta_{g_j}T_h^j\Big)\|A_j\|_2^2.
    \end{aligned}
    \end{eqnarray}
\end{proposition}

\begin{proof}
For any $x,y\in\mathbb{R}^n$, using the fundamental theorem of integral calculus, we get
\begin{eqnarray}\label{ineq8:propo5:final:diffe:Lipschitz:weak}
    \nabla f(y)-\nabla f(x)=\int_0^1\nabla^2 f(x+t(y-x))(y-x)dt.
\end{eqnarray}
Since $\nabla^{2}f$ is bounded by $\delta_f$, \eqref{ineq8:propo5:final:diffe:Lipschitz:weak} implies that
\allowdisplaybreaks{
\begin{eqnarray}\label{ineq9:propo5:final:diffe:Lipschitz:weak}
    \|\nabla f(x)-\nabla f(y)\|&=&\Big\|\int_0^1\nabla^2 f(x+t(y-x))(x-y)dt\Big\|\nonumber\\
    &\leq&\int_0^1\|\nabla^2 f(x+t(y-x))\|_2\|x-y\|dt\\
    &\leq& \delta_f \|x-y\|.\nonumber
\end{eqnarray}}
Thus, $\nabla f$ is $\delta_{f}$-Lipschitz continuous.

Again,
\allowdisplaybreaks{
\begin{eqnarray}\label{ineq9.1:propo5:final:diffe:Lipschitz:weak}
    \|\tilde{g}_i(x)-\tilde{g}_i(y)\|&=&\Big\|\int_0^1J_{\tilde{g}_i}(x+t(y-x))(x-y)dt\Big\|\nonumber\\
    &\leq&\int_0^1\|J_{\tilde{g}_i}(x+t(y-x))\|_2\|x-y\|dt\\
    &\leq& \delta_{\tilde{g}_i}\|x-y\|.\nonumber
\end{eqnarray}}

For $x\in\mathcal{K}$,
\allowdisplaybreaks{
\begin{eqnarray}\label{ineq9.2:propo5:final:diffe:Lipschitz:weak}
    \|\nabla f(P_{C_i(x)}^f(x))-\nabla f(x)\|&\leq&\delta_f\|P_{C_i(x)}^f(x)-x\|\nonumber\\
    &\leq& \delta_{f}R_i.
\end{eqnarray}}

Now,
\begin{eqnarray}\label{ineq9.2.3:propo5:final:diffe:Lipschitz:weak}
    p_h^j(A_jy)-p_h^j(A_jx)=\int_0^1J_{p_h^j}(A_j(x+t(y-x)))A_j(y-x)dt .
\end{eqnarray}

Therefore,
\allowdisplaybreaks{
\begin{eqnarray}\label{ineq9.3:propo5:final:diffe:Lipschitz:weak}
\begin{aligned}
        \|p_h^j(A_jx)-p_h^j(A_jy)\|&=\Big\|\int_0^1J_{p_h^j}(A_j(x+t(y-x)))A_j(x-y)dt\Big\|\\
    &\leq\int_0^1\|J_{p_h^j}(A_j(x+t(y-x)))\|_2\|A_j\|_2\|x-y\|dt\\
    &\leq \delta_{p_h^j}\|A_j\|_2\|x-y\|.
\end{aligned}
\end{eqnarray}}

For any $x,y\in\mathcal{K}$, using our hypotheses, \eqref{ineq9.2:propo5:final:diffe:Lipschitz:weak}, and Proposition \ref{propo1:tran:breg:lipschitz} (applied with the constants $L_f:=\delta_f$ and $L_g:=\delta_{\tilde g_i}$, so that $x\mapsto P^f_{C_i(x)}(x)$ is Lipschitz continuous with constant $\delta_{\tilde{g}_i}+\frac{\delta_f}{\sigma}(\delta_{\tilde{g}_i}+1)$), we obtain
\allowdisplaybreaks{
\begin{eqnarray}\label{ineq10:propo5:final:diffe:Lipschitz:weak}
    &\|&J_{\tilde{g}_i}(x)^T(\nabla f(P_{C_i(x)}^f(x))-\nabla f(x))-J_{\tilde{g}_i}(y)^T(\nabla f(P_{C_i(y)}^f(y))-\nabla f(y))\|\nonumber\\
    &\leq& \|J_{\tilde{g}_i}(x)^T(\nabla f(P_{C_i(x)}^f(x))-\nabla f(P_{C_i(y)}^f(y)))\|
    +\|(J_{\tilde{g}_i}(x)-J_{\tilde{g}_i}(y))^T(\nabla f(P_{C_i(y)}^f(y))-\nabla f(y)\|\nonumber\\
    &&\hspace{1cm}+\|J_{\tilde{g}_i}(x)^T(\nabla f(x)-\nabla f(y))\|\\
    &\leq& \delta_{\tilde{g}_i}\delta_f\Big(\delta_{\tilde{g}_i}+\frac{\delta_f}{\sigma}(\delta_{\tilde{g}_i}+1)\Big)\|x-y\|+\Lambda_{\tilde{g}_i} \delta_f R_i\|x-y\|+\delta_{\tilde{g}_i}\delta_f\|x-y\|\nonumber\\
    &=&\Big(\delta_{\tilde{g}_i}\delta_f\Big(1+\frac{\delta_f}{\sigma}\Big)(\delta_{\tilde{g}_i}+1)+\Lambda_{\tilde{g}_i} \delta_fR_i\Big)\|x-y\|\nonumber
\end{eqnarray}}
and
\allowdisplaybreaks{
\begin{eqnarray}\label{ineq11:propo5:final:diffe:Lipschitz:weak}
\begin{aligned}
    &\|\nabla^2 f(x)(P_{C_i(x)}^f(x)-x)-\nabla^2 f(y)(P_{C_i(y)}^f(y)-y)\|\\
    &\leq \|(\nabla^2 f(x)-\nabla^2 f(y))(P_{C_i(x)}^f(x)-x)\|+\|\nabla^2 f(y)\|_2\|P_{C_i(x)}^f(x)-x-P_{C_i(y)}^f(y)+y\|\\
    &\leq \Lambda_f R_i \|x-y\|+\delta_f\Big(1+\Big(\delta_{\tilde{g}_i}+\frac{\delta_f}{\sigma}(\delta_{\tilde{g}_i}+1)\Big)\Big)\|x-y\|\\
    &=\Big(\Lambda_f R_i +\delta_f\Big(1+\Big(\delta_{\tilde{g}_i}+\frac{\delta_f}{\sigma}(\delta_{\tilde{g}_i}+1)\Big)\Big)\Big)\|x-y\|.
\end{aligned}    
\end{eqnarray}}

Similarly, for any $x,y\in\mathcal{K}$, using our hypotheses, \eqref{ineq9.3:propo5:final:diffe:Lipschitz:weak}, and Proposition \ref{propo2:tran:breg:lipschitz} (applied with $\tilde g:=g_j$, $q:=p_h^j$, $L_{\tilde g}:=\delta_{g_j}$ and $L_q:=\delta_{p_h^j}$), we derive
\allowdisplaybreaks{
\begin{eqnarray}\label{ineq12:propo5:final:diffe:Lipschitz:weak}
\begin{aligned}
    &\|J_{p_h^j}(A_jx)^T(\nabla g_j(P_{Q_h^j(A_jx)}^{g_j}(A_jx))-\nabla g_j(A_jx))-J_{p_h^j}(A_jy)^T(\nabla g_j(P_{Q_h^j(A_jy)}^{g_j}(A_jy))-\nabla g_j(A_jy))\|\\&
    \hspace{0.01cm}\leq\|J_{p_h^j}(A_jx)^T(\nabla g_j(P_{Q_h^j(A_jx)}^{g_j}(A_jx))-\nabla g_j(P_{Q_h^j(A_jy)}^{g_j}(A_jy)))\|\\&
    \hspace{0.5cm}+\|(J_{p_h^j}(A_jx)-J_{p_h^j}(A_jy))^T(\nabla g_j(P_{Q_h^j(A_jy)}^{g_j}(A_jy))-\nabla g_j(A_jy))\|\\&
    \hspace{0.5cm}+\|J_{p_h^j}(A_jx)^T(\nabla g_j(A_jx)-\nabla g_j(A_jy))\|\\&
    \hspace{0.01cm}\leq \delta_{p_h^j}\delta_{g_j}\Big(\delta_{p_h^j}+\frac{\delta_{g_j}}{\sigma_j}(\delta_{p_h^j}+1)\Big)\|A_j\|_2\|x-y\|+\Lambda_{p_h^j}\delta_{g_j}T_h^j\|A_j\|_2\|x-y\|\\&
    \hspace{0.5cm}+\delta_{p_h^j}\delta_{g_j}\|A_j\|_2\|x-y\|\\&
    \hspace{0.01cm}=\Big(\delta_{p_h^j}\delta_{g_j}\Big(1+\frac{\delta_{g_j}}{\sigma_j}\Big)(\delta_{p_h^j}+1)+\Lambda_{p_h^j}\delta_{g_j}T_h^j\Big)\|A_j\|_2\|x-y\|
\end{aligned}    
\end{eqnarray}}
and
\allowdisplaybreaks{
\begin{eqnarray}\label{ineq13:propo5:final:diffe:Lipschitz:weak}
\begin{aligned}
    &\hspace{-1cm}\|\nabla^2 g_j(A_jx)(P_{Q_h^j(A_jx)}^{g_j}(A_jx)-A_jx)-\nabla^2 g_j(A_jy)(P_{Q_h^j(A_jy)}^{g_j}(A_jy)-A_jy)\|\\
    &\leq \|(\nabla^2 g_j(A_jx)-\nabla^2 g_j(A_jy))(P_{Q_h^j(A_jx)}^{g_j}(A_jx)-A_jx)\|\\&
    \hspace{0.01cm}+\|\nabla^2 g_j(A_jy)\|\|P_{Q_h^j(A_jx)}^{g_j}(A_jx)-A_jx-P_{Q_h^j(A_jy)}^{g_j}(A_jy)+A_jy\|\\
    &\leq \Lambda_{g_j} T_h^j \|A_j\|_2\|x-y\|+\delta_{g_j}\Big(1+\Big(\delta_{p_h^j}+\frac{\delta_{g_j}}{\sigma_j}(\delta_{p_h^j}+1)\Big)\Big)\|A_j\|_2\|x-y\|\\
    &= \Big(\Lambda_{g_j} T_h^j +\delta_{g_j}\Big(1+\Big(\delta_{p_h^j}+\frac{\delta_{g_j}}{\sigma_j}(\delta_{p_h^j}+1)\Big)\Big)\Big)\|A_j\|_2\|x-y\|.
\end{aligned}
\end{eqnarray}}

Based on our hypotheses, by Proposition \ref{propo4:final:diffe:sum}, we get
\allowdisplaybreaks{
\begin{eqnarray}\label{ineq14:propo5:final:diffe:Lipschitz:weak}
    \begin{aligned}
        \nabla p_{f,g}(x)&=\sum_{i=1}^{r}\alpha_{i}\Big( J_{\tilde{g}_i}(x)^T(\nabla f(P_{C_i(x)}^f(x))-\nabla f(x))-\nabla^2 f(x)(P_{C_i(x)}^f(x)-x)\Big)\\&\hspace{0.5cm}+\sum_{j=1}^{t}\sum_{h=1}^{\ell_{j}}\beta_{h}^{j}\Big(A_j^TJ_{p_h^j}(A_jx)^T(\nabla g_j(P_{Q_h^j(A_jx)}^{g_j}(A_jx))\\&\hspace{1cm}-\nabla g_j(A_jx))-A_j^T\nabla^2 g_j(A_jx)(P_{Q_h^j(A_jx)}^{g_j}(A_jx)-A_jx)\Big).
    \end{aligned}
\end{eqnarray}}

Therefore, using \eqref{ineq10:propo5:final:diffe:Lipschitz:weak}-\eqref{ineq13:propo5:final:diffe:Lipschitz:weak}, we obtain
\allowdisplaybreaks{
\begin{eqnarray}\label{ineq15:propo5:final:diffe:Lipschitz:weak}
    &&\hspace{-1cm}\|\nabla p_{f,g}(x)-\nabla p_{f,g}(y)\|\nonumber\\
    &&\leq\sum_{i=1}^r\alpha_i\|J_{\tilde{g}_i}(x)^T(\nabla f(P_{C_i(x)}^f(x))-\nabla f(x))-J_{\tilde{g}_i}(y)^T(\nabla f(P_{C_i(y)}^f(y))-\nabla f(y))\|\nonumber\\&&
    \hspace{0.01cm}+\sum_{i=1}^{r}\alpha_{i}\|\nabla^2 f(x)(P_{C_i(x)}^f(x)-x)-\nabla^2 f(y)(P_{C_i(y)}^f(y)-y)\|\nonumber\\&&
    \hspace{0.01cm}+\sum_{j=1}^{t}\sum_{h=1}^{\ell_{j}}\beta_{h}^{j}\|A_j\|_2\|\nabla^2 g_j(A_jx)(P_{Q_h^j(A_jx)}^{g_j}(A_jx)-A_jx)\nonumber\\&&
    \hspace{3cm}-\nabla^2 g_j(A_jy)(P_{Q_h^j(A_jy)}^{g_j}(A_jy)-A_jy)\|\nonumber\\&&
    \hspace{0.01cm}+\sum_{j=1}^{t}\sum_{h=1}^{\ell_{j}}\beta_{h}^{j}\|A_j\|_2\|J_{p_h^j}(A_jx)^T(\nabla g_j(P_{Q_h^j(A_jx)}^{g_j}(A_jx))-\nabla g_j(A_jx))\\&&
    \hspace{3cm}-J_{p_h^j}(A_jy)^T(\nabla g_j(P_{Q_h^j(A_jy)}^{g_j}(A_jy))-\nabla g_j(A_jy))\|\nonumber\\
    &&\leq \bigg(\sum_{i=1}^r\alpha_i\Big(\delta_{\tilde{g}_i}\delta_f\Big(1+\frac{\delta_f}{\sigma}\Big)(\delta_{\tilde{g}_i}+1)+\Lambda_{\tilde{g}_i} \delta_fR_i\Big)\nonumber\\&&
    \hspace{0.2cm}+\sum_{i=1}^{r}\alpha_{i}\Big(\Lambda_f R_i +\delta_f\Big(1+\Big(\delta_{\tilde{g}_i}+\frac{\delta_f}{\sigma}(\delta_{\tilde{g}_i}+1)\Big)\Big)\Big)\nonumber\\&&
    \hspace{0.2cm}+\sum_{j=1}^{t}\sum_{h=1}^{\ell_{j}}\beta_{h}^{j}\Big(\Lambda_{g_j} T_h^j +\delta_{g_j}\Big(1+\Big(\delta_{p_h^j}+\frac{\delta_{g_j}}{\sigma_j}(\delta_{p_h^j}+1)\Big)\Big)\Big)\|A_j\|_2^2\nonumber\\&&
\hspace{0.2cm}+\sum_{j=1}^{t}\sum_{h=1}^{\ell_{j}}\beta_{h}^{j}\Big(\delta_{p_h^j}\delta_{g_j}\Big(1+\frac{\delta_{g_j}}{\sigma_j}\Big)(\delta_{p_h^j}+1)+\Lambda_{p_h^j}\delta_{g_j}T_h^j\Big)\|A_j\|_2^2\Big)\|x-y\|.\nonumber
\end{eqnarray}}
\end{proof}

\begin{remark}[Notation for Lipschitz constants]\label{rem:Lip:notation}
In Propositions~\ref{propo1:tran:breg:lipschitz}--\ref{propo2:tran:breg:lipschitz}, $L_f$, $L_{\tilde g}$, $L_g$ and $L_q$ denote Lipschitz constants of $\nabla f$, $\nabla\tilde g$ and of the translation operators $g$, $q$, respectively. In Proposition~\ref{propo5:final:diffe:Lipschitz:weak}, the bounds $\delta_f$, $\delta_{g_j}$, $\delta_{\tilde g_i}$, $\delta_{p_h^j}$ on the second derivatives and Jacobians play the role of these first-order constants (by the mean value inequality, $\nabla f$ is $\delta_f$-Lipschitz, $\tilde g_i$ is $\delta_{\tilde g_i}$-Lipschitz, and so on), whereas $\Lambda_f$, $\Lambda_{g_j}$, $\Lambda_{\tilde g_i}$, $\Lambda_{p_h^j}$ are Lipschitz constants of $\nabla^2f$, $\nabla^2g_j$, $J_{\tilde g_i}$, $J_{p_h^j}$. Finally, $L_{\nabla p_{f,g}}$ in \eqref{ineq6:propo5:final:diffe:Lipschitz} is a \emph{global} Lipschitz constant of $\nabla p_{f,g}$, used in Algorithms~\ref{Algo:GSPGA} and~\ref{Algo:CPGA}, while $L_{\mathcal{K}}$ in \eqref{ineq6:propo5:final:diffe:Lipschitz:weak} is a Lipschitz constant of $\nabla p_{f,g}$ on $\mathcal{K}$ only, used in Algorithm~\ref{Algo:GPSPGA}; in general $L_{\mathcal{K}}\neq L_{\nabla p_{f,g}}$.
\end{remark}

\subsubsection{The generalized projected simultaneous projection gradient algorithm and its convergence}
In this subsection, we assume that the set $\mathcal{K}$ defined in Proposition \ref{propo5:final:diffe:Lipschitz:weak} is a nonempty convex subset of $\mathbb{R}^n$; by Lemma~\ref{lem:convex:components}, this holds, e.g., in the Euclidean case when the maps $x\mapsto x-\tilde g_i(x)$ and $x\mapsto A_jx-p_h^j(A_jx)$ satisfy (a) or (b) of that lemma. The closedness of $\mathcal{K}$ follows from the hypotheses of Proposition \ref{propo5:final:diffe:Lipschitz:weak}, and the solution set $\varTheta$ of the problem \eqref{MSMSICFP} satisfies $\varTheta\subseteq \mathcal{K}$. We introduce a projected gradient-type algorithm for solving the problem \eqref{MSMSICFP}.

\begin{breakablealgorithm}\caption{Generalized Projected Simultaneous Projection Gradient Algorithm (GPSPGA)}
\label{Algo:GPSPGA} \begin{algorithmic}[1]\vspace{-0.2cm}
\State \textbf{Input}: $x^{0}\in\mathcal{K}$.
\State Choose $\lambda\in(0,1/L_{\mathcal{K}}]$.
\For{$k=0,1,2,\dots$}
\State Compute
\begin{eqnarray}\label{ineq1:grad:breg:gen:Algo:GPSPGA}
\begin{aligned}
x^{k+1}:=P_\mathcal{K}(x^k-\lambda \nabla p_{f,g}(x^k)).
\end{aligned}
\end{eqnarray}
\EndFor
\end{algorithmic}
\end{breakablealgorithm}

\begin{theorem}\label{MainThm1:Algo:GPSPGA}
Assume that all the hypotheses of Proposition \ref{propo5:final:diffe:Lipschitz:weak} hold and that the set $\mathcal{K}$ is nonempty and convex. Additionally, let $p_{f,g}$ be convex on $\mathcal{K}$, and suppose the solution set $\varTheta$ of \eqref{MSMSICFP} is nonempty. Let $L:=L_{\mathcal{K}}$ be the constant in \eqref{ineq6:propo5:final:diffe:Lipschitz:weak} and $\lambda\in(0,1/L]$. Then any
sequence $\{x^{k}\}_{k=0}^{\infty}$ generated by Algorithm \ref{Algo:GPSPGA}
converges in norm to a point in $\varTheta$, and
\begin{equation}\label{eq:rate:GPSPGA}
p_{f,g}(x^k)\leq\frac{\mathrm{dist}^2(x^0,\varTheta)}{2\lambda k}\;\;\text{for all}\;\;k\geq1.
\end{equation}
\end{theorem}

\begin{proof}
The set $\mathcal{K}$ is closed, because $x\mapsto P^f_{C_i(x)}(x)$ and $x\mapsto P^{g_j}_{Q_h^j(A_jx)}(A_jx)$ are continuous (Corollary~\ref{cor1:tran:breg:cont}), and $\varTheta\subseteq\mathcal{K}$ since every solution has zero residuals. Denote $g_k:=\nabla p_{f,g}(x^k)$, and let $z\in\mathcal{K}$. Since $x^k,x^{k+1}\in\mathcal{K}$ and $\mathcal{K}$ is convex, the segment $[x^k,x^{k+1}]$ lies in $\mathcal{K}$, where $\nabla p_{f,g}$ is $L$-Lipschitz; hence the descent lemma (whose proof only uses the Lipschitz continuity of the gradient along this segment) gives
$p_{f,g}(x^{k+1})\leq p_{f,g}(x^k)+\langle g_k,x^{k+1}-x^k\rangle+\frac{L}{2}\|x^{k+1}-x^k\|^2$.
By convexity, $p_{f,g}(x^k)\leq p_{f,g}(z)+\langle g_k,x^k-z\rangle$. By the characterization \eqref{proj:inq1} of $x^{k+1}=P_{\mathcal{K}}(x^k-\lambda g_k)$,
$\langle g_k,x^{k+1}-z\rangle\leq\frac1\lambda\langle x^k-x^{k+1},x^{k+1}-z\rangle=\frac{1}{2\lambda}\big(\|x^k-z\|^2-\|x^{k+1}-z\|^2-\|x^{k+1}-x^k\|^2\big)$.
Combining these three inequalities and using $L\leq1/\lambda$, we obtain, for every $z\in\mathcal{K}$,
\begin{equation}\label{eq:key:GPSPGA}
p_{f,g}(x^{k+1})\leq p_{f,g}(z)+\frac{1}{2\lambda}\big(\|x^k-z\|^2-\|x^{k+1}-z\|^2\big).
\end{equation}
Taking $z=x^k$ shows that $\{p_{f,g}(x^k)\}$ is nonincreasing. Taking $z\in\varTheta$ (so that $p_{f,g}(z)=0$) shows that $\{x^k\}$ is Fej\'er monotone with respect to $\varTheta$. Summing \eqref{eq:key:GPSPGA} with $z=P_\varTheta(x^0)$ (note that $\varTheta=\arg\min p_{f,g}$ is closed and convex) from $0$ to $k-1$ and using the monotonicity of $\{p_{f,g}(x^k)\}$, we obtain $k\,p_{f,g}(x^k)\leq\frac{1}{2\lambda}\mathrm{dist^2}(x^0,\varTheta)$, which is \eqref{eq:rate:GPSPGA} (compare with the classical analysis of the projected gradient method, e.g., \cite[Chapter 10]{ABEC2017}). Consequently, $p_{f,g}(x^k)\to0$, every accumulation point of the bounded sequence $\{x^k\}$ belongs to $\mathcal{K}$ and is a zero of $p_{f,g}$, that is, it belongs to $\varTheta$, and the Fej\'er monotonicity implies that the whole sequence converges to it.
\end{proof}

\begin{remark}
Since $\nabla p_{f,g}$ is Lipschitz continuous only on $\mathcal{K}$, Lemma~\ref{Baillon-Haddad} is not available, which is why the step sizes are restricted to $(0,1/L]$. Moreover, $P_{\mathcal{K}}$ is in general not explicit; the proof above remains valid if $\mathcal{K}$ is replaced in Algorithm~\ref{Algo:GPSPGA} by any computable nonempty, closed, and convex set $\mathcal{K}'$ with $\varTheta\subseteq\mathcal{K}'\subseteq\mathcal{K}$ (e.g., a box known to contain the solutions).
\end{remark}

\section{Numerical experiments}\label{Sec:numerics}

In this section, we illustrate the theoretical results and compare the proposed algorithms with implicit projection methods from the literature. All computations were performed in Python (NumPy/SciPy); the code that reproduces every number and figure of this section is provided in the accompanying file \texttt{msms\_icfp\_experiments.py}. The following methods are compared.

\begin{itemize}
\item \textbf{GSPGA}: Algorithm~\ref{Algo:GSPGA} with the step size $\lambda=1/L$, where $L$ is either the Lipschitz constant $L_{\nabla p_{f,g}}$ given by Proposition~\ref{propo5:final:diffe:Lipschitz} (denoted by GSPGA($L_{\rm th}$)) or a sharper (numerically computed or exact) constant.
\item \textbf{CPGA}: Algorithm~\ref{Algo:CPGA} with the $q$-steering sequence $\sigma_k=\sigma_0(1+k/q)^{-0.6}$, where $\sigma_0<2/L$ with $L$ as in Theorem~\ref{MainThm1:Algo:CPGA}.
\item \textbf{GPBCPA}: Algorithm~\ref{Algo:GPBCPA} (only for constant sets, i.e., for \eqref{MSMSCFP}).
\item \textbf{SIM-F} (simultaneous, frozen sets): at the current iterate $x^k$, every implicit set is treated as the fixed convex set $C_i(x^k)$ (respectively $Q_h^j(A_jx^k)$), and a simultaneous projection step is taken,
\begin{equation}\label{eq:SIMF}
x^{k+1}=x^k+\lambda\Big(\sum_{i=1}^r\alpha_i\big(P_{C_i(x^k)}(x^k)-x^k\big)+\sum_{j=1}^t\sum_{h=1}^{\ell_j}\beta_h^jA_j^T\big(P_{Q_h^j(A_jx^k)}(A_jx^k)-A_jx^k\big)\Big).
\end{equation}
This is the simultaneous projection scheme for the ICFP of the type studied in \cite{YCEN2016}, with the split constraints treated in the CQ fashion of \cite{CBYR2002,YCEN2005}. We use $\lambda=1/(\sum_i\alpha_i+\sum_{j,h}\beta_h^j\|A_j\|_2^2)$, which guarantees convergence when all the sets are constant.
\item \textbf{SEQ-F} (sequential, frozen sets): the iterate is successively replaced by $P_{C_i(x)}(x)$, $i=1,\dots,r$, and by $x+A_j^T\big(P_{Q_h^j(A_jx)}(A_jx)-A_jx\big)/\|A_j\|_2^2$ for all pairs $(j,h)$; one sweep through all the sets counts as one iteration. This is the sequential implicit projection method of \cite{MRGH2013,YCEN2016}.
\end{itemize}
Both frozen-set methods are implemented in their basic form, without the relaxation parameters and other refinements of the original algorithms in \cite{YCEN2016,MRGH2013}; they serve to isolate the effect of freezing the implicit sets.

In the Euclidean case $f=\frac12\|\cdot\|^2$, the contribution of the set $C_i(x)=\tilde g_i(x)+\varOmega_i$ to $\nabla p_{f,g}$ is, by Proposition~\ref{propo3:tran:breg:diff},
\begin{equation}\label{eq:grad:vs:frozen}
\nabla\varphi_i(x)=\big(I-J_{\tilde g_i}(x)\big)^T\big(x-P_{C_i(x)}(x)\big),
\end{equation}
whereas the frozen-set methods use the direction $x-P_{C_i(x)}(x)$. The factor $(I-J_{\tilde g_i}(x))^T$ accounts for the motion of the set $C_i(x)$ with $x$. When $\tilde g_i$ is expansive (for instance, if $J_{\tilde g_i}$ has eigenvalues larger than one), the frozen direction may point away from the solution set, while $-\nabla\varphi_i$ is always a descent direction for $p_{f,g}$.

\subsection{An academic example with expansive translations}\label{subsec:toy}

\begin{example}\label{Example1:numar}
Let $r=3$, $t=2$, $\ell_1=\ell_2=2$, $\mathbb{R}^n=\mathbb{R}^2$, $\mathbb{R}^{m_1}=\mathbb{R}^{m_2}=\mathbb{R}$, $A_1=(1\; 1)$, $A_2=(1\;-1)$, so that $A_1x=x_1+x_2$, $A_2x=x_1-x_2$ and $\|A_j\|_2=\sqrt{2}$. Set $\alpha_i=1$, $\beta_h^j=1$,
\begin{eqnarray}
   f(x)=\frac12\|x\|^2,\qquad g_1(y)=g_2(y)=\frac12 y^2,
\end{eqnarray}
\begin{eqnarray}
    \tilde g_1(x)=
\begin{pmatrix}
2x_1\\
2x_2
\end{pmatrix},
\qquad
\tilde g_2(x)
=
\begin{pmatrix}
x_1+\sqrt{1+x_2^2}\\
2x_1+x_2
\end{pmatrix},
\qquad
\tilde g_3(x)
=
\begin{pmatrix}
3x_1\\
2x_2
\end{pmatrix},
\end{eqnarray}
\begin{eqnarray}
\Omega_1=[-2,0]\times[-2,0],
\qquad
\Omega_2=[-2,0]\times[-3,0],
\qquad
\Omega_3=[-3,0]\times[-2,0],
\end{eqnarray}
\begin{eqnarray}
    p_1^1(y)=2y,\;\; p_2^1(y)=3y,\;\;p_1^2(y)=2y,\;\; p_2^2(y)=3y,
\end{eqnarray}
\begin{eqnarray}
    \tilde\Omega_1^1=[-3,0],\;\;\tilde\Omega_2^1=[-5,0],\;\;\tilde\Omega_1^2=[-2,2],\;\;\widetilde\Omega_2^2=[-3,3],
\end{eqnarray}
and $C_i(x)=\tilde g_i(x)+\Omega_i$, $Q_h^j(y)=p_h^j(y)+\tilde\Omega_h^j$.
All the hypotheses of Proposition~\ref{propo5:final:diffe:Lipschitz} are satisfied: $f$ and $g_j$ are $C^2$ and $1$-strongly convex Bregman functions, $\tilde{g}_i\in C^1(\mathbb{R}^2;\mathbb{R}^2)$, $p_h^j\in C^1(\mathbb{R};\mathbb{R})$, and the sets $\Omega_i$, $\tilde\Omega_h^j$ are compact boxes. Since $x\in C_i(x)$ if and only if $x-\tilde g_i(x)\in\Omega_i$ (and similarly for $Q_h^j$), a direct computation shows that the solution set of \eqref{MSMSICFP} is the polygon
\begin{equation}\label{eq:toy:Theta}
\varTheta=\big\{x\in\mathbb{R}^2\;\big\vert\;0\leq x_1\leq1.5,\;0\leq x_2\leq\sqrt3,\;x_1+x_2\leq2.5,\;|x_1-x_2|\leq1.5\big\};
\end{equation}
in particular, $(1,1)\in\varTheta$. Under the given data,
$\phi_i(x):=D_f\big(P_{C_i(x)}^f(x),x\big)=\frac{1}{2} \mathrm{dist}^2\big(x-\tilde{g}_i(x),\Omega_i\big)$ and
$\psi_h^j(x):=D_{g_j}\big(P_{Q_h^j(A_jx)}^{g_j}(A_jx),A_jx\big)=\frac{1}{2} \mathrm{dist}^2\big(A_jx-p_h^j(A_jx),\tilde\Omega_h^j\big)$, so that
\begin{equation*}
p_{f,g}(x)=\phi_1(x)+\phi_2(x)+\phi_3(x)+\psi_1^1(x)+\psi_2^1(x)+\psi_1^2(x)+\psi_2^2(x).
\end{equation*}
All seven components are convex by Lemma~\ref{lem:convex:components} (the maps $x-\tilde g_1(x)$, $x-\tilde g_3(x)$, and $A_jx-p_h^j(A_jx)$ are linear, and $x-\tilde g_2(x)=(-\sqrt{1+x_2^2},-2x_1)$ has a concave first component bounded above by $-1\leq0$ and a linear second component; see Appendix~\ref{Appendix:Example1} for the explicit form of the components), hence $p_{f,g}$ is convex, and $\nabla p_{f,g}$ is Lipschitz continuous with the constant $L_{\rm th}:=L_{\nabla p_{f,g}}=57+\sqrt{22}\approx61.69$ obtained from Proposition~\ref{propo5:final:diffe:Lipschitz} (see Appendix~\ref{Appendix:Example1}). Hence Theorems~\ref{MainThm:2}, \ref{thm:rate:convex}, and \ref{MainThm1:Algo:CPGA} apply.
\end{example}

\paragraph{Validation.} The gradient formula of Proposition~\ref{propo4:final:diffe:sum} agrees with central finite differences up to a relative error of $1.0\times10^{-9}$ at $200$ random points. Sampling the (almost everywhere defined) Hessian of $p_{f,g}$ at $4000$ random points of $[-20,20]^2$ gives the estimate $L\approx19.0$, i.e., the bound $L_{\rm th}$ is conservative by a factor of about $3.2$; this is expected, since $L_{\rm th}$ is a sum of the constants of the individual components. Note that a sampled value is only an estimate (from below) of the true Lipschitz constant; the step $\lambda=1/19$ satisfies the condition $\lambda<2/L_{\nabla p_{f,g}}$ of Theorem~\ref{MainThm:2} as long as the true constant is smaller than $38$. For CPGA we use $\sigma_0=1/24$, where $24=\max_{i,j,h}L_{i,j,h}$ is computed from the constants listed in Appendix~\ref{Appendix:Example1}.

\begin{figure}[htbp]
\centering
\includegraphics[width=0.6\textwidth]{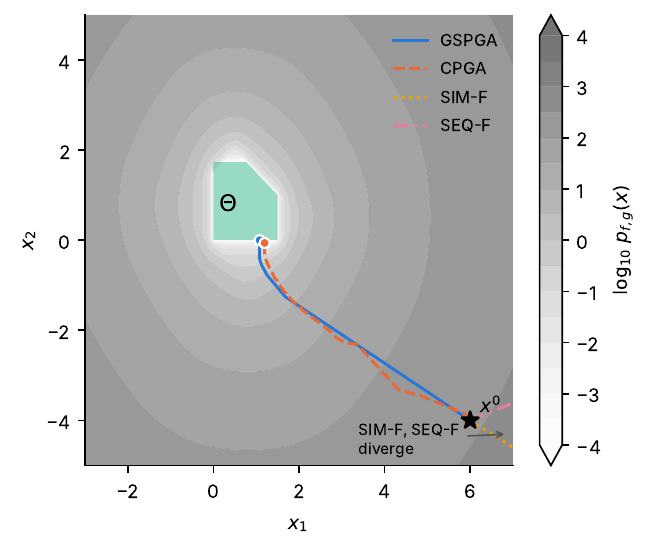}
\caption{Example~\ref{Example1:numar}: level sets of $\log_{10}p_{f,g}$, the solution set $\varTheta$ (shaded), and the first $400$ iterates of the methods starting from $x^0=(6,-4)$. Both frozen-set methods leave the displayed region.}\label{fig:toy:traj}
\end{figure}

\begin{figure}[htbp]
\centering
\includegraphics[width=\textwidth]{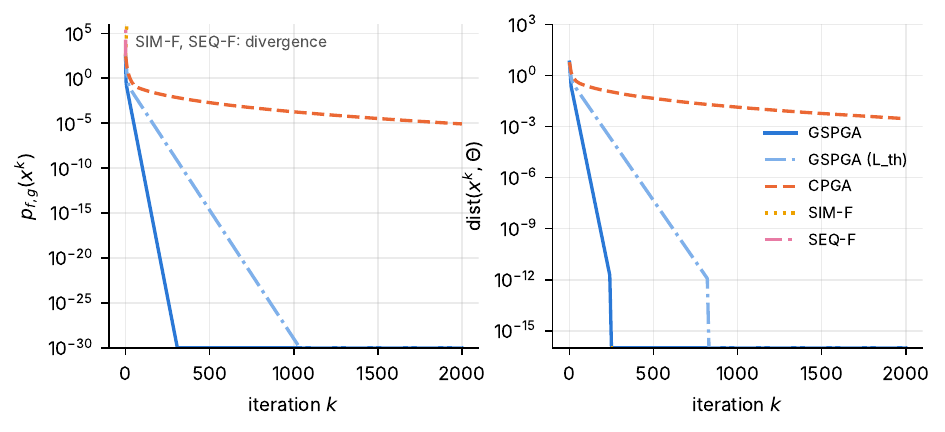}
\caption{Example~\ref{Example1:numar}: proximity values $p_{f,g}(x^k)$ (left) and distances $\mathrm{dist}(x^k,\varTheta)$ (right) for $x^0=(6,-4)$. The curves of the frozen-set methods end when the iterates overflow.}\label{fig:toy:conv}
\end{figure}

\paragraph{Behaviour of the iterates.} Figure~\ref{fig:toy:traj} shows the trajectories starting from $x^0=(6,-4)$, and Figure~\ref{fig:toy:conv} shows the proximity values and the distances to $\varTheta$. GSPGA and CPGA move towards $\varTheta$, while both frozen-set methods diverge. The reason is the expansive translation $\tilde g_1(x)=2x$: for $x_1>2$, the frozen displacement of the first set is $\big(P_{C_1(x)}(x)-x\big)_1=x_1-2>0$, which pushes the iterate away from $\varTheta$, whereas, by \eqref{eq:grad:vs:frozen}, $-\big(\nabla\phi_1(x)\big)_1=-(x_1-2)<0$. GSPGA with $\lambda=1/L$ reaches the solution set after a few hundred iterations, and with the conservative step $\lambda=1/L_{\rm th}$ it needs about three times more iterations. After reaching the boundary of $\varTheta$, the iterates remain fixed, since $\nabla p_{f,g}=0$ on $\varTheta$. CPGA converges much more slowly because of the diminishing steering steps, in agreement with the fact that Theorem~\ref{MainThm1:Algo:CPGA} provides no rate.

\begin{figure}[htbp]
\centering
\includegraphics[width=\textwidth]{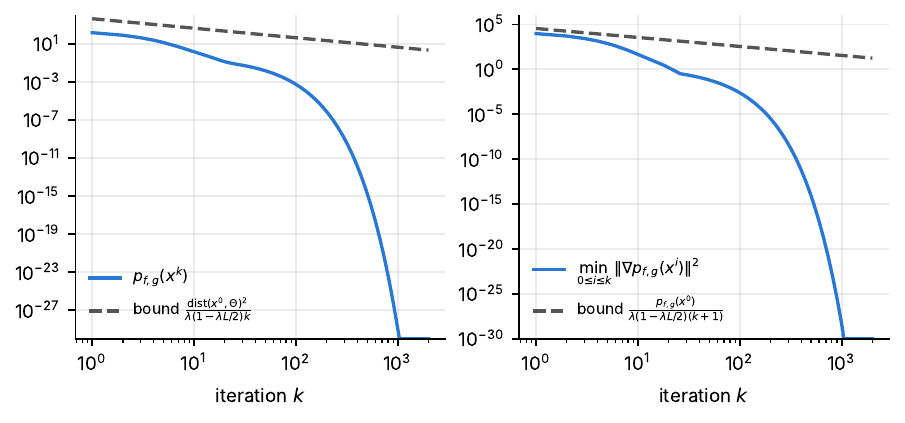}
\caption{Example~\ref{Example1:numar}, GSPGA with $\lambda=1/L_{\rm th}$: comparison of $p_{f,g}(x^k)$ with the bound \eqref{eq:rate:convex:a} (left) and of $\min_{0\leq i\leq k}\|\nabla p_{f,g}(x^i)\|^2$ with the bound \eqref{eq:rate:nonconvex} (right).}\label{fig:toy:rate}
\end{figure}

\paragraph{Convergence rates.} Figure~\ref{fig:toy:rate} compares the observed values with the bounds of Theorems~\ref{thm:rate:convex}(a) and \ref{thm:rate:nonconvex} for $\lambda=1/L_{\rm th}$. Both bounds hold along the whole run. The bounds are $O(1/k)$, whereas the observed decay is geometric (the curves become vertical in the log--log scale), which is the behaviour predicted by Theorem~\ref{thm:rate:convex}(c) when the quadratic growth condition \eqref{eq:QG} holds; the observed geometric decay indicates that \eqref{eq:QG} holds for this example.

\paragraph{Robustness with respect to the starting point.} We ran each method from $300$ starting points drawn uniformly from $[-10,10]^2$, with $3000$ iterations (sweeps for SEQ-F). A run is called successful if $\mathrm{dist}(x^{3000},\varTheta)\leq10^{-5}$. The results are reported in Table~\ref{tab:toy}. GSPGA is successful from every starting point. CPGA always approaches $\varTheta$ (it is convergent by Theorem~\ref{MainThm1:Algo:CPGA}), but within $3000$ iterations it attains the accuracy $10^{-5}$ only in about one third of the runs. The frozen-set methods succeed only when $x^0$ is (essentially) already in $\varTheta$.

\begin{table}[htbp]
\centering
\caption{Example~\ref{Example1:numar}: results from $300$ random starting points in $[-10,10]^2$ ($3000$ iterations).}\label{tab:toy}
\begin{tabular}{lccc}
\toprule
Method & success & median its.\ to $p_{f,g}\leq10^{-12}$ & ms/run$^*$\\
\midrule
GSPGA ($\lambda=1/L$, $L\approx19$) & 100.0\% & 125 & 523\\
GSPGA($L_{\rm th}$) & 100.0\% & 420 & 530\\
CPGA & 32.7\% & 746 & 309\\
SIM-F & 0.0\% & -- & 411\\
SEQ-F & 0.3\% & -- & 393\\
\bottomrule
\end{tabular}

\smallskip
{\footnotesize $^*$Including the evaluation of $p_{f,g}$ at every iterate.}
\end{table}

\subsection{Application: cooperative sensor network localization}\label{subsec:sensor}

The ICFP was originally motivated by the cooperative positioning problem in wireless sensor networks \cite{MRGH2013}; see also \cite{BIYE2004,PATW2005} for background on network localization. A number $N$ of sensors with unknown positions $x_1,\dots,x_N\in\mathbb{R}^2$ and $M$ anchors with known positions $a_1,\dots,a_M$ are deployed in the region $[0,1]^2$. Two nodes within the communication radius $R$ of each other measure their distance, and the measurements provide upper bounds: $\|x_i-x_j\|\leq\rho_{ij}$ for every sensor--sensor link $\{i,j\}\in\mathcal{E}$ and $\|x_i-a_k\|\leq\rho_{ik}^a$ for every sensor--anchor link $(i,k)\in\mathcal{E}_a$. The feasible region of a sensor depends on the unknown positions of its neighbours, which is precisely the implicit structure of \eqref{MSMSICFP}.

\paragraph{Formulation as an MSMS-ICFP.} Let $x=(x_1,\dots,x_N)\in\mathbb{R}^{2N}$. For every link $e=\{i,j\}\in\mathcal{E}$ with $i<j$, define the translation-type operator
\begin{equation}\label{eq:sensor:C}
C_e(x):=g_e(x)+\varOmega_e,\qquad g_e(x):=E_ix_j,\qquad \varOmega_e:=\{z\in\mathbb{R}^{2N}\mid\|z_i\|\leq\rho_{ij}\},
\end{equation}
where $E_i:\mathbb{R}^2\to\mathbb{R}^{2N}$ places a vector in the $i$-th block. Then $x\in C_e(x)$ if and only if $\|x_i-x_j\|\leq\rho_{ij}$. The deployment region gives the constant set $C_0:=[0,1]^{2N}$. The anchor constraints are split constraints: with $A_i x:=x_i$ ($A_i:\mathbb{R}^{2N}\to\mathbb{R}^2$, $\|A_i\|_2=1$), we require $A_ix\in Q_{ik}:=a_k+B(0,\rho_{ik}^a)$ for $(i,k)\in\mathcal{E}_a$. Here $g_e$ is linear, $\varOmega_e$ is closed and convex but unbounded (covered by Proposition~\ref{propo3:tran:breg:diff}), and we take all the Bregman functions Euclidean and all the weights equal to one. Then
\begin{equation}
p_{f,g}(x)=\frac12\sum_{\{i,j\}\in\mathcal{E}}\big[\|x_i-x_j\|-\rho_{ij}\big]_+^2+\frac12\sum_{(i,k)\in\mathcal{E}_a}\big[\|x_i-a_k\|-\rho^a_{ik}\big]_+^2+\frac12\mathrm{dist}^2\big(x,[0,1]^{2N}\big),
\end{equation}
which is convex by Lemma~\ref{lem:convex:components}(a), and, by \eqref{eq:grad:vs:frozen}, the link $e=\{i,j\}$ contributes $r_e$ to the $i$-th block and $-r_e$ to the $j$-th block of $\nabla p_{f,g}(x)$, where $r_e:=(x_i-x_j)-P_{B(0,\rho_{ij})}(x_i-x_j)$. Writing $p_{f,g}(x)=\frac12\sum_\ell\mathrm{dist}^2(B_\ell x,S_\ell)$ with linear maps $B_\ell$, one sees that $\nabla p_{f,g}$ is Lipschitz continuous with constant $L=\lambda_{\max}\big(\sum_\ell B_\ell^TB_\ell\big)=\lambda_{\max}\big((\mathcal{L}+D_a+I)\otimes I_2\big)$, where $\mathcal{L}$ is the Laplacian of the graph $(\{1,\dots,N\},\mathcal{E})$ and $D_a$ is the diagonal matrix of the numbers of anchor links. We also compute the constant $L_{\rm th}$ of Proposition~\ref{propo5:final:diffe:Lipschitz}, obtained by summing the Lipschitz constants of the components: $2\sqrt2$ for each link ($L_1=M_1=\sqrt2$), and $1$ for each anchor constraint and for the box (constant sets, for which the Jacobian term vanishes). Since $x\in C_e(x)$ if and only if $x_i-x_j\in B(0,\rho_{ij})$, the problem can also be written as an MSMS-CFP with $A_ex:=x_i-x_j$; we use this reformulation for GPBCPA, with $K=4$ blocks and the dominating function $u=\frac c2\|\cdot\|^2$ (see the discussion after Definition~\ref{def1:thm1:MSMSCFP}). In particular, this application could also be treated as an explicit CFP; its purpose here is to compare the implicit-set methods under controlled conditions, and to show the effect of freezing the implicit sets.

The frozen-set method SIM-F applied to \eqref{eq:sensor:C} moves only the $i$-th sensor for the link $\{i,j\}$, since $C_e(x^k)$ restricts only the $i$-th block. Consequently, a fixed point of SIM-F only requires the residuals acting on each sensor to cancel, which does not imply feasibility. If every link is listed twice, as $(i,j)$ and $(j,i)$, then the SIM-F step coincides with the GSPGA step (we verified numerically that both sequences agree to machine precision); in other words, GSPGA performs this symmetrization automatically through the Jacobian term in \eqref{eq:grad:vs:frozen}.

\paragraph{Setup.} We take $N=60$ sensors drawn uniformly in $[0,1]^2$, $M=9$ anchors placed at the nine points of the grid $\{0.1,0.5,0.9\}^2$, and $R=0.3$, which yields $|\mathcal{E}|=354$ sensor--sensor links and $|\mathcal{E}_a|=96$ sensor--anchor links. The measured distances are $\hat d=d+\varepsilon$, where $d$ is the true distance and $\varepsilon$ is Gaussian noise with standard deviation $0.02$, truncated to $[-0.06,0.06]$, and we set $\rho=\hat d+0.06\geq d$, so that the true configuration belongs to $\varTheta$. We obtain $L=23.01$ and $L_{\rm th}=1098.3$. All the methods start from the same configuration, in which each sensor is placed at the center $(0.5,0.5)$ of the region plus a random perturbation of size $10^{-2}$. Iterations are measured in passes over all the sets: one GSPGA or SIM-F iteration, one sweep of CPGA or SEQ-F, and one cycle of GPBCPA through the $K$ blocks. Each method is run for a fixed budget of passes; to compare the methods at a common accuracy, Table~\ref{tab:sensor} also reports the first pass at which the maximal constraint violation drops below $10^{-8}$.

\begin{figure}[htbp]
\centering
\includegraphics[width=\textwidth]{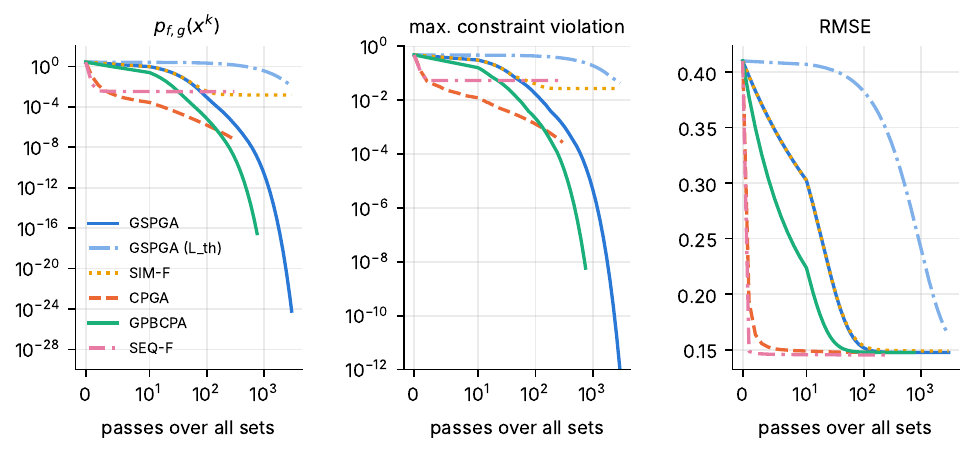}
\caption{Sensor network localization: proximity value (left), maximal constraint violation (middle), and RMSE (right) versus the number of passes over all the sets.}\label{fig:sensor:conv}
\end{figure}

\begin{figure}[htbp]
\centering
\includegraphics[width=0.5\textwidth]{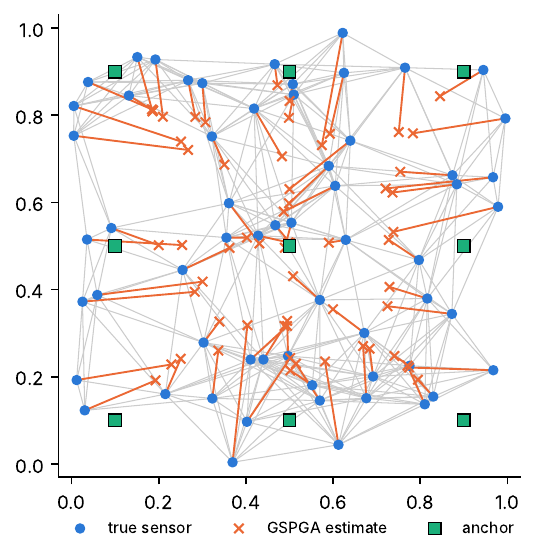}
\caption{Sensor network localization: true sensor positions, anchors, communication links (grey), and the feasible configuration computed by GSPGA; each estimate is joined to the corresponding true position.}\label{fig:sensor:net}
\end{figure}

\paragraph{Results.} Table~\ref{tab:sensor} and Figure~\ref{fig:sensor:conv} report the proximity value, the maximal constraint violation, and the root-mean-square error $\mathrm{RMSE}=\big(\frac1N\sum_i\|x_i-x_i^{\rm true}\|^2\big)^{1/2}$. GSPGA with $\lambda=1/L$ and GPBCPA find a configuration satisfying all the constraints up to a negligible violation, and the linear decay of $p_{f,g}$ is again visible. GPBCPA, which uses the larger block-wise step $1/c$ with $c=9.67$, reaches the accuracy $10^{-8}$ after 723 passes, compared with 1796 passes for GSPGA. CPGA decreases the proximity function quickly in the first sweeps, but then slows down owing to the steering steps. GSPGA with the conservative constant $L_{\rm th}$ is far from convergence after $3000$ iterations, which shows the practical importance of sharp Lipschitz constants (cf. Remark~\ref{rem:rates}(iii)). The frozen-set methods SIM-F and SEQ-F stagnate at infeasible configurations (maximal violations of order $10^{-2}$), as anticipated above. Their RMSE values are close to those of the feasible solutions (SEQ-F even attains a slightly smaller value), but the computed configurations are not consistent with the measurements. The resulting localization is shown in Figure~\ref{fig:sensor:net}. Note that the convex feasibility formulation only uses upper bounds on the distances, so that every point of $\varTheta$ is consistent with the data and the RMSE of a feasible point reflects the size of $\varTheta$ rather than the quality of the algorithm; improving the accuracy, e.g., by superiorization of the feasibility-seeking iterations with respect to a suitable target function, is a natural topic for further research.

\begin{table}[htbp]
\centering
\caption{Sensor network localization ($N=60$, $M=9$): number of passes, first pass with maximal violation $\leq10^{-8}$ (``--'': not reached), and values after the last pass.}\label{tab:sensor}
\small
\begin{tabular}{lcccccc}
\toprule
Method & passes & viol.\ $\leq10^{-8}$ & $p_{f,g}$ & max. violation & RMSE & time (s)\\
\midrule
GSPGA ($\lambda=1/L$) & 3000 & 1796 & $6.2\times10^{-25}$ & $8.4\times10^{-13}$ & 0.1475 & 0.58\\
GSPGA($L_{\rm th}$) & 3000 & -- & $1.1\times10^{-2}$ & $4.3\times10^{-2}$ & 0.1630 & 0.55\\
CPGA & 300 & -- & $6.0\times10^{-8}$ & $2.7\times10^{-4}$ & 0.1475 & 1.36\\
GPBCPA & 750 & 723 & $2.9\times10^{-17}$ & $5.9\times10^{-9}$ & 0.1475 & 0.26\\
SIM-F & 3000 & -- & $1.5\times10^{-3}$ & $2.7\times10^{-2}$ & 0.1488 & 0.68\\
SEQ-F & 300 & -- & $3.2\times10^{-3}$ & $5.4\times10^{-2}$ & 0.1453 & 0.46\\
\bottomrule
\end{tabular}
\end{table}

\section{Conclusion}
\label{Sec:conclusion}

In this work, we have developed the theory for the multi-split multiple-set implicit convex feasibility problem (introduced in \cite{SDYCEN2026}) within the Bregman projection framework.
The main contributions of this work are an explicit formula for implicit Bregman projections onto translation-type sets, their continuity and Lipschitz continuity, a gradient formula for the associated generalized proximity function, and simultaneous, cyclic, and block-iterative algorithms under both global and local assumptions.
We have established the convergence of the sequences generated by these algorithms and, for the simultaneous and projected schemes, sublinear $O(1/k)$ rates, together with a linear rate under a quadratic growth condition.
The block-iterative Bregman algorithm developed for the multi-split multiple-sets convex feasibility problem suggests a new direction for extending such methods to more general split problems.
The numerical experiments show that, because the gradient of the proximity function takes the motion of the implicit sets into account (through the Jacobians of the translation operators), the proposed methods converge in situations where frozen-set implicit projection methods diverge or stagnate.
Several questions remain open: the convergence of a projected cyclic variant of Algorithm~\ref{Algo:CPGA} under the local assumptions of Proposition~\ref{propo5:final:diffe:Lipschitz:weak}, verifiable sufficient conditions for the quadratic growth condition \eqref{eq:QG}, rates for the cyclic and block-iterative schemes, the extension of the analysis to general (non-translation-type) point-to-set operators, and the combination of the proposed feasibility-seeking methods with superiorization in order to improve, e.g., the localization accuracy in the application of Section~\ref{subsec:sensor}.

\bigskip{}
\noindent\textbf{Acknowledgments.} This work is supported in part by the Israel Science Foundation (ISF), Grant No. 992/23.
The authors would like to express their special thanks to Professor Charles L. Byrne, Department of Mathematics and Statistics at University of Massachusetts, Lowell, for his valuable comments and suggestions, which helped to improve an earlier version of this work.

%===================================================
% \section*{Declarations}

% \textbf{Funding.} This work is supported in part by the Israel Science Foundation (ISF), Grant No. 992/23.

% \bigskip{}
% \noindent\textbf{Declaration of interests.} The authors declare that they have no known competing financial interests

% \bigskip{}
% \noindent\textbf{Disclosure statement.} No potential conflict of interest was reported by the authors.

% % \bigskip{}
% % \noindent\textbf{Competing interests.} The authors declare that they have no competing interests.

% \bigskip{}
% \noindent\textbf{Data availability statement.} No data sets were generated or analysed in this study other than the synthetic data described in Section~\ref{Sec:numerics}. The Python code that reproduces all numerical results and figures is available from the corresponding author upon reasonable request.

%===================================================================

%===========================================================================
\newpage
\appendix
\section*{Appendix}

\section{Details for Example \ref{Example1:numar}}\label{Appendix:Example1}

Write $[s]_+:=\max\{0,s\}$, so that $\mathrm{dist}(s,[a,b])=[s-b]_++[a-s]_+$ for $s\in\mathbb{R}$, and recall that $|[a]_+-[b]_+|\leq|a-b|$. By Proposition~\ref{propo1:tran:breg}, $P_{C_i(x)}^f(x)-x=P_{\Omega_i}(s_i)-s_i$ with $s_i:=x-\tilde g_i(x)$, where
$s_1=(-x_1,-x_2)$, $s_2=(-\sqrt{1+x_2^2},-2x_1)$ and $s_3=(-2x_1,-x_2)$; similarly, $P_{Q_h^j(A_jx)}^{g_j}(A_jx)-A_jx=P_{\tilde\Omega_h^j}(u_h^j)-u_h^j$ with $u_h^j:=A_jx-p_h^j(A_jx)$, i.e., $u_1^1=-(x_1+x_2)$, $u_2^1=-2(x_1+x_2)$, $u_1^2=-(x_1-x_2)$ and $u_2^2=-2(x_1-x_2)$. Hence,
\begin{eqnarray*}
\phi_1(x)&=&\frac12\left[[-x_1]_+^2+[x_1-2]_+^2+[-x_2]_+^2+[x_2-2]_+^2\right],\\
\phi_2(x)&=&\frac{1}{2}\left[\left[\sqrt{1+x_2^2}-2\right]_+^2+\left(-[-2x_1]_+ +[2x_1-3]_+\right)^2\right],\\
\phi_3(x)&=&\frac12\left[[-2x_1]_+^2+[2x_1-3]_+^2+[-x_2]_+^2+[x_2-2]_+^2\right],
\end{eqnarray*}
and
\begin{eqnarray*}
\psi_1^1(x)&=&\frac12\left[[-x_1-x_2]_+^2+[x_1+x_2-3]_+^2\right],\\
\psi_2^1(x)&=&\frac12\left[[-2(x_1+x_2)]_+^2+[2(x_1+x_2)-5]_+^2\right],\\
\psi_1^2(x)&=&\frac12\left[[-x_1+x_2-2]_+^2+[x_1-x_2-2]_+^2\right],\\
\psi_2^2(x)&=&\frac12\left[[-2(x_1-x_2)-3]_+^2+[2(x_1-x_2)-3]_+^2\right].
\end{eqnarray*}

Finally,
\begin{equation*}
\begin{aligned}2p_{f,g}(x)={}&[-x_1]_+^2+[x_1-2]_+^2+[-x_2]_+^2+[x_2-2]_+^2\\&
+\left[\sqrt{1+x_2^2}-2\right]_+^2+\left(-[-2x_1]_+ +[2x_1-3]_+\right)^2\\&
+[-2x_1]_+^2+[2x_1-3]_+^2+[-x_2]_+^2+[x_2-2]_+^2\\&
+[-x_1-x_2]_+^2+[x_1+x_2-3]_+^2\\&
+[-2(x_1+x_2)]_+^2+[2(x_1+x_2)-5]_+^2\\&
+[-x_1+x_2-2]_+^2+[x_1-x_2-2]_+^2\\&
+[-2(x_1-x_2)-3]_+^2+[2(x_1-x_2)-3]_+^2.
\end{aligned}
\end{equation*}

Now, we obtain
\begin{eqnarray*}
J_{\tilde g_1}(x)^T\left(\nabla f\left(P_{C_1(x)}^f(x)\right)-\nabla f(x)\right)&=&
\begin{pmatrix}
2\left(-[-x_1]_+ + [x_1-2]_+\right)\\[1mm]
2\left(-[-x_2]_+ + [x_2-2]_+\right)
\end{pmatrix},\\
J_{\tilde g_2}(x)^T\left(\nabla f\left(P_{C_2(x)}^f(x)\right)-\nabla f(x)\right)&=&
\begin{pmatrix}
\left[\sqrt{1+x_2^2}-2\right]_+
-2[-2x_1]_+
+2[2x_1-3]_+
\\[3mm]
\dfrac{x_2}{\sqrt{1+x_2^2}}
\left[\sqrt{1+x_2^2}-2\right]_+
-[-2x_1]_+
+[2x_1-3]_+
\end{pmatrix},\\
J_{\tilde g_3}(x)^T\left(\nabla f\left(P_{C_3(x)}^f(x)\right)-\nabla f(x)\right)&=&
\begin{pmatrix}
3\left(-[-2x_1]_+ + [2x_1-3]_+\right)\\[2mm]
2\left(-[-x_2]_+ + [x_2-2]_+\right)
\end{pmatrix}
\end{eqnarray*}
are $2,\sqrt{22},6$-Lipschitz, respectively (for the second map, write $a(x_2):=[\sqrt{1+x_2^2}-2]_+$ and $b(x_1):=-[-2x_1]_++[2x_1-3]_+$; then $|b'|\leq2$, $|a'|\leq1$, and, where $a>0$, $\big|\frac{d}{dx_2}\frac{x_2a(x_2)}{\sqrt{1+x_2^2}}\big|=1-2(1+x_2^2)^{-3/2}\leq1$, so the Frobenius norm of the almost everywhere defined Jacobian is at most $\sqrt{16+1+4+1}=\sqrt{22}$), and
\begin{eqnarray*}
\nabla^2 f(x)\left(P_{C_1(x)}^f(x)-x\right)&=&
\begin{pmatrix}
-[-x_1]_+ + [x_1-2]_+\\[1mm]
-[-x_2]_+ + [x_2-2]_+
\end{pmatrix},\\
\nabla^2 f(x)\left(P_{C_2(x)}^f(x)-x\right)&=&
\begin{pmatrix}
\left[\sqrt{1+x_2^2}-2\right]_+
\\[2mm]
-[-2x_1]_+ + [2x_1-3]_+
\end{pmatrix},\\
\nabla^2 f(x)\left(P_{C_3(x)}^f(x)-x\right)&=&
\begin{pmatrix}
-[-2x_1]_+ + [2x_1-3]_+\\[2mm]
-[-x_2]_+ + [x_2-2]_+
\end{pmatrix}
\end{eqnarray*}
are $1,2,2$-Lipschitz, respectively.
Similarly, we obtain
\begin{eqnarray*}
J_{p_1^1}(A_1x)^T\left(\nabla g_1\left(P_{Q_1^1(A_1x)}^{g_1}(A_1x)\right)-\nabla g_1(A_1x)\right)&=&
2\left(-[-(x_1+x_2)]_++[x_1+x_2-3]_+\right),\\
J_{p_2^1}(A_1x)^T\left(\nabla g_1\left(P_{Q_2^1(A_1x)}^{g_1}(A_1x)\right)-\nabla g_1(A_1x)\right)&=&
3\left(-[-2(x_1+x_2)]_++[2(x_1+x_2)-5]_+\right),\\
J_{ p_1^2}(A_2x)^T\left(\nabla g_2\left(P_{Q_1^2(A_2x)}^{g_2}(A_2x)\right)-\nabla g_2(A_2x)\right)&=&
2\left(-[-(x_1-x_2)-2]_++[x_1-x_2-2]_+\right),\\
J_{p_2^2}(A_2x)^T\left(\nabla g_2\left(P_{Q_2^2(A_2x)}^{g_2}(A_2x)\right)-\nabla g_2(A_2x)\right)&=&
3\left(-[-2(x_1-x_2)-3]_++[2(x_1-x_2)-3]_+\right)
\end{eqnarray*}
are $2\sqrt{2},6\sqrt{2},2\sqrt{2},6\sqrt{2}$-Lipschitz, respectively, and
\begin{eqnarray*}
\nabla^2g_1(A_1x)\left(P_{Q_1^1(A_1x)}^{g_1}(A_1x)-A_1x\right)&=&
-[-(x_1+x_2)]_++[x_1+x_2-3]_+,\\
\nabla^2g_1(A_1x)\left(P_{Q_2^1(A_1x)}^{g_1}(A_1x)-A_1x\right)&=&
-[-2(x_1+x_2)]_++[2(x_1+x_2)-5]_+,\\
\nabla^2g_2(A_2x)\left(P_{Q_1^2(A_2x)}^{g_2}(A_2x)-A_2x\right)&=&
-[-(x_1-x_2)-2]_++[x_1-x_2-2]_+,\\
\nabla^2g_2(A_2x)\left(P_{Q_2^2(A_2x)}^{g_2}(A_2x)-A_2x\right)&=&
-[-2(x_1-x_2)-3]_++[2(x_1-x_2)-3]_+
\end{eqnarray*}
are $\sqrt{2}, 2\sqrt{2}, \sqrt{2}, 2\sqrt{2}$-Lipschitz, respectively.


\begin{thebibliography}{10}

\bibitem{RAHA1989} Aharoni, R., Censor, Y.: Block-iterative projection methods for parallel computation of solutions to convex feasibility problems. \textit{Linear Algebra and its Applications} {\bf 120}, 165--175 (1989)

\bibitem{HHBA1996} Bauschke, H.H., Borwein, J.M.: On projection algorithms
for solving convex feasibility problems. \textit{Society for Industrial and Applied Mathematics Review} {\bf 38}(3), 367--426 (1996)

\bibitem{HHBA2006} Bauschke, H.H., Combettes, P.L., Noll, D.: Joint
minimization with alternating Bregman proximity operators. \textit{Pacific Journal of Optimization} {\bf 2}(3), 401--424 (2006)

\bibitem{HHBA2011} Bauschke, H.H., Combettes, P.L.: \textit{Convex
Analysis and Monotone Operator Theory in Hilbert Spaces} (Second Edition).
Springer, New York (2017)

\bibitem{ABEC2017} Beck, A.: \textit{First-Order Methods in Optimization}. Society for Industrial and Applied Mathematics, Philadelphia, PA (2017)

% \bibitem{ABEC2014} \sd{Beck, A.: \textit{Introduction to Nonlinear Optimization: Theory, Algorithms, and Applications with MATLAB}. Society for Industrial and Applied Mathematics (2014)}

\bibitem{BIYE2004} Biswas, P., Ye, Y.: Semidefinite programming for ad hoc wireless sensor network localization. In: \textit{Proceedings of the 3rd International Symposium on Information Processing in Sensor Networks (IPSN)}, pp. 46--54 (2004)

\bibitem{LMBR1965} Bregman, L.M.: Finding the common point of convex sets by the method of successive projection. \textit{Doklady Akademii Nauk SSSR} {\bf 162}(3), 487--490 (1965)

\bibitem{LMBR1967} Bregman, L.M.: The relaxation method of finding the common point of convex sets and its application to the solution of problems in convex programming. \textit{USSR Computational Mathematics and Mathematical Physics} {\bf 7}(3), 200--217 (1967)

\bibitem{BUTN1990} Butnariu, D., Censor, Y.: On the behavior of a
block-iterative projection method for solving convex feasibility problems.
\textit{International Journal of Computer Mathematics} {\bf 34}(1--2), 79--94 (1990)

\bibitem{RSBU2021} Burachik, R.S., Dao, M.N., Lindstrom, S.B.: Generalized
Bregman envelopes and proximity operators. \textit{Journal of Optimization Theory and Applications} {\bf 190}(3), 744--778 (2021)

\bibitem{CBYR2001} Byrne, C., Censor, Y.: Proximity function minimization
using multiple Bregman projections, with applications to split feasibility
and Kullback--Leibler distance minimization. \textit{Annals of Operations Research}
{\bf 105}, 77--98 (2001)

\bibitem{CBYR22001} Byrne, C.: Bregman-Legendre multidistance projection algorithms for convex feasibility and optimization. \textit{Studies in Computational Mathematics} {\bf 8}, 87--99 (2001)

\bibitem{CBYR2002} Byrne, C.: Iterative oblique projection onto convex sets and the split feasibility problem. \textit{Inverse Problems} {\bf 18}(2), 441--453 (2002)

\bibitem{BYRN2004} Byrne, C.: A unified treatment of some iterative
algorithms in signal processing and image reconstruction. \textit{Inverse
Problems} {\bf 20}(1), 103--120 (2004)

% \bibitem{CBYR2024} Byrne, C.: Multi-distance block-iterative Bregman-retraction algorithms: Part II. To appear (2024).

\bibitem{ACEG2012} Cegielski, A.: \textit{Iterative Methods for Fixed Point Problems in Hilbert Spaces}. Springer, Heidelberg, Germany (2012)

\bibitem{YCEN1988} Censor, Y.: Parallel applications of block-iterative methods in medical imaging and radiation therapy. \textit{Mathematical Programming}, {\bf 42}, 307--325 (1988)

\bibitem{SDYCEN2026} Censor, Y., Dey, S.: Multi-split multiple-sets implicit convex feasibility seeking as split inverse problems. \textit{Communicated} (2026) (\href{https://arxiv.org/abs/2610.39991}{arXiv:2610.39991})

\bibitem{YCEN2005} Censor, Y., Elfving, T., Kopf, N., Bortfeld, T.:
The multiple-sets split feasibility problem and its applications for
inverse problems. \textit{Inverse Problems} {\bf 21}(6), 2071--2084
(2005)

\bibitem{YCEN2016} Censor, Y., Gibali, A., Lenzen, F., Schnörr, C.:
The implicit convex feasibility problem and its application to adaptive
image denoising. \textit{Journal of Computational Mathematics} {\bf 34}(6), 610--625
(2016)

\bibitem{YCEN2004} Censor, Y., Pierro, A.R.D., Zaknoon, M.: Steered sequential projections for the inconsistent convex feasibility problem. \textit{ Nonlinear Analysis: Theory, Methods $\&$ Applications} {\bf 59}(3), 385--405 (2004)

\bibitem{YCEN22012} Censor, Y., Gibali, A., Reich, S.: Algorithms
for the split variational inequality problem. \textit{Numerical Algorithms}
{\bf 59}(2), 301--323 (2012)

\bibitem{YCEN1982} Censor, Y., Lent, A.: Cyclic subgradient projections.
\textit{Mathematical Programming} {\bf 24}(1), 233--235 (1982)

\bibitem{YCEN1997} Censor, Y., Zenios, S.A.: \textit{Parallel Optimization:
Theory, Algorithms, and Applications}. Oxford University Press, New
York (1997)

\bibitem{DEUT1992} Deutsch, F.: The method of alternating orthogonal projections. In: Singh, S.P. (ed.) \textit{Approximation Theory, Spline Functions and Applications}, NATO ASI Series C, Kluwer Academic Publishers, Dordrecht. {\bf 356}, 105--121 (1992)

\bibitem{SDEY2024} Dey, S., Izuchukwu, C., Taiwo, A., Reich, S.:
New iterative algorithms for solving split variational inclusions.
\textit{Journal of Global Optimization} {\bf 91}(3), 587--609 (2025) %(https://doi.org/10.1007/s10898-024-01444-7)

\bibitem{SDEY22024} Dey, S., Reich, S.: Levitin--Polyak well-posedness
of split multivalued variational inequalities. \textit{Optimization}
{\bf 74}(13), 3081--3111 (2025) %(https://doi.org/10.1080/02331934.2024.2358408)

\bibitem{SDEYSR2026} Dey, S., Reich, S.: Weak quasi-convexity: first- and second-order conditions, monotonicity, and quadratic functions. \textit{Communicated} (2026) (\href{https://arxiv.org/abs/2611.13878}{arXiv:2611.13878})

\bibitem{TLDI2026} Dinh, T.L., Jansen, G.S.M., Luke, D.R.: Cyclic relaxed Douglas-Rachford splitting for inconsistent nonconvex feasibility. \textit{Journal of Optimization Theory and Applications} {\bf 209}(2), 1--27 (2026) (Article no. 50)

% \bibitem{DOLI1984} Dolidze, Z.O.: Convergence of an analogue of the gradient method. \textit{Èkonomika i Matematicheskie Metody} {\bf 20}(4), 755--758 (1984)

\bibitem{EREM1969} Eremin, I.I.: Fejér mappings and convex programming.
\textit{Siberian Mathematical Journal} {\bf 10}(5), 762--772 (1969)

\bibitem{FUKU1992} Fukushima, M.: Equivalent differentiable optimization
problems and descent methods for asymmetric variational inequality
problems. \textit{Mathematical Programming} {\bf 53}(1--3), 99--110 (1992)


\bibitem{MRGH2013} Gholami, M.R., Tetruashvili, L., Ström, E.G.,
Censor, Y.: Cooperative wireless sensor network positioning via implicit
convex feasibility. \textit{IEEE Transactions on Signal Processing} {\bf 61}(23),
5830--5840 (2013)

\bibitem{GUBI1967} Gubin, L.G., Polyak, B.T., Raik, E.V.: The method
of projection for finding the common point of convex sets. \textit{USSR Computational Mathematics and Mathematical Physics} {\bf 7}(6), 1--24 (1967)

\bibitem{HANS1988} Han, S.P.: A successive projection method. \textit{Mathematical Programming} {\bf 40}(1--3), 1--14 (1988)

\bibitem{SKAR1969} Karamardian, S.: The nonlinear complementarity problem with applications, part 2. \textit{Journal of Optimization Theory and Applications} {\bf 4}(3), 167--181 (1969)

\bibitem{KUBO2009} Kubota, K., Fukushima, M.: Gap function approach
to the generalized Nash equilibrium problem. \textit{Journal of Optimization Theory and Applications} {\bf 144}(3), 511--531 (2009)

\bibitem{JLI2019} Li, J.: Split equilibrium problems for related
games and applications to economic theory. \textit{Optimization} {\bf 68}(6),
1203--1222 (2019)

\bibitem{GLOP2012} Lopez, G., Martin-Marquez, V., Wang, F., Xu, H.K.:
Solving the split feasibility problem without prior knowledge of matrix
norms. \textit{Inverse Problems} {\bf 28}(8), 1--18 (2012) (Article no. 085004)

\bibitem{EMAS2007} Masad, E., Reich, S.: A note on the multiple-set
split convex feasibility problem in Hilbert space. \textit{Journal of Nonlinear and Convex Analysis} \textbf{8}(3), 367--371 (2007)

\bibitem{AMOU2014} Moudafi, A., Thakur, B.S.: Solving proximal split
feasibility problems without prior knowledge of operator norms. \textit{Optimization Letters} {\bf 8}(7), 2099--2110 (2014)

\bibitem{AMOU2011} Moudafi, A.: Split monotone variational inclusion.
\textit{Journal of Optimization Theory and Applications} {\bf 150}(2), 275--283 (2011)

\bibitem{NNG2019} Necoara, I., Nesterov, Yu., Glineur, F.: Linear convergence of first order methods for non-strongly convex optimization. \textit{Mathematical Programming} {\bf 175}(1--2), 69--107 (2019)

\bibitem{NEST2018} Nesterov, Y.: \textit{Lectures on Convex Optimization} (Second Edition). Springer, Cham (2018)

\bibitem{YNES2011} Nesterov, Y., Scrimali, L.: Solving strongly monotone
variational and quasi-variational inequalities. \textit{Discrete and Continuous Dynamical Systems} {\bf 31}(4), 1383--1396 (2011)

\bibitem{PATW2005} Patwari, N., Ash, J.N., Kyperountas, S., Hero, A.O., Moses, R.L., Correal, N.S.: Locating the nodes: cooperative localization in wireless sensor networks. \textit{IEEE Signal Processing Magazine} {\bf 22}(4), 54--69 (2005)

\bibitem{ARDP1992} Pierro, A.R.D.: An extended decomposition through formalization in product spaces. Technical Report 27/92, IMECC, Universidade Estadual de Campinas, 1992.

\bibitem{SREI22022} Reich, S., Tuyen, T.M.: A new approach to solving
split equality problems in Hilbert spaces. \textit{Optimization} {\bf 71}(15),
4423--4445 (2022)

\bibitem{RTRO2009} Rockafellar, R.T., Wets, R.J.-B.: \textit{Variational Analysis} (Third Edition). Springer, Berlin (2009)

\bibitem{ASHA2009} Shapiro, A., Dentcheva, D., Ruszczyński, A.: \textit{Lectures on Stochastic Programming: Modeling and Theory} (Third Edition). Society for Industrial and Applied Mathematics (2021)

\bibitem{YSHE2015} Shehu, Y., Cai, G., Iyiola, O.S.: Iterative approximation of solutions for proximal split feasibility problems. \textit{Fixed Point Theory and Applications} {\bf 2015}(1), 1--18 (2015) (Article no. 123)

\bibitem{XWAN2026} Wang, X., Hu, Y., Li, C., Yang, X.: Convergence rates of subgradient algorithms for convex feasibility problems. \textit{Applied Mathematics $\&$ Optimization} {\bf 94}, 1--32 (2026) (Article no. 76)

\end{thebibliography}
\end{document}